\documentclass[11pt]{article}

\usepackage[dvipsnames]{xcolor}
\usepackage{graphicx}
\usepackage{subcaption}
\usepackage[margin=1.0in]{geometry}
\usepackage[utf8]{inputenc}
\usepackage[T1]{fontenc}
\usepackage{pdfpages}
\usepackage[
    colorlinks=true,
    linkcolor=RoyalPurple,
    citecolor=Blue,
    urlcolor=Violet,
    linktoc=page
]{hyperref}

\usepackage{url}
\usepackage{amsfonts}
\usepackage[numbers, sort]{natbib}
\usepackage{amsmath, amssymb, amsfonts, amsthm}
\usepackage{mathrsfs}
\usepackage{booktabs}
\usepackage{aliascnt}
\usepackage[capitalise]{cleveref}

\usepackage{algorithm}
\usepackage{algpseudocode}
\usepackage{nicefrac}
\usepackage{bbm}
\usepackage{microtype}
\usepackage{tikz}
\usepackage[textsize=scriptsize]{todonotes}
\makeatletter
\renewcommand{\@todonotes@drawLineToRightMargin}{%
  \if@todonotes@line
    \begin{tikzpicture}[remember picture,overlay]
      \coordinate (mdTextAnchor) at ([yshift=-0.2cm]inText);
      \draw[draw=\@todonotes@currentlinecolor,line width=0.4pt]
        (mdTextAnchor) -- (inNote.west |- mdTextAnchor);
    \end{tikzpicture}%
  \fi
}
\renewcommand{\@todonotes@drawLineToLeftMargin}{%
  \if@todonotes@line
    \begin{tikzpicture}[remember picture,overlay]
      \coordinate (mdTextAnchor) at ([yshift=-0.2cm]inText);
      \draw[draw=\@todonotes@currentlinecolor,line width=0.4pt]
        (mdTextAnchor) -- (inNote.east |- mdTextAnchor);
    \end{tikzpicture}%
  \fi
}
\makeatother
\usepackage{enumitem}
\usepackage{mathtools, nccmath, extarrows}
\usepackage{bm}
\usepackage{scalerel}
\usepackage{comment}

\usepackage[most]{tcolorbox}
\usepackage{varwidth}
 \newtcbox{\mybox}{colback=blue!5,
	colframe=blue!30!black, center, enhanced, varwidth upper}

\usepackage{stackengine}
\mathtoolsset{showonlyrefs} %
\stackMath

\DeclareFontEncoding{LS1}{}{}
\DeclareFontSubstitution{LS1}{stix}{m}{n}

\makeatletter
  \newcommand*\textmathversion{\csname textmv@\math@version\endcsname}
  \newcommand*\textmv@normal{m}
  \newcommand*\textmv@bold{b}
\makeatother

\newtheorem{theorem}{Theorem}
\newaliascnt{lemma}{theorem}
\newtheorem{lemma}[lemma]{Lemma}
\aliascntresetthe{lemma}
\newaliascnt{corollary}{theorem}
\newtheorem{corollary}[corollary]{Corollary}
\aliascntresetthe{corollary}
\newaliascnt{proposition}{theorem}
\newtheorem{proposition}[proposition]{Proposition}
\aliascntresetthe{proposition}
\newaliascnt{definition}{theorem}
\newtheorem{definition}[definition]{Definition}
\aliascntresetthe{definition}
\newaliascnt{assumption}{theorem}
\newtheorem{assumption}[assumption]{Assumption}
\aliascntresetthe{assumption}
\newaliascnt{remark}{theorem}
\newtheorem{remark}[remark]{Remark}
\aliascntresetthe{remark}
\newaliascnt{claim}{theorem}

\aliascntresetthe{claim}
\newaliascnt{fact}{theorem}

\aliascntresetthe{fact}
\newaliascnt{example}{theorem}
\newtheorem{example}[example]{Example}
\aliascntresetthe{example}
\newaliascnt{property}{theorem}

\aliascntresetthe{property}
\newaliascnt{obs}{theorem}

\aliascntresetthe{obs}
\newaliascnt{conjecture}{theorem}

\aliascntresetthe{conjecture}

\crefname{theorem}{Theorem}{Theorems}
\Crefname{theorem}{Theorem}{Theorems}
\crefname{lemma}{Lemma}{Lemmas}
\Crefname{lemma}{Lemma}{Lemmas}
\crefname{corollary}{Corollary}{Corollaries}
\Crefname{corollary}{Corollary}{Corollaries}
\crefname{proposition}{Proposition}{Propositions}
\Crefname{proposition}{Proposition}{Propositions}
\crefname{definition}{Definition}{Definitions}
\Crefname{definition}{Definition}{Definitions}
\crefname{assumption}{Assumption}{Assumptions}
\Crefname{assumption}{Assumption}{Assumptions}
\crefname{remark}{Remark}{Remarks}
\Crefname{remark}{Remark}{Remarks}
\crefname{claim}{Claim}{Claims}
\Crefname{claim}{Claim}{Claims}
\crefname{fact}{Fact}{Facts}
\Crefname{fact}{Fact}{Facts}
\crefname{example}{Example}{Examples}
\Crefname{example}{Example}{Examples}
\crefname{property}{Property}{Properties}
\Crefname{Property}{Property}{Properties}
\crefname{obs}{Observation}{Observations}
\Crefname{obs}{Observation}{Observations}
\crefname{conjecture}{Conjecture}{Conjectures}
\Crefname{conjecture}{Conjecture}{Conjectures}

\newcommand{\mdrev}[1]{{\color{black}#1}}
\newcommand{\Ustar}{U^{\star}}
\newcommand{\UstarT}{U^{\star\top}}
\newcommand{\ustar}{u^{\star}}
\newcommand{\Spm}{S_d^{(\pm)}}
\newcommand{\Norbit}{\mathcal{N}_{\mathrm{orb}}}
\newcommand{\Ustarntilde}{\widetilde{U}_n^{\star}}
\newcommand{\dframe}{d_{\mathrm{frame}}}

\newcommand{\F}{\mathcal{F}}

\reversemarginpar

\newcommand{\uF}[1]{\underline{\mathbb{F}}}
\newcommand{\one}{\mathbf{1}}

\newcommand{\Ac}[0]{\mathcal{A}}
\newcommand{\Bc}[0]{\mathcal{B}}

\newcommand{\Ec}[0]{\mathcal{E}}

\newcommand{\Ic}[0]{\mathcal{I}}

\newcommand{\Nc}[0]{\mathcal{N}}

\newcommand{\Pc}[0]{\mathcal{P}}

\newcommand{\Sc}[0]{\mathcal{S}}
\newcommand{\Tc}[0]{\mathcal{T}}

\newcommand{\cE}{\mathcal{E}}

\newcommand{\RR}{{\mathbf R}}

\newcommand{\Span}{\operatorname{span}}

\newcommand{\Tr}{\operatorname{tr}}

\newcommand{\dist}{\operatorname{dist}}

\newcommand{\EE}{\operatorname{\mathbb E}}
\newcommand{\Var}{\operatorname{Var}}

\newcommand{\PP}{\operatorname{\mathbb P}}

\newcommand{\supp}{\operatorname{supp}}

\newcommand{\argmin}{\operatornamewithlimits{argmin}}
\newcommand{\argmax}{\operatornamewithlimits{argmax}}
\newcommand{\bigO}{\mathcal{O}}
\newcommand{\littleo}{\vcenter{\hbox{\scalebox{0.6}{$\mathcal{O}$}}}}

\newcommand{\conv}{\operatorname{conv}}

\newcommand{\diam}{\operatorname{diam}}

\DeclarePairedDelimiter{\dotp}{\langle}{\rangle}

\def\shortdisplay{\setlength{\abovedisplayskip}{5pt}%
	\setlength{\belowdisplayskip}{5pt}%
	\setlength{\abovedisplayshortskip}{2pt}%
	\setlength{\belowdisplayshortskip}{2pt}}

\let\oldselectfont\selectfont
\def\selectfont{\oldselectfont\shortdisplay}

\newcommand{\name}{\text{orthogonal K-SVD with Riemannian GD}}
\newcommand{\Csw}{C_{\mathrm{sw}}}
\newcommand{\Cmar}{C_{\mathrm{mar}}}
\newcommand{\Cgeom}{C_{\mathrm{geom}}}
\newcommand{\Sswitch}{\Sc_{\mathrm{switch}}}
\newcommand{\Cstat}{C_{\mathrm{stat}}}

\newcommand{\Enet}{\Ec_{\mathrm{net}}}
\newcommand{\Equasi}{\Ec_{\mdrev{\mathrm{coarse}}}}

\newcommand{\Eloc}{\Ec_{\mathrm{loc}}}

\newcommand{\Egood}{\Ec_{\mathrm{good}}}
\newcommand{\Efar}{\Ec_{\mathrm{far}}}
\newcommand{\Hess}{\operatorname{Hess}}

\newcommand{\Is}{\Ic_{\mathrm{s}}}
\newcommand{\Ins}{\Ic_{\mathrm{ns}}}
\newcommand{\rhomu}{\rho_{\mu}}

\newcommand{\etau}{\varepsilon_{\tau}}
\newcommand{\Srho}{\Sc_{\rho}}
\newcommand{\Srhomu}{\Sc_{\rhomu}}

\newcommand{\mfour}{M_{4}}
\newcommand{\mtwo}{M_{2}}

\newcommand{\orth}{\mathrm{O}}

\newcommand{\Rfar}{\mathcal{R}_{\mathrm{far}}}
\newcommand{\Rnear}{\mathcal{R}_{\mathrm{near}}}
\newcommand{\Capp}{\mathrm{Cap}}
\newcommand{\Exp}{\operatorname{Exp}}
\newcommand{\poly}{\mathrm{poly}}
\newcommand{\Fr}{\mathrm{Fr}}
\newcommand{\opp}{\mathrm{op}}
\newcommand{\tube}{\operatorname{Tube}}

\definecolor{mdrevision}{RGB}{120,81,169}
\long\def\mdrevised#1{{\color{mdrevision}#1}}
\expandafter\let\csname md_revised\endcsname\mdrevised

\newcommand{\Rgrad}{\nabla^{\mathrm{R}}}
\newcommand{\Egrad}{\nabla^{\mathrm{E}}}

\begin{document}

\title{Wasserstein approximation of measures via orthogonal frames: \\ statistical and computational complexities}
\author{
  Mateo D\'iaz\thanks{Department of Statistics \& Data Science Institute, University of Chicago, Chicago, IL 60637. MD and DLC were partially supported by NSF awards CCF 2442615 and DMS 2502377, and a Sloan Research Fellowship.} \and
  Nicol\'as Garc\'ia Trillos\thanks{Department of Statistics, University of Wisconsin--Madison, Madison, WI, 53706. NGT and CY were partially supported by NSF-DMS award 2236447.} \and
  Daniel L\'opez-Casta\~no\thanks{Department of Applied Mathematics and Statistics, Johns Hopkins University, Baltimore, MD 21218. DLC was partially supported by a Fulbright--Minciencias Scholarship. 
\\\indent \textbf{Authors are listed in alphabetical order by family name.}} \and
  Congwei Yang\footnotemark[2]
   }
\date{}

\maketitle

\begin{abstract}
Dimension-reduction methods such as PCA identify a subspace capturing most of the variation in the data, but not a distinguished basis within it. Selecting a basis in which the samples are sparsely represented is a recurring task in unsupervised learning, addressed by methods in factor analysis, dictionary learning, and subspace clustering. We pose this task as Wasserstein approximation of the data distribution by measures supported on the axes of an orthogonal frame, which, in turn, reduces to a
nonsmooth optimization problem over the orthogonal group $\orth(d)$. We show that, under mild assumptions, both empirical risk minimizers and a locally initialized stochastic gradient method recover the population-optimal frame at the near-parametric rate $\widetilde{\bigO}_d(n^{-1/2})$, where $n$ is the sample size. The main difficulty is that, although the population objective is smooth, the empirical objective is not---this mismatch precludes standard covering arguments. %
We overcome this issue by replacing Lipschitz continuity of the empirical gradient with a coarsely Lipschitz bound, a tool that may be of independent interest. Numerical experiments support our theoretical findings.

\end{abstract}

\section{Introduction}

Approximating a probability measure by one with prescribed structure is a fundamental problem in statistics and machine learning. Readily, this amounts to solving an optimization problem of the form $\inf_{\nu \in \F} D(\mu,\nu),$ where $\mu$ is the target probability measure, $\F$ is a class of structured probability measures, and $D$ is a discrepancy between probability measures. The choice of $\F$ specifies the structure sought in the approximation, while $D$ quantifies the resulting approximation error. This framework encompasses, among other examples, measure quantization, approximation by lattice-supported measures, minimum-distance estimation, and generative modeling \cite{GrafLuschgy2000,HammKhurana2023,Bassetti,ArjovskyChintalaBottou2017}.

The $2$-Wasserstein distance $W_2$ is particularly well suited to structured approximations supported on lower-dimensional sets and therefore not absolutely continuous with respect to the target measure. {Unlike divergences that require absolute continuity, the Wasserstein distance remains finite even when its input measures do not share the same support.} To give a concrete and widely used example, if $\mu$ is centered and $\F$ is the class of probability measures supported on arbitrary $k$-dimensional linear subspaces, minimizing $W_2^2(\mu,\nu)$ over $\nu \in \F$ recovers PCA: the optimal subspace is spanned by the top $k$ eigenvectors of the covariance matrix of $\mu$ and the optimal $\nu$ is the projection of $\mu$ onto this subspace.

\begin{figure}[ht]
    \centering
    \includegraphics[width=\linewidth]{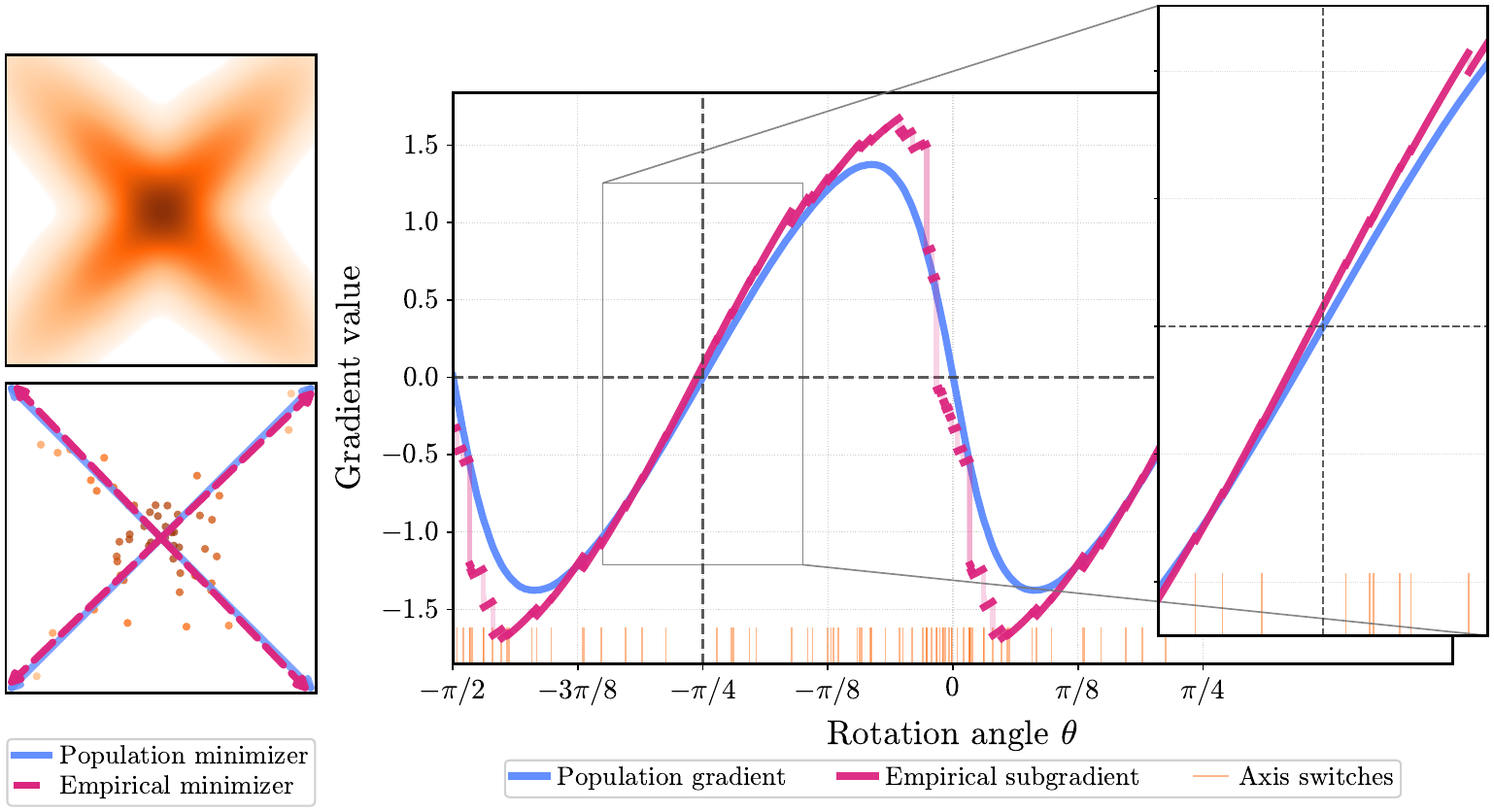}
    \caption{
\small
\textbf{Population and empirical frame recovery in $d=2$.} When $d=2$, frames are associated with a rotation angle $\theta$.
\textit{Top left:} Heat map of the population Gaussian mixture.
\textit{Bottom left:} An empirical sample from this population, together with the frames associated with the population minimizer $\theta^{\star}$ and the empirical minimizer $\theta_{n}^{\star}$.
\textit{Right:} The Riemannian gradient of the population objective $\Rgrad\Ec(U)$ (blue) and the empirical subgradient set (magenta), plotted against the rotation angle. Orange ticks mark jump discontinuities. These switches are less frequent near the population minimizer
$\theta^\star=-\pi/4$.}
        \label{fig:WFA1-sum_axis}
\end{figure}

\paragraph{The problem} We focus on approximating a target measure by probability measures supported on the axes of an orthonormal frame. {Throughout, we use $d\geq2$ to denote the dimension of the ambient space. An orthogonal frame is parameterized by an orthogonal matrix $U=[u_1,\ldots,u_d]\in\orth(d)$; we note that we can only identify a frame with $U$ up to signed permutations. Any such $U$ defines a collection of probability measures that are supported on its associated frame, namely} \( \F_U:=\left\{\nu\in\Pc_2(\RR^d)\,\middle|\, \supp(\nu)\subseteq\bigcup_{j=1}^d\Span\{u_j\}\right\}. \) Given a target measure $\mu$, our idealized goal is to solve the population approximation problem
\begin{equation}\label{eq:generic_approx}
    \min_{\nu\in\F} W^2_2(\mu,\nu) \qquad \text{with} \qquad \F = \bigcup_{U \in \orth(d)} \F_U.
\end{equation}
{A solution consists of an optimal measure $\nu^{\star}$ together with an orthogonal matrix $U^{\star}$ such that $\nu^{\star}\in\F_{U^{\star}}$. As is common in applications, we do not assume direct access to $\mu$; instead, we access a finite i.i.d.\ sample $x_{1},\dots,x_{n}\sim\mu$. Writing $\mu_n$ for the empirical measure associated with this sample, we replace $\mu$ in \eqref{eq:generic_approx} by $\mu_n$ and consider the empirical problem $\inf_{\nu\in\F}W_2^2(\mu_n,\nu)$ with minimizers $\nu^{\star}_{n}$ and $U^{\star}_{n}$. Figure~\ref{fig:WFA1-sum_axis} illustrates the population and empirical problems in dimension two. Our goals are twofold: $(i)$ on the statistical front, we wish to quantify how accurately an optimal empirical frame $U_n^\star$ estimates the population-optimal frame $U^\star$ as a function of the sample size, and $(ii)$ on the computational front, we would like to analyze the use of SGD for estimating $U^\star$ from finite samples.}

\paragraph{Motivation} Our work is motivated by a phenomenon occurring in several applications: many latent-variable and spectral methods identify a low-dimensional subspace without identifying a distinguished basis within it. This rotational ambiguity appears in factor analysis, where factor scores and loadings can be rotated simultaneously without changing the fitted model \cite{thurstone1935,thurstone1947,browne2001analytic,rohe2023}, and in dictionary learning and sparse coding, where one seeks a basis in which observations have sparse coefficient vectors \cite{kSVD,SurveySparseRepresentation,mairal2010online,spielmanwangwright2012,sunquwright2017}; as we shall soon discuss (see our section on related work below), our problem \eqref{eq:generic_approx} can actually be written as a special instance in a family of optimization problems arising in \textit{orthogonal dictionary learning}. We study, from a statistical perspective, a deliberately strong geometric form of this structure: after an unknown orthogonal transformation, each represented vector concentrates near a single coordinate axis and therefore has one dominant coordinate. Depending on the application, these vectors may represent factor scores of observations, loading vectors of features, or coefficients in an orthogonal dictionary model. Classical criteria such as Varimax \cite{kaiser1958,rohe2023} and Quartimax \cite{neuhaus1954quartimax} select a basis by rewarding coordinate concentration. Our Wasserstein formulation instead selects the frame that minimizes the transportation cost from the distribution to a measure supported on its axes. One could consider a more general model where the measures are supported on unions of orthogonal subspaces (see the discussion in the related work section below). Nonetheless, we favor this simple structure because it already exposes an interesting theoretical tension that we explain next.

\paragraph{The challenge} To state our main contributions, it is convenient to first discuss some of the structural properties of \eqref{eq:generic_approx}. At first glance, this appears to be an infinite-dimensional problem since it requires minimizing over the measure $\nu$. However, the Wasserstein distance enables a key simplification: minimization over measures can be reduced to minimization over frames. Indeed, for each fixed frame $U$, an optimal measure over $\F_U$ admits an explicit closed-form characterization. Formally, let $P_U$ be any measurable nearest-axis projection satisfying
\begin{equation}\label{eq:projection_operator}
P_U(x) \in
\left\{
    (u_j^\top x)u_j :
    j \in \argmin_{i\in[d]}
    \bigl\|x-(u_i^\top x)u_i\bigr\|^2
\right\}.
\end{equation}
Then the pushforward measure $P_U\#\mu$ minimizes $W_2^2(\mu,\nu)$ over $\nu\in\F_U$. Notice that the constraint set is $\F_U$, rather than $\F$; see Proposition~\ref{prop:nearest-axis-wasserstein} in the appendix for a formal proof. Hence, the population and empirical versions of \eqref{eq:generic_approx} reduce, respectively, to minimizing the following objectives on $\orth(d)$:
\begin{equation}\label{eq:population_empirical_objectives}
{
\Ec(U):=\EE_{x\sim\mu}\ell(x,U) \qquad \text{and}
\qquad \Ec_n(U):=\frac1n\sum_{i=1}^n\ell(x_i,U),}
\end{equation}
where $\ell(x,U):= \|x - P_{U}(x)\|^{2}=\min_{j\in[d]}\|x-(u_j^\top x)u_j\|^2.$
Thus, our formulation reduces to a finite-dimensional stochastic optimization problem over the compact Riemannian manifold $\orth(d)$, giving us access to a mature algorithmic toolbox \cite{AbsMahSep2008,boumal2023intromanifolds}.

Interestingly, these objectives violate standard assumptions used in stochastic optimization or in the analysis of  empirical risk minimizers and their generalization guarantees. Traditional proof techniques assume that the loss function $\ell(x,\cdot)$ is smooth, so both $\Ec$ and $\Ec_n$ are smooth as well. This facilitates transferring properties of the population gradient, which are typically simpler to establish, to the empirical gradient. Here, however, there is a regularity mismatch: under fairly mild assumptions, $\Ec$ is smooth, whereas both $\ell(x,\cdot)$ and $\Ec_n$ are only piecewise smooth. Indeed, as the frame varies, the axis nearest to a given sample point $x$ changes whenever two projection residuals in $\ell(x,\cdot)$ become equal, creating a ``cusp'' in the objective and a jump in the gradient. Thus, there are at most $\binom{d}{2}=\bigO(d^2)$ switching surfaces per sample. As a result, the empirical gradient has jumps and is not Lipschitz continuous. Population averaging can smooth out these discontinuities. Figure~\ref{fig:WFA1-sum_axis} illustrates these observations with a simple two-dimensional problem.

\paragraph{Statistical complexity} Our first goal is to show that the distance between an empirical minimizer $\Ustar_{n}$ and the population minimizer $\Ustar$ decreases at a sufficiently fast rate as $n$ grows. To this end, let us introduce the natural distance to measure the fit between two frames:
\begin{equation}\label{eq:dframe}
\dframe(U,V):=\min_{R\in\Spm}d_g(UR,V),
\end{equation}
where $\Spm$ is the set of signed permutation matrices and $d_g$ is the geodesic distance on $\orth(d)$. Thus, $\dframe$ aligns column signs and order before comparing the two orthogonal matrices. Let $G_{n}(U)$ denote an empirical subgradient selection of $\Ec_n$ and let $G(U):=\Rgrad\Ec(U)$ denote the population gradient.
Standard arguments to establish a bound on $\dframe(\Ustar_{n}, \Ustar)$ follow a straightforward three-step strategy: (1) establish pointwise gradient concentration, i.e., {for each fixed $U$} we have that $\|G_{n}(U)-G(U)\| = \bigO_{p}(n^{-1/2})$; (2) use Lipschitz continuity of the empirical gradient in tandem with an ${\varepsilon}$-net argument to extend this to uniform concentration, i.e., $\sup_{U}\|G_{n}(U)-G(U)\| = \bigO_{p}(n^{-1/2})$; and (3) use strong convexity of the population objective near a minimizer to translate gradient concentration into a bound on $\dframe(\Ustar_{n}, \Ustar)$. In our setting, the first and third {steps} follow from standard arguments, but the {second} one fails due to the lack of Lipschitzness of $G_{n}(U)$.

{Nonetheless, there is an important observation that allows us to bypass this difficulty: the jumps of the empirical subgradient $G_n$ near the minimizer are small, so $G_n$ obeys a local coarsely Lipschitz bound, i.e., it is Lipschitz up to an additive constant \cite{mosher2003stable,kapovichsardar2022trees}. Let us quickly summarize how we establish this property through ``switching samples.'' For two nearby frames $U_1,U_2$, let $\mathcal I_{\mathrm{sw}}(U_1,U_2)$ contain the indices of samples whose nearest axis switches. The individual subgradient bounds give}

\begin{equation}\label{eq:coarse-lip}
\begin{aligned}
\|G_n(U_1)-G_n(U_2)\|_{\Fr}
&\leq
\underbrace{\left(\frac{2}{n}
\sum_{i\notin\mathcal I_{\mathrm{sw}}(U_1,U_2)}\|x_i\|_2^2\right)\|U_1-U_2\|_{\Fr}}_{\text{Lipschitz term}}+
\underbrace{\frac{4}{n}
\sum_{i\in\mathcal I_{\mathrm{sw}}(U_1,U_2)}\|x_i\|_2^2}_{\text{switching term}}.
\end{aligned}
\end{equation}
The first term is proportional to the frame distance, while the second is controlled by the number and size of the switching samples. A switch requires the normalized sample direction to lie near a low-dimensional region of the space. One of our main assumptions ensures that the number of switching points is small and, since each contribution is weighted by $1/n$, their aggregate jump decreases with $n$. This yields a coarsely Lipschitz substitute for ordinary Lipschitz continuity and allows the pointwise concentration bound to be extended from the net to the empirical minimizer. Leveraging these ideas, we establish the following guarantee.

\begin{theorem}[Empirical minimizer convergence -- informal version of Theorem~\ref{thm:convergence}]\label{thm:informal-statistical}
    Let $X \sim \mu$ and assume that $\|X\|^{2}$ is a $K$-sub-exponential random variable and that there exists a constant $\Cmar > 0$ such that the probability of falling within width $h$ of a relevant switching region grows at most like $\Cmar \cdot h$. Further, suppose that the population objective $\cE(\cdot)$ is geodesically $m$-strongly convex and has an $L$-Lipschitz Riemannian gradient near a minimizer $U^{\star}.$  Assume also that, for some $\rho>0$, the sublevel set $\{U\in\orth(d):\Ec(U)\leq\Ec(U^\star)+\rho\}$ is contained in the signed-permutation copies of this neighborhood.
    Then, for all large $n$, the following bound holds with high probability,
\[
\dframe(\Ustar_n,\Ustar)
\lesssim \frac{L+K d+K\Cmar d^2}{m}\cdot\frac{\log n}{\sqrt n}.
\]
\end{theorem}
Let us comment on three aspects of this result. First, the $K\Cmar d^2$ term arises from the number of switching regions. Second, we state this version of the theorem to favor generality and explicitly decouple the fundamental ingredients required for a bound---we expect it to cover a wide variety of settings. We instantiate it for Gaussian mixtures concentrating near a frame and show that the constants exhibit dimension dependencies; in particular, $\Cmar=\bigO(d^{1/2})$, $K\asymp1$, and $m\asymp L\asymp d^{-1}$, which yields an upper bound of $\widetilde{\bigO}(\sqrt{d^{7}/n})$. We do not expect the dependence on $d$ to be optimal.
Third, the $n^{-1/2}$ rate is the parametric rate for estimation in finite dimensions and bypasses the exponential sample complexity required to approximate $\mu$ in $W_2$ from samples, which has a rate of roughly $n^{-1/d}$ \cite{panaretos2019,weedbach2019}.

 \paragraph{A computational guarantee} Next, we turn to the problem of numerically approximating $U^{\star}$. Our approach to this problem takes a slightly different perspective from the one used for our statistical results. Instead of analyzing the empirical minimizer, we analyze an online stochastic gradient descent (SGD) algorithm that uses a fresh independent batch of $B$ samples at each iteration. Given a sample budget $n$, we run $T$ iterations with $B=\lfloor n/T\rfloor$, and thus use at most $n$ samples. Formally, at each iteration we update
 \begin{equation}\label{eq:rie-SGD}
U_{t+1}=\Exp_{U_t}\!\left(-\eta g(U_t,z_t)\right),
\end{equation}
where $z_t=\{x_{t,1},\ldots,x_{t,B}\}$ is a fresh batch of samples, $g(U_t,z_t)$ is an estimator of the population Riemannian gradient $\Rgrad\Ec(U_t)$ akin to the subgradient $G_{n}$ we use in \eqref{eq:coarse-lip}, $\eta>0$ is the step size, and $\Exp_{U_t}$ is the Riemannian exponential map on $\orth(d)$; see Section~\ref{sec:main_results} for the precise constructions. As our second contribution, we establish a local convergence rate.

\begin{theorem}[Stochastic gradient descent -- informal version of Theorem~\ref{thm:SGD_convergence_thm_prob}]\label{thm:informal-sgd}
    Let $X\sim\mu$ and assume that $\|X\|^2$ is $K$-subexponential. Suppose that the population objective $\Ec(\cdot)$ is geodesically $m$-strongly convex with $L$-Lipschitz gradient near a population minimizer $\Ustar$. Run iteration~\eqref{eq:rie-SGD} with
    $
    \eta=\frac{2}{m+L},
    T=\left\lceil \frac{L}{m}\log n\right\rceil,$ $B=\lfloor n/T\rfloor,$ and an initial point $U_0$ that is sufficiently close to $\Ustar$.
    Then, for all large $n$, the following bound holds with high probability,
    \[
    \dframe(U_T,\Ustar)
    \lesssim \left({d_g(U_0,\Ustar)}
    +\frac{K\sqrt L}{m^{3/2}}\right)\frac{\log n}{\sqrt n}.
    \]
\end{theorem}
A few observations are in order. First, unlike Theorem~\ref{thm:informal-statistical}, the proof does not rely on coarse Lipschitz continuity or uniform control of switching boundaries. Instead, the use of fresh batches reduces the analysis to contraction of the smooth population objective together with pointwise concentration of the stochastic-gradient oracle. Second, general local and high-probability guarantees are available for SGD and Riemannian stochastic optimization~\cite{FehrmanGessJentzen2020,pmlr-v235-dodd24a}, but they do not directly cover our setting. We discuss this in more detail in Section~\ref{sec:RelatedWork}. Third, the theorem attains the parametric \(n^{-1/2}\) dependence, up to logarithmic factors, matching the rate of the empirical minimizer. Although the abstract bound has no explicit dependence on \(d\), its geometric and distributional constants do. For a balanced Gaussian-mixture model, \(K\asymp 1\) and \(m\asymp L\asymp d^{-1}\), so the bound specializes to $\widetilde{\bigO}\left({d}/{\sqrt{n}}\right).$ Fourth, this guarantee is local. We expect that a two-stage procedure combining a spectral initializer with SGD could yield a global recovery guarantee, potentially at the cost of additional samples and stronger model-specific assumptions, as in analogous analyses of related problems \cite{pmlr-v162-wang22r, davis2025clustering}. This is beyond the scope of this paper. 

\paragraph{Outline} The rest of the paper is organized as follows. The remainder of this section covers related literature. %
Section~\ref{sec:background} introduces notation and necessary technical background. Section~\ref{sec:prob_setting} outlines the assumptions for the analysis and introduces a statistical model that we prove satisfies these requirements. Section~\ref{sec:main_results} introduces our main theoretical results. Section~\ref{sec:analysis} presents the core ideas underlying our technical arguments. Section~\ref{sec:experiments} closes the paper with numerical examples illustrating the performance of the method discussed in the paper in high-dimensional settings. Long technical proofs are deferred to the appendix.

\subsection{Related work}
\label{sec:RelatedWork}
\paragraph{Wasserstein approximation}
 There is a vast body of work dealing with $W_2$-approximation of a measure $\mu$ over a structured family.
Optimal quantization \cite{GrafLuschgy2000} considers measures supported on at most $N$ points, yielding sharp $\bigO(N^{-1/d})$ rates under regularity conditions, though this suffers from the curse of dimensionality as accuracy degrades exponentially with dimension. Extensions of this framework include Wasserstein barycenters \cite{AguehCarlier2011}, Wasserstein $k$-means \cite{zhuang2022wasserstein}, lattice approximations \cite{HammKhurana2023}, and {projection onto Gaussian and related location-scale families} \cite{WangWei2026}. In statistics, the problem \eqref{eq:generic_approx} can also be used for inference. Indeed, $\F$ may be taken to be a parametric family of distributions and $\mu$ an empirical measure of samples from an unknown element in the family \cite{Bassetti,Belili,Bernton}. In that context, the problem \eqref{eq:generic_approx} is used to define estimators for the unknown distribution that, under some regularity assumptions, can be shown to be asymptotically sensitivity-efficient \cite{2511.07414}.
Computational considerations have also motivated tractable alternatives to full $W_2$ transport; we mention three representative examples. Wasserstein GANs use $W_1$ to promote gradient stability \cite{ArjovskyChintalaBottou2017,Gulrajani2017}. For Gaussian mixtures, restricting the class of couplings yields mixture-like transport distances \cite{ChenGeorgiouTannenbaum2018,DelonDesolneux2020}. Sliced Wasserstein distances take a different route by averaging one-dimensional transport costs over projection directions \cite{Bonnotte2013,Kolouri2018,ManoleBalakrishnanWasserman2022,ParkSlepcev2025}.

\paragraph{{Latent coordinate representations}}
{Several methods seek a distinguished coordinate system for a latent or learned representation. In factor analysis, rotational ambiguity motivates post-PCA criteria such as Varimax and Quartimax, which favor nearly sparse loading patterns \cite{kaiser1958,neuhaus1954quartimax}; recent work gives conditions under which Varimax recovers latent structure \cite{rohe2023}, albeit at a slow logarithmic rate. Dictionary learning pursues a related goal by learning atoms for which observations admit sparse coefficient vectors \cite{kSVD,SurveySparseRepresentation,olshausenfield1997,mairal2010online,spielmanwangwright2012,sunquwright2017}.  Independent Component Analysis instead selects coordinates through statistical independence \cite{HyvarinenOja2000}. Although these approaches differ in their models and identification principles, they share the aim of selecting meaningful axes that are not determined by subspace information alone. The method explored here provides a complementary geometric criterion: it selects an orthogonal frame by minimizing the Wasserstein transportation cost from the target $\mu$ to its projection onto the frame's axes.}

\paragraph{Dictionary learning and projective clustering}
The reconstruction loss in \eqref{eq:population_empirical_objectives} arises in orthogonal dictionary learning, since our loss function $\ell(x, U)$, which appears in \eqref{eq:population_empirical_objectives}, can be written as $\ell(x,U)=\min_{\|z\|_0\le1}\|x-Uz\|^2$ for $U\in\orth(d)$. Indeed, written in this form our loss function can be seen as an instance of the family of loss functions discussed in \cite{kSVD,rusu2016fast,Sezer2,Sezer1}, which consider losses of the form $\min_{\|z\|_0\le s}\|x-Uz\|^2$ for some sparsity level parameter $s$; this more general objective can be naturally associated to analogous Wasserstein approximation problems, but, as was mentioned earlier, we will restrict ourselves to the case $s=1$ for simplicity.
Algorithms for orthogonal sparse reconstruction include Householder-based transforms \cite{rusu2016fast} and geodesic-flow updates \cite{schuetze2015learning}. For a prescribed sparsity level $s$ as above, \cite{ravishankar2020analysis} establishes local linear convergence of alternating minimization to the generating dictionary under a noiseless sparse model and suitable initialization; we highlight that in this paper we analyze an online Riemannian SGD scheme under relatively general assumptions (in particular, without assuming the noiseless setting). Local identifiability of the K-SVD reconstruction criterion proposed in \cite{kSVD} (which does not force atoms in the dictionary to be orthogonal) is studied in \cite{MauerPontil,schnass2014identifiability}, while other generalization bounds for dictionary learning are derived in \cite{gribonval2015sample,pmlr-v80-georgogiannis18a,JMLR:v12:vainsencher11a}. The work \cite{bucarelli2023generalization}  controls the reconstruction risk in projective clustering, another problem whose family of loss functions includes our loss when subspaces are assumed to be $d$ mutually orthogonal lines through the origin. The generalization bounds in all previously mentioned papers control excess risk, i.e., the difference between the population \textit{objective function} at the population and empirical minimizers. Recovery guarantees have also been established for different sparsity criteria, including entrywise $\ell^4$-norm maximization \cite{zhai2020complete} and $\ell_0$-\textit{penalized} reconstruction \cite{liang2025simple}; the latter work also analyzes online learning settings. Broadly speaking, those papers prove convergence of empirical risk \textit{minimizers} toward corresponding population minimizers at a rate of roughly $\bigO_d(n^{-1/4})$, a rate that can be obtained by combining excess risk bounds with forms of local strong convexity of the population objective. In contrast, we obtain a $\widetilde{\bigO}_d(n^{-1/2})$ convergence rate, which is in closer agreement with known lower bounds for the estimation of  dictionaries from data (not necessarily orthogonal atoms); see \cite{MinimaxDictionary}. In summary, our theoretical analysis complements in new and strong ways the existing statistical guarantees in the literatures of dictionary learning and sparse coding. We refer the reader to \cite{SurveySparseRepresentation} for a more complete survey of methodologies in dictionary learning and some of the optimization schemes used to implement them.

\paragraph{Subspace methods} Subspace clustering seeks to partition data drawn from a union of low-dimensional subspaces; see \cite{Vidal2011} for an overview. Self-expressive methods, notably Sparse Subspace Clustering (SSC) and its robust variants \cite{ElhamifarVidal2013,SoltanolkotabiElhamifar2014,WangXu2016}, provide a natural comparison for our problem. These methods accommodate nonorthogonal subspaces of dimension greater than one, whereas we seek a collection of mutually orthogonal one-dimensional axes. For the noiseless problem, optimization-based alternatives include $K$-subspaces, for which local convergence and recovery guarantees are known \cite{pmlr-v162-wang22r}. In the setting we study, robust subspace clustering methods offer complementary tradeoffs. Existing pipelines \cite{ElhamifarVidal2013,SoltanolkotabiElhamifar2014} construct an $n\times n$ coefficient or affinity matrix and apply spectral clustering. Greedy and stochastic variants substantially mitigate the resulting computational burden \cite{YouRobinsonVidal2016,ChenLiYou2020}, but still rely on sample-to-sample representations that scale with $n$. By contrast, the frame formulation optimizes over a $d\times d$ orthogonal matrix and can process fresh batches without storing pairwise affinities between all samples, reducing the memory footprint in the common scenario in which $n\gg d$.

\paragraph{Population smoothness with empirical nonsmoothness}
The idea of exploiting population-level smoothness when the empirical objective lacks it is not new. The paper \cite{jin2018local} studies local minima of a possibly nonsmooth empirical risk through its uniform proximity to a smooth population risk. Related statistical work shows that nonsmooth empirical criteria can nevertheless exhibit classical behavior: earlier studies established asymptotic normality for empirical \(k\)-means centers and least-absolute-deviation estimators \cite{pollard1982central,pollard1991asymptotics}, while subsequent work derived nonasymptotic excess-distortion bounds for vector quantization under margin conditions \cite{levrard2015nonasymptotic}. A complementary line of work considers settings in which even the population objective may be nonsmooth, through general graphical or uniform convergence results for empirical subdifferentials and problem-specific analyses of nonsmooth statistical landscapes \cite{davisdrusvyatskiy2022graphical,ruan2025uniform,davisdrusvyatskiypaquette2020nonsmooth,diaz2019nonsmooth}. Finally, the bounded subgradient variation framework of \cite{diakonikolas2026optimization} is related in spirit to the SGD scheme analyzed here: both replace Lipschitz continuity of the gradient with control over how subgradients vary locally.
\paragraph{Stochastic Riemannian optimization}
{SGD theory on manifolds is fairly mature, spanning almost-sure convergence guarantees \cite{Bonnabel2013}, global complexity bounds for smooth and nonsmooth geodesically convex problems on Hadamard manifolds \cite{pmlr-v49-zhang16b}, local rates for averaged iterates near geodesically strongly convex minima \cite{pmlr-v75-tripuraneni18a}, and high-probability guarantees \cite{pmlr-v235-dodd24a}. Nonsmooth stochastic subgradient methods have also been studied on the Stiefel manifold \cite{LiChenDengQuZhuSo2021}. A result in \cite{FehrmanGessJentzen2020}, closely related to Theorem~\ref{thm:informal-sgd}, provides a general local guarantee that does not require empirical smoothness, albeit in the Euclidean setting. These results do not directly imply our computational guarantee. The global guarantees require Hadamard manifolds, whereas \(\orth(d)\) is non-negatively curved. The results in \cite{LiChenDengQuZhuSo2021} do not apply to our setting since the sample objective is not weakly convex. Moreover, under our subexponential assumption (Assumption \ref{ass:subexp-norms}) the switching branch-gradient oracle need not be almost surely bounded or samplewise Lipschitz as required by~\cite{pmlr-v75-tripuraneni18a, pmlr-v235-dodd24a}.}

\section{Preliminaries}\label{sec:background}
We collect the notation and background used in the sequel. We refer the interested reader to the excellent manuscripts \cite{Villani2009, boumal2023intromanifolds, clarke1990optimization} for more details.
\paragraph{Notation} For $U=[u_1,\ldots,u_d]\in\orth(d)$, write $\operatorname{Axes}(U):=\{\operatorname{span}(u_j):1\leq j\leq d\}$ for its collection of unoriented axes. For $x,y\in\RR^d$, $\dist(x,y):=\|x-y\|_2$ is the Euclidean distance, extended to a set $A\subseteq\RR^d$ by $\dist(x,A):=\inf_{y\in A}\dist(x,y)$; we distinguish this from the geodesic distance $d_g$ on $\orth(d)$ introduced below.

\paragraph{Wasserstein distance} Let $\mathcal{P}_2(\RR^d)$ denote the space of  probability measures with finite second moment. Given a Borel map $T\colon\RR^d\to\RR^d$ and a probability measure $\mu$, the pushforward under $T$ is defined by $(T\#\mu)(B):=\mu(T^{-1}(B))$. For $\mu,\nu\in\mathcal{P}_2(\RR^d)$, the squared 2-Wasserstein distance is then defined by 
\begin{equation}\label{W_2-dist}
    W_2^2(\mu, \nu) = \inf_{\gamma\in\Pi(\mu,\nu)} \int_{\RR^d\times\RR^d} \|x-y\|^2\,d\gamma(x,y),
\end{equation}
where $\Pi(\mu,\nu)$ denotes the set of couplings of $\mu$ and $\nu$. A minimizer of this problem is called an \emph{optimal transport plan}.

\paragraph{Orthogonal group as a Riemannian manifold} The \textit{orthogonal group} $\orth(d) = \{U \in \RR^{d \times d} : U^\top U = I\}$ is an embedded submanifold of $\RR^{d \times d}$ with \textit{tangent space} at each point $U \in \orth(d)$ given by $T_U \orth(d) = \{Z \in \RR^{d \times d} \colon U^\top Z + Z^\top U = 0\}$. We equip $\orth(d)$ with the Riemannian metric induced by the Frobenius inner product, whose associated norm is $\|Z\|_F = \sqrt{\Tr(Z^\top Z)}$. The corresponding \textit{geodesic distance} $d_g(U_1, U_2)$ is the length of a shortest geodesic connecting $U_1,U_2\in\orth(d)$, while the \textit{exponential map} $\Exp_U \colon T_U \orth(d) \to \orth(d)$ follows the geodesic starting at $U$ with velocity $Z$ for unit time, explicitly given by $\Exp_U(Z) = U\exp(U^\top Z)$, where $\exp$ denotes the matrix exponential. Whenever the minimizing geodesic from $U$ to $V$ is unique and $\Exp_U^{-1}(V)$ is well-defined, we write $\operatorname{Log}_U(V):=\Exp_U^{-1}(V)$ for the Riemannian logarithm and $\mathsf P_{V\to U}\colon T_V\orth(d)\to T_U\orth(d)$ for parallel transport along that geodesic from $V$ to $U$. Finally, for a smooth function $f\colon\orth(d)\to\RR$, the projection of an ambient matrix $Z$ onto $T_U\orth(d)$ is given by $\Pi_U(Z):=\tfrac12(Z-UZ^\top U)$, yielding the \textit{Riemannian gradient} $\Rgrad f(U):=\Pi_U(\Egrad f(U))$, where $\Egrad f(U)$ is the Euclidean gradient. For further details, see \cite{AbsMahSep2008,boumal2023intromanifolds}.

\paragraph{Nonsmooth analysis}
For a locally Lipschitz function $f\colon\orth(d)\to\RR$, its \textit{Riemannian Clarke subdifferential} at $U\in\orth(d)$ is defined by
\[
\partial_Rf(U):=\conv\!\left\{\lim_{k\to\infty}\Rgrad f(U_k):U_k\to U,\ f\text{ differentiable at }U_k\right\},
\]
where we use the ambient matrix identification of tangent spaces: if $Z_k\in T_{U_k}\orth(d)$, $U_k\to U$, and $Z_k\to Z$, then $U_k^\top Z_k+Z_k^\top U_k=0$ implies $Z\in T_U\orth(d)$. Appendix~\ref{app:subdiff} establishes the active-branch formula for the empirical loss subdifferential used throughout the paper. In particular, every local minimizer satisfies $0\in\partial_Rf(U)$. Moreover, if $f(UR)=f(U)$, then $\partial_Rf(UR)=\partial_Rf(U)R$, so the first-order optimality condition is independent of the chosen signed-permutation representative. This invariance holds for all risks considered in the paper.

\section{Assumptions}\label{sec:prob_setting} In this section, we state our main assumptions on the distribution $\mu$ and the population objective $\Ec$ in \eqref{eq:population_empirical_objectives}. We also show that these assumptions are satisfied by mixtures of Gaussians concentrated near frames.
We start with two  mild  assumptions on $\mu$.
\begin{assumption}[\textbf{Bounded fourth moments}]
\label{ass:concentration} The measure of interest $\mu$ has bounded fourth moments, i.e., $\mu\in\Pc_4(\RR^d)$. We write $\mtwo:=\EE_{x\sim\mu}[\|x\|_2^2]$ and $\mfour:=\EE_{x\sim\mu}[\|x\|_2^4]<\infty$.
\end{assumption}

\begin{assumption}[\textbf{Sub-exponential squared norms}]
\label{ass:subexp-norms}
For $x\sim\mu$, we have
\[
K :=\inf\left\{s>0:\EE\!\left[\exp\!\left(\frac{\|x\|^2}{s}\right)\right]\leq2\right\}<\infty.
\]
\end{assumption}

Sub-exponentiality implies a finite fourth moment \cite{vershynin2018high}. Turning to the next two assumptions, we emphasize a subtle distinction: although our gradient-based algorithm operates in the space of orthogonal matrices $\orth(d)$, we are truly interested in the space of frames---i.e., orthogonal matrices modulo signed permutations $\orth(d)/\Spm.$ This distinction is rather mild since both population and empirical objectives are invariant under signed column permutations: $\Ec(UR)=\Ec(U)$ and $\Ec_n(UR)=\Ec_n(U)$ for every $R\in\Spm$. Thus, most conclusions on $\orth(d)$ descend to the space of frames cleanly. Nonetheless, it requires us to be careful about the way we state the following two assumptions. For the first, let $\Ustar \in \argmin_{U\in \orth(d)} \Ec(U)$ be a minimizer of the population objective.

\begin{assumption}[\textbf{Local regularity of the objective}]
\label{ass:strong-convexity} There exists an open neighborhood $\Nc\subseteq\orth(d)$ of $\Ustar$ such that $\dframe(U,\Ustar)=d_g(U,\Ustar)$ for every $U\in\Nc$. For every $U,V\in\Nc$, there is a unique minimizing geodesic from $U$ to $V$, its image lies in $\Nc$, and
$\operatorname{Log}_U(V)$ is well-defined and depends smoothly on $(U,V)$. The objective $\Ec$ is $C^2$ and geodesically $m$-strongly convex on $\Nc$, i.e.,

\[
\Ec(V)\geq\Ec(U)+\bigl\langle\Rgrad\Ec(U),\operatorname{Log}_U(V)\bigr\rangle_U+\frac{m}{2}\cdot d_g(U,V)^2\qquad \text{for all } U,V\in\Nc.
\]
Further, its Riemannian gradient is intrinsically $L$-Lipschitz on $\Nc$, i.e.,
\[
\bigl\|\Rgrad\Ec(U)-\mathsf P_{V\to U}\Rgrad\Ec(V)\bigr\|_{\Fr}
\leq L\cdot d_g(U,V)\qquad \text{for all } U,V\in\Nc.
\]
\end{assumption}
This assumption is standard in the literature and allows us to bound the objective locally between two quadratic functions; it will be key to control the stability of the empirical risk minimizer and the online-SGD output. The next assumption considers the orbit of $\Nc$, namely,
\[
\Norbit:=\bigcup_{R\in\Spm}\Nc R.
\]

\begin{assumption}[\textbf{Sublevel set containment}]
\label{ass:sublevel-containment}
There exists $\rho > 0$ such that the sublevel set
\(
\Srho:=\{U\in\orth(d)\colon \Ec(U)\leq \Ec(\Ustar)+\rho\}
\) is contained in $\Norbit.$
\end{assumption}
This assumption is key for global identifiability of $\Ustar$ up to sign permutations. Next, we impose a quantitative assumption ensuring that the number of samples lying near a switching region is small for any frame near $\Ustar$. To do so, we need to introduce some additional notation.
Define the outer geodesic radius of $\Nc$ about $\Ustar$ by
\[
r_{\Nc}:=\sup_{U\in\Nc}d_g(U,\Ustar).
\]
Moreover, define the set $\Capp(u)$ around $u\in S^{{d-1}}$ to be the pair of spherical caps given by $\Capp(u) := \{ v \in S^{d-1}: |\dotp{u,v}| \geq \cos(2r_{\Nc})\}.$ For any line $\ell = \Span(u)$, define $\Capp(\ell):=\Capp(u).$
Equipped with these notions, we define the sets of nearby and far away directions as
\[
  \Rnear(\Nc):=\bigcup_{U\in\Nc}\bigcup_{\ell\in\operatorname{Axes}(U)}\Capp(\ell),
\qquad\text{and}\qquad \Rfar(\Nc):=S^{d-1}\setminus\Rnear(\Nc).
\]
These sets are invariant under signed permutations. For a fixed frame $U$, the nearest-axis assignment can change only by crossing a pairwise tie set \(\left\{v\in S^{d-1}:|\dotp{u_i,v}|=|\dotp{u_j,v}|\right\}\). For $i<j$, this set is the union of two great subspheres with unit normals $(u_i\pm u_j)/\sqrt{2}$.
The following set collects all such normals generated by $\Nc$ together with their limits,
\[
\mathcal W_{\Nc}:=
\overline{\left\{\frac{u_i\pm u_j}{\sqrt2}:U=[u_1,\ldots,u_d]\in\Nc\right\}},
\]
where the bar denotes closure.

\begin{assumption}[\textbf{Far-region switching margin}]
\label{ass:switching-margin}
The following inequality holds \[
c_{\Nc}:=\cos(2r_{\Nc})-2r_{\Nc}>\frac{1}{\sqrt2}.
\]
There also exist $\eta_0\in(0,1/2]$ and $\Cmar>0$ such that, for every
$w\in S^{d-1}$ and $0<h\leq\eta_0$ with
\(
\dist(w,\mathcal W_{\Nc})\leq\eta_0,
\)
we have
\[
\PP_{x\sim \mu}\!\left(x\neq0,\ \frac{x}{\|x\|}\in\Rfar(\Nc), \text{ and }\
\left|\left\langle\frac{x}{\|x\|},w\right\rangle\right|\leq h\right)
\leq \Cmar h.
\]
\end{assumption}
Intuitively, this assumption precludes $\mu$ from giving too much weight to regions in directions that are at roughly 45 degrees from the axes of $U^\star$. Notice that the at-most-linear growth in $h$ is not essential for consistency; see Appendix~\ref{app:erm-proofs} (Remark~\ref{rem:margin-exponent}). We close this section by exhibiting an example of a family of distributions satisfying all of the above assumptions. %
\begin{example}[Gaussian Convolution Model]\label{ex:gaussian-convolution}
    Fix $\Ustar\in\orth(d)$. Let $J$ satisfy $\PP(J=j)=p_j>0$, where $\sum_{j=1}^d p_j=1$, and, conditional on $J=j$, set $Y=S_j\Ustar_j$ with $S_j\sim\lambda_j$ on $\RR$. Equivalently, with $\iota_j\colon\RR\to\RR^d$ defined by $\iota_j(s)=s\Ustar_j$, the signal law is
    \[
    \nu=\sum_{j=1}^d p_j(\iota_j)_\#\lambda_j.
    \]
    Assume that every $\lambda_j$ has a Lebesgue density bounded by a common constant $M<\infty$, set $a_j:=p_j\EE[S_j^2]$, and assume further that
    \[
    a_{\min}:=\min_{1\leq j\leq d}a_j>0,
    \qquad a_{\max}:=\max_{1\leq j\leq d}a_j,
    \qquad \text{and} \qquad
    \max_{1\leq j\leq d}\|S_j^2\|_{\psi_1}<\infty.
    \]
    Let $Z\sim\Nc(0,I_d)$ be independent of $Y$ and set $X=Y+\tau Z\sim\mu$ for $\tau>0$. Under a signed permutation of the representative, the same law $\nu$ is described by relabeling the pairs $(p_j,\lambda_j)$ and pushing $\lambda_j$ forward by $s\mapsto-s$ when the corresponding column sign changes. The displayed bounds are invariant under this transformation.
\end{example}

To instantiate the example, one can take $\lambda_j=\Nc(0,\sigma_j^2)$ for some $\sigma_j^2>0$ and $p_j=1/d$, so the signal law is a mixture of Gaussians along the axes of $\Ustar$. The next proposition shows that, in this case, our assumptions hold uniformly on the dimension.

\begin{proposition}[Gaussian-mixtures satisfy our assumptions]
\label{prop:balanced-gaussian-mixture-uniform}
Fix variance bounds $0<\sigma_-^2\leq\sigma_+^2<\infty$ and a noise budget constant $b_+<\infty$. There exist $d_0\geq4$ and $r\in(0,1/8]$, depending only on $\sigma_-^2,\sigma_+^2$, and $b_+$, such that, for every $d\geq d_0$, the following holds. In Example~\ref{ex:gaussian-convolution}, take $p_j=1/d$, $S_j\sim\Nc(0,\sigma_j^2)$ with $\sigma_-^2\leq\sigma_j^2\leq\sigma_+^2$, and $0<\tau^2\leq b_+/d$. Then $\arg\min_{U\in\orth(d)}\Ec(U)=\Ustar\Spm$, and Assumptions~\ref{ass:concentration}--\ref{ass:switching-margin} hold on $\Nc=B_r(\Ustar)$ with
\[
m=\frac{\sigma_-^2}{2d},
\qquad
L=\frac{4\sigma_+^2}{d},
\qquad\text{and}\qquad
\rho=\frac{\sigma_-^2}{8d}\min\Bigl\{\frac12,\frac{2r^2}{\pi^2}\Bigr\}.
\]
The remaining constants satisfy
\[
\frac{\sigma_-^2}{\log2}\leq K\leq\frac{16}{3}\sigma_+^2+4b_+,
\qquad\text{and}\qquad
\Cmar\leq \Cgeom(r,\sigma_-)\,\tau d,
\]
where \(\Cgeom(r,\sigma_-)\) is an explicit dimension-independent constant determined by the geometry of the neighborhood \(\Nc=B_r(\Ustar)\) and the uniform bound \(\sigma_-\) on the component densities.
\end{proposition}

In particular, under the scaling $\tau\asymp d^{-1/2}$, the proposition yields $m\asymp L\asymp\rho\asymp d^{-1}$, $K\asymp1$, and $\Cmar=\bigO(d^{1/2})$. A natural question is whether the general model of Example~\ref{ex:gaussian-convolution} also satisfies our assumptions for all sufficiently small noise levels $\tau>0$. We answer this affirmatively in Appendix~\ref{appendix:gaussian-proofs}, albeit for fixed $d$.

\section{Main results}\label{sec:main_results} In this section, we state our two main results, which include high-probability guarantees for convergence of the empirical minimizer (Theorem~\ref{thm:convergence}, formal version of Theorem~\ref{thm:informal-statistical}) and the online Riemannian SGD algorithm (Theorem~\ref{thm:SGD_convergence_thm_prob}, formal version of Theorem~\ref{thm:informal-sgd}). All errors are measured after optimally aligning signs and column order.

\paragraph{Empirical risk minimization}
We start by making Theorem~\ref{thm:informal-statistical} precise, with explicit constants and a finite-sample threshold on $n$. The proof, including the origin of the constants below, is given in Subsection~\ref{subsec:empirical-nonsmooth}, with supporting lemmas proved in Appendix~\ref{app:erm-proofs}. For each $k\geq1$, let $\omega_k$ denote the surface area of the unit sphere $S^{k-1}\subset\RR^k$.

\begin{theorem}[Convergence rate of empirical minimizer]\label{thm:convergence}
Suppose Assumptions~\ref{ass:concentration}, \ref{ass:subexp-norms}, \ref{ass:strong-convexity}, \ref{ass:sublevel-containment}, and~\ref{ass:switching-margin} hold, and fix $\beta\in(0,1)$. Let $\rho>0$ be such that $\Srho\subseteq\Norbit$ per Assumption~\ref{ass:sublevel-containment}, and set $\rho_0:=\min\{\rho,K\}$ and $\Delta_{\Nc}:=c_{\Nc}-\sqrt{1-c_{\Nc}^2}>0$. There exist absolute constants $C_B,\Cstat>1$ such that, for every sample size $n$ satisfying
\[
n>\max\!\left\{16\eta_0^{-2},\,4\Delta_{\Nc}^{-2},\,\frac{K^2}{C_B\mtwo^2}\log(8n),\,\Cstat\frac{K^2}{\rho_0^2}\left[d^2\log\!\left(\frac{\Cstat\sqrt d\,K}{\rho_0}\right)+\log\!\left(\frac{\Cstat}{\beta}\right)\right]\right\},
\]
with probability at least
\(
1 - \beta - \frac{3}{n} - 2(d-1)\pi^{d-2}\omega_d\omega_{d-1}^{-1}n^{-(d-1)},
\)
we have
\[
\dframe(\Ustar_n,\Ustar) \;\leq\; \frac{1}{m}\left[ A_1 \cdot \frac{\log n}{\sqrt{n}} \;+\; A_2 \cdot \frac{(\log n)^{3/2}}{n^{3/4}} \;+\; A_3 \cdot \frac{(\log n)^2}{n} \right].
\]
For absolute constants $C_0,C_1>0$ independent of $d$, set \(\Gamma:=C_1(2K+2C_0\sqrt{\mfour})\) and \(\Lambda_d:=2+\frac{d(d-1)}{4}+\frac{d(d-1)\log(1 + 4r_{\Nc})}{2}\); the constants above are then
\begin{align*}
A_1 &:= 2\left (\sqrt{8\mfour}+\sqrt2\,\Gamma\right )\sqrt{\Lambda_d}+L(1+r_{\Nc})+4\mtwo
+48K\Cmar d(d-1), \\
A_2 &:= 72Kd(d-1)
\sqrt{\Cmar d},\quad \text{and} \\
A_3 &:= 4\Gamma\Lambda_d+108Kd^2(d-1).
\end{align*}
\end{theorem}

The explicit three-term bound requires no relation between $d$ and $n$ beyond the displayed thresholds. If $n\gtrsim d(\log n)^2$ and all problem constants, including $m^{-1}$ and $\Cmar$, are dimension-uniform, then its three terms are jointly $\bigO(d^{5/2}\log n/\sqrt n)$; for fixed $d$, the leading switching contribution is $\bigO(\Cmar d^2\log n/\sqrt n)$. Otherwise $A_1,A_2,A_3$ give the dimension dependence.

\paragraph{Online-stochastic Riemannian gradient descent for estimation}
Our second result makes Theorem~\ref{thm:informal-sgd} precise. To this end, let $U_0$ be deterministic, fix a batch size $B$ and a horizon $T$, and let each $z_t$ contain $B$ fresh i.i.d. samples, independently across iterations. Given a sample budget $n$, we take $B=\lfloor n/T\rfloor$, ensuring that the algorithm uses $TB\leq n$ samples. Riemannian SGD on $\orth(d)$ again follows the update~\eqref{eq:rie-SGD},
\begin{equation}\label{SGD_update}
    U_{t+1}=\Exp_{U_t}\!\left(-\eta g(U_t,z_t)\right),
\end{equation}
where the gradient estimator $g$ is constructed as follows. For $A=[a_1,\ldots,a_d]$, define $\bar\Ec(A):=\int(\|x\|^2-\max_j\langle a_j,x\rangle^2)\,d\mu(x)$. Given $U=[u_1,\ldots,u_d]$, let
\begin{equation}\label{eq:active-set}
    J(x;U):=\operatorname*{arg\,max}_{1\leq j\leq d}|\langle u_j,x\rangle|^2,\qquad\text{and}\qquad h(x,U):=-\frac{2}{|J(x;U)|}\sum_{j\in J(x;U)}\langle u_j,x\rangle x e_j^\top,
\end{equation}
so that, under Assumptions~\ref{ass:concentration} and~\ref{ass:strong-convexity}, $\EE_{x\sim\mu}[h(x,U)]=\nabla\bar\Ec(U)$ for every $U\in\Nc$ (Lemma~\ref{lem:diff-implies-tie-null}). Accordingly, averaging $h$ over a batch $z=\{x_1,\ldots,x_B\}$ and projecting onto $T_U\orth(d)$ yields an estimator of $\Rgrad\Ec(U)$. More precisely, setting $\bar g(U,z):=B^{-1}\sum_{i=1}^B h(x_i,U)$, we obtain 
\begin{equation}\label{grad_estimator}
    g(U,z)=\tfrac12\bigl(\bar g(U,z)-U\bar g(U,z)^{\top}U\bigr).
\end{equation}
Its population mean $g_\infty(U):=\EE_z[g(U,z)]$ is well-defined for every $U\in\orth(d)$ under Assumption~\ref{ass:concentration} alone and satisfies $g_\infty(U)=\Rgrad\Ec(U)$ for $U\in\Nc$. Appendix~\ref{app:SGD} (Lemma~\ref{lem:sgd-unbiased}) establishes this identity and bounds the second moment of the centered estimator $g(U,z)-g_\infty(U)$. Moreover, the identity $h(x,UR)=h(x,U)R$ implies the equivariance relations $g(UR,z)=g(U,z)R$ and $g_\infty(UR)=g_\infty(U)R$. Consequently, changing the signed-permutation representative of $U_0$ merely changes every iterate by the same right factor; see Appendix~\ref{app:SGD}.

With this construction, the following result establishes a last-iterate guarantee with prescribed high probability. In particular, the good-noise event ensures that the trajectory remains in the geodesically strongly convex ball while controlling its terminal error, thereby avoiding an additional escape-probability term. The proof is given in Section~\ref{subsubsec:sgd-high-prob-proof}, with supporting lemmas established in Appendix~\ref{app:SGD}.

\begin{theorem}[High-probability online-SGD convergence]\label{thm:SGD_convergence_thm_prob}
    Suppose Assumptions~\ref{ass:concentration} and \ref{ass:strong-convexity} hold. Choose $r_c>0$ such that $\Bc(\Ustar,r_c)\subseteq\Nc$, and let $r_{\mathrm{init}}:=d_g(U_0,\Ustar)<r_c$.
    Define $\kappa:=\tfrac{L}{m}$, and $\eta:=\tfrac{2}{m+L}$. For a fixed $c_{\mathrm{iter}}\geq1/2$, take
    \(
    T:=\left\lceil c_{\mathrm{iter}}\kappa\log n\right\rceil
    \) and \(
    B:=\left\lfloor\frac{n}{T}\right\rfloor,
    \)
    and assume the sample budget satisfies $n\geq\max\{2T, e^2\}$. Fix $\delta\in(0,1)$. The following statements hold.
    \begin{enumerate}[label=\arabic*., leftmargin=2.2em, labelsep=0.6em]
        \item (\textbf{Finite fourth moment}) Suppose further that $ 2\sqrt{\frac{\mfour T}{B\delta}}\leq mr_c$. Then, with probability at least $1-\delta$, all iterates $U_0,\ldots,U_T$ remain in $\Bc(\Ustar,r_c)$ and
        \begin{equation}\label{eq:thm2_poly}
        \dframe(U_T,\Ustar)
        \leq\frac{r_{\mathrm{init}}}{\sqrt n} +\frac{2\sqrt2}{m}\sqrt{\frac{\mfour T^2}{n\delta}}.
        \end{equation}
        \item (\textbf{Sub-exponential tails}) Now, suppose that Assumption~\ref{ass:subexp-norms} also holds. Then there is a universal constant $C > 0$ such that, whenever
        \(CK\left(\sqrt{\frac{\log(2T/\delta)}{B}}
        +\frac{\log(2T/\delta)}{B}\right) \leq mr_c\), we have that, with probability at least $1-\delta$, all iterates remain in $\Bc(\Ustar,r_c)$ and
        \begin{equation}\label{eq:thm2_subexp}
        \dframe(U_T,\Ustar)
        \leq\frac{r_{\mathrm{init}}}{\sqrt n}
        +\frac{CK}{m}\left(
        \sqrt{\frac{T\log(2T/\delta)}{n}}
        +\frac{T\log(2T/\delta)}{n}\right).
        \end{equation}
    \end{enumerate}
\end{theorem}
An expected-error counterpart to this result is stated and proved in Appendix~\ref{app:sgd-expected}.

\section{Analysis}\label{sec:analysis}
In this section, we prove Theorems~\ref{thm:convergence} and~\ref{thm:SGD_convergence_thm_prob}, both of which must address the nonsmoothness of the empirical risk $\Ec_n$, whose nearest-axis assignment changes discontinuously across switching boundaries, even though the population risk $\Ec$ is smooth near $\Norbit$.
\subsection{Empirical non-smooth minimization}\label{subsec:empirical-nonsmooth}

We first prove Theorem~\ref{thm:convergence}, for which we formalize the geometry of ties and switching to upgrade pointwise subgradient concentration to uniform control. Combined with localization, strong convexity, and the vanishing of $G_n$ at empirical minimizers, this uniform control yields the desired result.

\begin{definition}[Ambiguous points]
\label{def:ambiguous}
A point $x$ is \emph{ambiguous} at $U$ if $|J(x;U)|\geq2$. The corresponding set of ambiguous directions is \(\Ac(U):=\{v\in S^{d-1}:|J(v;U)|\geq2\}.\) For each labeled representative, we adopt the tie-breaking convention $j^\star(x;U):=\min J(x;U)$.
\end{definition}
For $j\in J(x;U)$, the Riemannian gradient of the active smooth branch is given by
\begin{equation}\label{eq:riemannian-gradient-fi}
  \gamma_j(x,U):=\Pi_U\!\left[-2(u_j^\top x)xe_j^\top\right]
  =(u_j^\top x)\bigl(u_j(U^\top x)^\top-xe_j^\top\bigr).
\end{equation}
Recall that $h(x,U)$ from \eqref{eq:active-set} averages the Euclidean branch gradients over $J(x;U)$. Accordingly, at an unambiguous point, the Riemannian gradient of $f_x(U):=\ell(x,U)$ is $\Rgrad f_x(U)=\gamma_{j^\star(x;U)}(x,U)$, whereas at a tie its Clarke subdifferential is the convex hull of the active gradients $\gamma_j(x,U)$; see Appendix~\ref{app:subdiff}. Moreover, by Lemma~\ref{lem:diff-implies-tie-null}, every nonzero sample is unambiguous almost surely at each fixed $U\in\Nc$, while all active gradients vanish when $x=0$.

For uniform analysis on $\Nc$, we select the canonical branch $g_{\mathrm{can}}(x,U):=\gamma_{j^\star(x;U)}(x,U)$ and define its empirical average $G_n(U):=n^{-1}\sum_i g_{\mathrm{can}}(x_i,U)$. By Proposition~\ref{prop:clarke-formulas}, $G_n(U)\in\partial_R\Ec_n(U)$. Although this tie-breaking rule depends on the labeling, it extends to $\Nc R$ through $g_R(x,UR):=g_{\mathrm{can}}(x,U)R$. Furthermore, since nonzero ties occur with probability zero at each fixed $U\in\Nc$, all active branch selections agree almost surely there, making the pointwise bounds selection-independent.

\subsubsection{Geometry of switching regions}
To analyze the variation of empirical subgradients across $\Nc$, we introduce the notion of switching points.

\begin{definition}[Switching points]
\label{def:switching-points}
For $U_1,U_2\in\Nc$, a sample $x_i$ \emph{switches} if $j^\star(x_i;U_1)\neq j^\star(x_i;U_2)$. When $\|U_1-U_2\|_F\leq\delta$, set
\[
N_\delta(U_1,U_2):=\#\{i:j^\star(x_i;U_1)\neq j^\star(x_i;U_2)\}.
\]
\end{definition}

The geometry of these switches is controlled by the ambiguous directions introduced above. In particular, a nonzero sample can switch only if its normalized direction lies within $(\pi/\sqrt2)\delta$ of $\Ac(U_1)$, which is itself contained in the union of the $d(d-1)$ signed pairwise bisector spheres; see Corollary~\ref{cor:switching-localization}. To obtain uniform control, at each scale $\delta$ we fix, before observing the data, a collection of packing centers $v_j$, each associated with a stripe such that these stripes cover the tubular neighborhoods around every signed bisector sphere; see Definition~\ref{def:universal-stripes}. Bernstein concentration then controls all stripes simultaneously, yielding a bound that applies even to data-dependent frame pairs without a union bound over sample-dependent switching sets.

\begin{theorem}[Switching concentration under a probability margin] \label{thm:switching-margin}
Under Assumptions~\ref{ass:strong-convexity} and~\ref{ass:switching-margin}, let $\Nc$ be the aligned component containing $\Ustar$, let $x_1,\ldots,x_n$ be i.i.d., and fix $\delta=n^{-1/2}$. Assume $n>\max\{16\eta_0^{-2},4\Delta_{\Nc}^{-2}\}$, where $\Delta_{\Nc}:=c_{\Nc}-\sqrt{1-c_{\Nc}^2}$ as in Theorem~\ref{thm:convergence}. Set $\mathcal J_n:=\{j:\dist(v_j,\mathcal W_{\Nc})\leq\delta\}$, $\displaystyle B_n:=4\Cmar\sqrt{n}+6\sqrt{\Cmar d\log n}\,n^{1/4}+9d\log n$, and $p_{\mathrm{sw}}(n):=2(d-1)\pi^{d-2}\omega_d\omega_{d-1}^{-1}n^{-(d-1)}$. With probability at least $1-p_{\mathrm{sw}}(n)$, every stripe indexed by $\mathcal J_n$ has far-region sample count at most $B_n$, and consequently
\[
N_\delta(U_1,U_2)\leq d(d-1)B_n,
\]
simultaneously for all $U_1,U_2\in\Nc$ satisfying $\|U_1-U_2\|_{F}\leq\delta$.
\end{theorem}

This switching bound is key for the coarsely Lipschitz property established next.

\subsubsection{Coarsely Lipschitz property and uniform convergence}

\begin{proposition}[Coarsely Lipschitz property of empirical subgradients]
\label{prop:quasi-lipschitz}
Under Assumptions~\ref{ass:concentration}, \ref{ass:subexp-norms}, \ref{ass:strong-convexity}, and~\ref{ass:switching-margin}, suppose $n$ satisfies the threshold in Theorem~\ref{thm:switching-margin} and $n\geq\dfrac{K^2}{C_B\mtwo^2}\log(8n)$, where $C_B>0$ is an absolute constant. With probability at least $1-2/n-p_{\mathrm{sw}}(n)$, simultaneously for all $U_1,U_2\in\Nc$ satisfying $\|U_1-U_2\|_F\leq n^{-1/2}$, we have
\begin{equation}\label{eq:quasi-lipschitz}
  \|G_n(U_1)-G_n(U_2)\|_{\Fr}
  \leq 4\mtwo\|U_1-U_2\|_F+\mathfrak J_n
  \quad\text{with}\quad
  \mathfrak J_n:=\frac{12Kd(d-1)B_n\log n}{n}.
\end{equation}
\end{proposition}

For nonswitching samples, the selected axis is unchanged and the contribution varies Lipschitzly with the frame. A switching sample may contribute $\bigO(\|x_i\|_2^2)$ independently of the frame distance, but Theorem~\ref{thm:switching-margin} controls how many such samples occur and the maximum-norm bound (Lemma~\ref{lem:max-norm}) controls their size. This produces a Lipschitz bound with a small additive residual. An $n^{-1/2}$-net then upgrades pointwise concentration to the uniform bound below; see Sections~\ref{app:quasilip} and~\ref{app:uniform-proof} of Appendix~\ref{app:erm-proofs}.

\begin{proposition}[Uniform empirical subgradient control]
\label{prop:uniform-convergence}
Under Assumptions~\ref{ass:concentration}, \ref{ass:subexp-norms}, \ref{ass:strong-convexity}, and~\ref{ass:switching-margin}, let $\Nc$ be the aligned component containing $\Ustar$, and assume the sample-size conditions of Proposition~\ref{prop:quasi-lipschitz}. Then there is an event $\Egood$ with $\PP(\Egood)\geq1-3/n-p_{\mathrm{sw}}(n)$ on which
\begin{align*}
    \sup_{U\in\Nc}\|\Rgrad\Ec(U)-G_n(U)\|_{\Fr}    &\leq\Bigl(2\left (\sqrt{8\mfour}+\sqrt2\,\Gamma\right )\sqrt{\Lambda_d}+L(1+r_{\Nc})+4\mtwo\Bigr)\sqrt{\tfrac{\log n}{n}}\\
    &\qquad+\mathfrak J_n+4\Gamma\Lambda_d\frac{\log n}{n}.
\end{align*}
Here $\Gamma:=C_1(2K+2C_0\sqrt{\mfour})$ and $\Lambda_d:=2+d(d-1)/4+d(d-1)\log(1+4r_{\Nc})/2$.
\end{proposition}

At an empirical minimizer, the selection controlled by this bound vanishes exactly. 

\begin{lemma}[Empirical local minimizers are smooth points]\label{lem:empirical-local-min}
    Fix any samples $x_1,\ldots,x_n$, and let $V$ be a local minimizer of $\Ec_n$ on $\orth(d)$. Then every nonzero sample is unambiguous at $V$, $\Ec_n$ is smooth on a neighborhood of $V$, and $\partial_R\Ec_n(V)=\{G_n(V)\}=\{0\}$.
\end{lemma}

The proof is given in Appendix~\ref{app:subdiff}. Although $\Ec_n$ is smooth near $V$, both $V$ and the radius of its smoothness neighborhood are sample dependent, and that radius cannot be bounded below independently of the sample. The switching analysis is therefore needed to obtain control of $G_n$ over $\Nc$ before evaluating it at the data-dependent minimizer.

\subsubsection{Localization of the empirical minimizer}

With uniform convergence established, we show that every empirical minimizer lies in $\Nc$ up to a signed permutation. It is low-risk for the population objective, so deterministic sublevel containment (Assumption~\ref{ass:sublevel-containment}) places it in $\Norbit$, whose disjoint decomposition then supplies the aligning permutation.

\begin{corollary}[Localization of empirical minimizer]
\label{cor:localization}
Under Assumptions~\ref{ass:subexp-norms}, \ref{ass:strong-convexity}, and~\ref{ass:sublevel-containment}, let $\rho>0$ satisfy $\Srho\subseteq\Norbit$. If $K=0$, the conclusion below holds surely. Otherwise set $\rho_0:=\min\{\rho,K\}$. There exists an absolute constant $\Cstat>1$ such that, for every $\beta\in(0,1)$, if
\[
n\geq\Cstat\frac{K^2}{\rho_0^2}
\left[d^2\log\!\left(\frac{\Cstat\sqrt d\,K}{\rho_0}\right)+\log\!\left(\frac{\Cstat}{\beta}\right)\right],
\]
then with probability at least $1-\beta$, every empirical minimizer $\Ustar_n$ belongs to $\Norbit$. Consequently, there is a unique $R_n\in\Spm$ such that the aligned $\Ustarntilde:=\Ustar_nR_n$ belongs to $\Nc$.
\end{corollary}

This threshold is exactly the fourth condition of Theorem~\ref{thm:convergence}.

\subsubsection{Proof of Theorem~\ref{thm:convergence}}\label{subsubsec:erm-proof}

Let $\Eloc$ denote the event in Corollary~\ref{cor:localization} on which every minimizer belongs to $\Norbit$, and $\Egood$ the event in Proposition~\ref{prop:uniform-convergence}. A union bound gives
\[
\PP(\Eloc\cap\Egood)\geq1-\beta-\frac{3}{n}
-2(d-1)\pi^{d-2}\frac{\omega_d}{\omega_{d-1}}n^{-(d-1)}.
\]
Fix an empirical minimizer $\Ustar_n$. On $\Eloc$, let $R_n$ be the unique signed permutation for which $\Ustarntilde:=\Ustar_nR_n\in\Nc$. By signed-permutation invariance of $\Ec_n$, the aligned representative $\Ustarntilde$ is also an empirical minimizer.

On $\Eloc$, the aligned minimizer $\Ustarntilde$ lies in $\Nc$, where $\Ec$ is geodesically $m$-strongly convex. Integrating the corresponding strong-monotonicity inequality along the minimizing geodesic from $\Ustar$ to $\Ustarntilde$ yields
\(
d_g(\Ustarntilde,\Ustar)\leq\frac{1}{m}\|\Rgrad\Ec(\Ustarntilde)\|_{\Fr}.
\)
Since $\Ustarntilde$ minimizes $\Ec_n$, Lemma~\ref{lem:empirical-local-min} gives $G_n(\Ustarntilde)=0$, so $\|\Rgrad\Ec(\Ustarntilde)\|_{\Fr}\leq\sup_{U\in\Nc}\|\Rgrad\Ec(U)-G_n(U)\|_{\Fr}$. On $\Eloc\cap\Egood$, Proposition~\ref{prop:uniform-convergence} bounds the right-hand side, and expanding $B_n$ in its switching term $\mathfrak J_n=12Kd(d-1)B_n\log n/n$ gives
\begin{align*}
48K\Cmar d(d-1)\frac{\log n}{\sqrt n} +72Kd(d-1)\sqrt{\Cmar d}\frac{(\log n)^{3/2}}{n^{3/4}} +108Kd^2(d-1)\frac{(\log n)^2}{n}.
\end{align*}
Using the trivial bounds $\sqrt{\tfrac{\log n}{n}}\leq\tfrac{\log n}{\sqrt n}$ and $\tfrac{\log n}{n}\leq\tfrac{(\log n)^2}{n}$ for the remaining terms gives
\[
\|\Rgrad\Ec(\Ustarntilde)\|_{\Fr}
\leq A_1\frac{\log n}{\sqrt n}+A_2\frac{(\log n)^{3/2}}{n^{3/4}}+A_3\frac{(\log n)^2}{n}.
\]
Combining this with strong convexity and $d_g(\Ustarntilde,\Ustar)=\dframe(\Ustar_n,\Ustar)$ proves the theorem. The bounds are uniform, so the conclusion holds for every empirical minimizer.
\qed

\subsection{Stochastic gradient descent}\label{subsec:stoch-grad-desc}

With the empirical minimizer localized, we turn to online SGD, which must keep every iterate inside the strongly convex component rather than localize a single minimizer. Equivariance of $g$ makes the induced frame iteration independent of the initial signed-permutation representative. Throughout, choose $r_c>0$ such that the closed geodesic ball $\Bc(\Ustar,r_c)$ is contained in $\Nc$, and align the deterministic initialization so that $d_g(U_0,\Ustar)<r_c$. The proof of Theorem~\ref{thm:SGD_convergence_thm_prob} has three steps: a dimension-free uniform stochastic-gradient error bound, a sharp local contraction, and an invariant-ball argument that closes the recursion along the entire trajectory.

\subsubsection{Dimension-free stochastic-gradient concentration}

Let $\F_t$ be the $\sigma$-algebra $\F_t:=\sigma(U_0,z_0,\ldots,z_{t-1})$. By the fresh-batch construction, $g(U_t,z_t)$ has conditional mean $g_\infty(U_t):=\EE_z[g(U_t,z)]$ with respect to $\F_t$; Appendix~\ref{app:SGD} shows $g_\infty(U)=\Rgrad\Ec(U)$ for $U\in\Nc$ (Lemma~\ref{lem:sgd-unbiased}).

\begin{lemma}[Dimension-free stochastic-gradient concentration]
\label{lemma:SG_bound}
For $U\in\orth(d)$ and $x\in\RR^d$, define
\[
\zeta(U,x):=\Pi_U h(x,U),
\qquad\text{and}\qquad
e_t:=g(U_t,z_t)-g_\infty(U_t).
\]
Under Assumption~\ref{ass:concentration},
\begin{equation}\label{eq:sg-second-moment}
\EE[\|e_t\|_{\Fr}^2\mid\F_t]\leq\frac{4\mfour}{B}.
\end{equation}
Consequently, for every $\delta\in(0,1)$,
\begin{equation}\label{eq:eps_B_delta_poly}
\PP\!\left(
\max_{0\leq t<T}\|e_t\|_{\Fr}
\leq2\sqrt{\frac{\mfour T}{B\delta}}
\right)\geq1-\delta.
\end{equation}
If Assumption~\ref{ass:subexp-norms} also holds, then there exists an absolute constant $C>0$ such that, for every $\delta\in(0,1)$,
\begin{equation}\label{eq:eps_B_delta_subexp}
\PP\!\left(
\max_{0\leq t<T}\|e_t\|_{\Fr}
\leq CK\left(
\sqrt{\frac{\log(2T/\delta)}{B}}
+\frac{\log(2T/\delta)}{B}\right)
\right)\geq1-\delta.
\end{equation}
\end{lemma}

Conditionally on $\F_t$, the $B$ centered summands $\zeta(U_t,x_{t,i})-g_\infty(U_t)$ lie in the fixed Hilbert space $T_{U_t}\orth(d)$. Their second moments give the finite-moment bound, while Hilbert-space Bernstein concentration gives the sub-exponential bound; a union bound over $0\leq t<T$ then makes each estimate uniform in time. A proof of Lemma~\ref{lemma:SG_bound} is given in Appendix~\ref{app:SGD}.

\subsubsection{Sharp local contraction and invariant-ball control}

Under Assumption~\ref{ass:strong-convexity}, the population gradient admits an averaged-Hessian representation, and $g(U_t,z_t)=\Rgrad\Ec(U_t)+e_t$ on $\Bc(\Ustar,r_c)$ (Lemma~\ref{lem:sgd-unbiased}); the exponential-map comparison then bounds the geodesic perturbation by the Frobenius norm, yielding the following one-step recursion.

\begin{lemma}[Sharp stochastic contraction]
\label{lem:sgd-sharp-contraction}
For $U_t\in\Bc(\Ustar,r_c)$, set
$r_t:=d_g(U_t,\Ustar)$. For
\[
q_\eta:=\max_{\lambda\in[m,L]}|1-\eta\lambda|,
\]
the update~\eqref{SGD_update} satisfies
\begin{equation}\label{eq:sharp-linear-recursion}
r_{t+1}\leq q_\eta r_t+\eta\|e_t\|_{\Fr}.
\end{equation}
For $0<\eta\leq2/(m+L)$,
$q_\eta=1-\eta m$. In particular, for
\(
\eta=\frac{2}{m+L}\) and
\(
q:=q_\eta=\frac{L-m}{L+m}
=\frac{\kappa-1}{\kappa+1},
\)
we have $1-q=\eta m$.
\end{lemma}

\begin{lemma}[Invariant strongly convex ball]
\label{lem:sgd-invariant-ball}
Let $0<\eta\leq2/(m+L)$; recall $q_\eta=1-\eta m$ from Lemma~\ref{lem:sgd-sharp-contraction}. If
$r_{\mathrm{init}}:=d_g(U_0,\Ustar)<r_c$ and
\[
\max_{0\leq t<T}\|e_t\|_{\Fr}\leq\varepsilon
\qquad\text{with}\qquad
\varepsilon\leq mr_c,
\]
then every $U_t$, $0\leq t\leq T$, belongs to
$\Bc(\Ustar,r_c)$ and
\begin{equation}\label{eq:invariant-ball-recursion}
r_t\leq q_\eta^tr_{\mathrm{init}}
+\frac{\varepsilon}{m}(1-q_\eta^t).
\end{equation}
\end{lemma}

Lemma~\ref{lem:sgd-sharp-contraction} separates population contraction from stochastic perturbation. Lemma~\ref{lem:sgd-invariant-ball} iterates that recursion and shows that a uniform error threshold $\varepsilon\leq mr_c$ keeps the trajectory in the same strongly convex ball. Proofs of Lemmas~\ref{lem:sgd-sharp-contraction} and~\ref{lem:sgd-invariant-ball} appear in Appendix~\ref{app:SGD}.

\subsubsection{Proof of Theorem~\ref{thm:SGD_convergence_thm_prob}}\label{subsubsec:sgd-high-prob-proof}
In the sub-exponential case, choose the universal constant $C$ in the theorem to be at least twice the concentration constant in Lemma~\ref{lemma:SG_bound}. Fix either regime and apply that lemma with
\[
\varepsilon:=
\begin{cases}
2\sqrt{\dfrac{\mfour T}{B\delta}}, & \text{in the finite-moment case},\\[1ex]
\dfrac{CK}{2}\left(\sqrt{\dfrac{\log(2T/\delta)}{B}}
+\dfrac{\log(2T/\delta)}{B}\right), & \text{in the sub-exponential case}.
\end{cases}
\]
With probability at least $1-\delta$, $\max_{0\leq t<T}\|e_t\|_{\Fr}\leq\varepsilon$. The corresponding assumption in the theorem implies $\varepsilon\leq mr_c$. Lemma~\ref{lem:sgd-invariant-ball} therefore keeps every iterate in $\Bc(\Ustar,r_c)$ and yields
\[
\dframe(U_T,\Ustar)=r_T
\leq q^Tr_{\mathrm{init}}+\frac{\varepsilon}{m}(1-q^T)
\leq q^Tr_{\mathrm{init}}+\frac{\varepsilon}{m}.
\]
Thus the same event controls both confinement and terminal error, with no additional escape event. Here the equality uses $U_T\in\Bc(\Ustar,r_c)\subseteq\Nc$. Moreover,
\[
q=\frac{\kappa-1}{\kappa+1}
\leq\exp\!\left(-\frac{2}{\kappa+1}\right),
\]
so $T=\lceil c_{\mathrm{iter}}\kappa\log n\rceil$, $c_{\mathrm{iter}}\geq1/2$, and $\kappa\geq1$ imply $q^T\leq n^{-1/2}$. Also, $n\geq2T$ gives $B=\lfloor n/T\rfloor\geq n/(2T)$. In the finite-moment regime, these estimates give
\[
\varepsilon
=2\sqrt{\frac{\mfour T}{B\delta}}
\leq2\sqrt2\sqrt{\frac{\mfour T^2}{n\delta}},
\]
which proves~\eqref{eq:thm2_poly}. In the sub-exponential regime, the same bound on $B$ gives
\[
\varepsilon
\leq CK\left(
\sqrt{\frac{T\log(2T/\delta)}{n}}
+\frac{T\log(2T/\delta)}{n}\right),
\]
which proves~\eqref{eq:thm2_subexp}, with the same constant $C$ in the basin condition and the terminal bound.
\qed

\section{Numerical Experiments}\label{sec:experiments}

In this section, we test whether the finite-sample guarantees of Section~\ref{sec:analysis} hold for stylized examples. A real-data illustration, on a text corpus not verified to satisfy our assumptions but exhibiting the same axis-concentration structure, is given in Appendix~\ref{app:realdata}. As noted in the introduction, our method is the $s=1$ case of fixed-sparsity orthogonal dictionary learning~\cite{lesage2005,rusu2016fast} (itself a constrained form of K-SVD~\cite{kSVD}), solved here by Riemannian gradient descent.

\subsection{Statistical convergence rates}
\label{sec:exp:rate}

Up to logarithmic factors, Theorem~\ref{thm:convergence} predicts an $n^{-1/2}$ frame-error rate. Defining the relative error as $\varepsilon = \dframe(\Ustar_n,\Ustar)/\|\Ustar\|_F$, this corresponds to a slope of $-1/2$ in a log-log plot of $\varepsilon$ versus $n$, or, equivalently, a $1/\sqrt{n}$ reference line. We sample data from an equal-weight mixture of $d$ components, each with mean zero and covariance assigning variance $1-t$ along coordinate axis $j$ and $t/(d-1)$ along all other axes, with $t = 0.05$. This instantiates the balanced Gaussian-mixture model of Example~\ref{ex:gaussian-convolution}, which Proposition~\ref{prop:balanced-gaussian-mixture-uniform} shows satisfies our assumptions for sufficiently large $d$. The double-exponential and Student-$t$ data have independent coordinates drawn from $\pm\mathrm{Exp}^{1.3}$ and $\pm t_3$ respectively; these distributions still concentrate near the frame's axes by construction; their non-Gaussian, heavier-tailed marginals instead test robustness of the convergence rate beyond our theoretical assumptions.

\begin{figure}[htpb!]
     \centering
     \includegraphics[width=\linewidth]{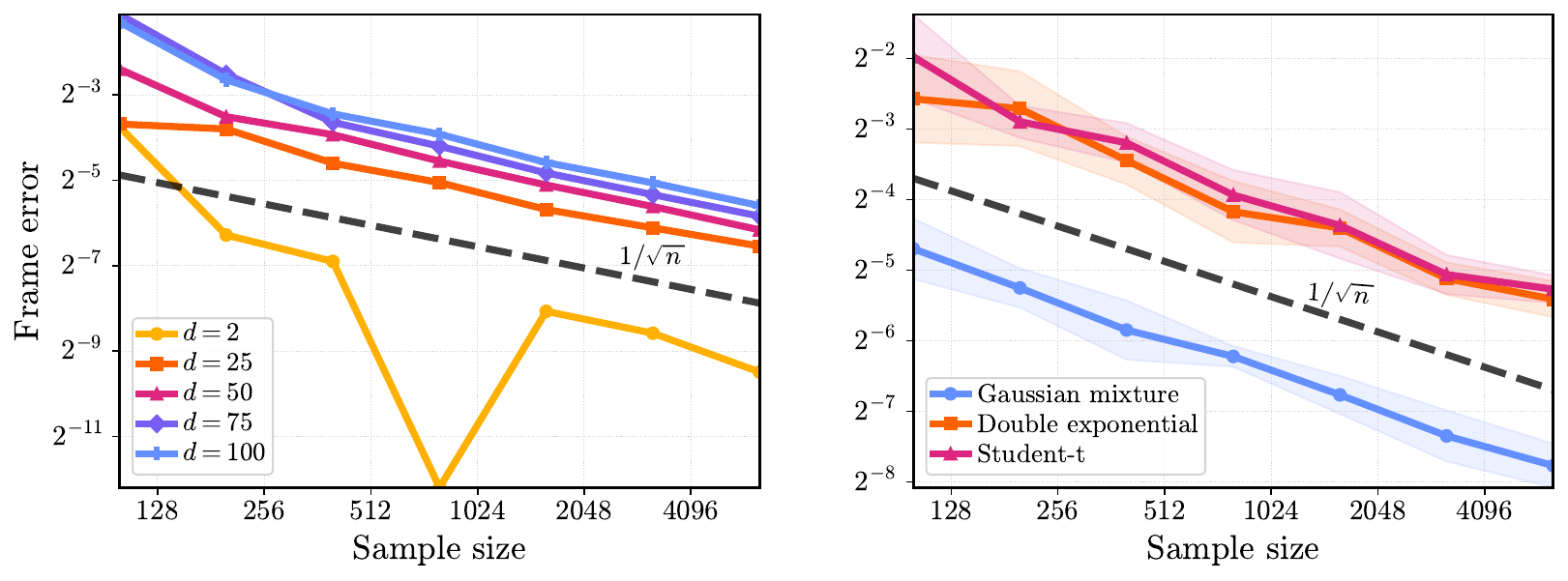}
     \caption{\small \textbf{Relative error versus sample size ($\log$-$\log$).} Relative error is defined as $\varepsilon = \dframe(\Ustar_n,\Ustar)/\|\Ustar\|_F$ after optimal sign and permutation alignment; the dashed line shows the $1/\sqrt{n}$ reference rate. \textit{Left:} Empirical rates for $d \in \{2, 25, 50, 75, 100\}$ under a Gaussian mixture. \textit{Right:} The same metric at fixed $d = 15$ across Gaussian mixture, double-exponential, and Student-$t$ ($\nu = 3$) distributions; shaded bands denote $\pm 1$ standard deviation over trials.}
     \label{fig:rates_combined}
\end{figure}
All curves decrease parallel to the $1/\sqrt n$ reference line across the tested sample sizes, supporting the predicted dependence on $n$. The vertical offsets show that the prefactor grows with $d$, but these plots do not identify its worst-case dimension dependence or isolate the distributional constants in Theorem~\ref{thm:convergence}.
The Student-$t$ distribution has an infinite fourth moment and violates Assumption~\ref{ass:concentration}. Nevertheless, all three curves are parallel to the $1/\sqrt n$ reference over the plotted range. These finite experiments support, but do not establish, robustness of the empirical slope to increasing dimension and distributional misspecification.

\subsection{Computational convergence}
\label{sec:exp:opt}

\begin{figure}[htpb!]
    \centering
    \includegraphics[width=\linewidth]{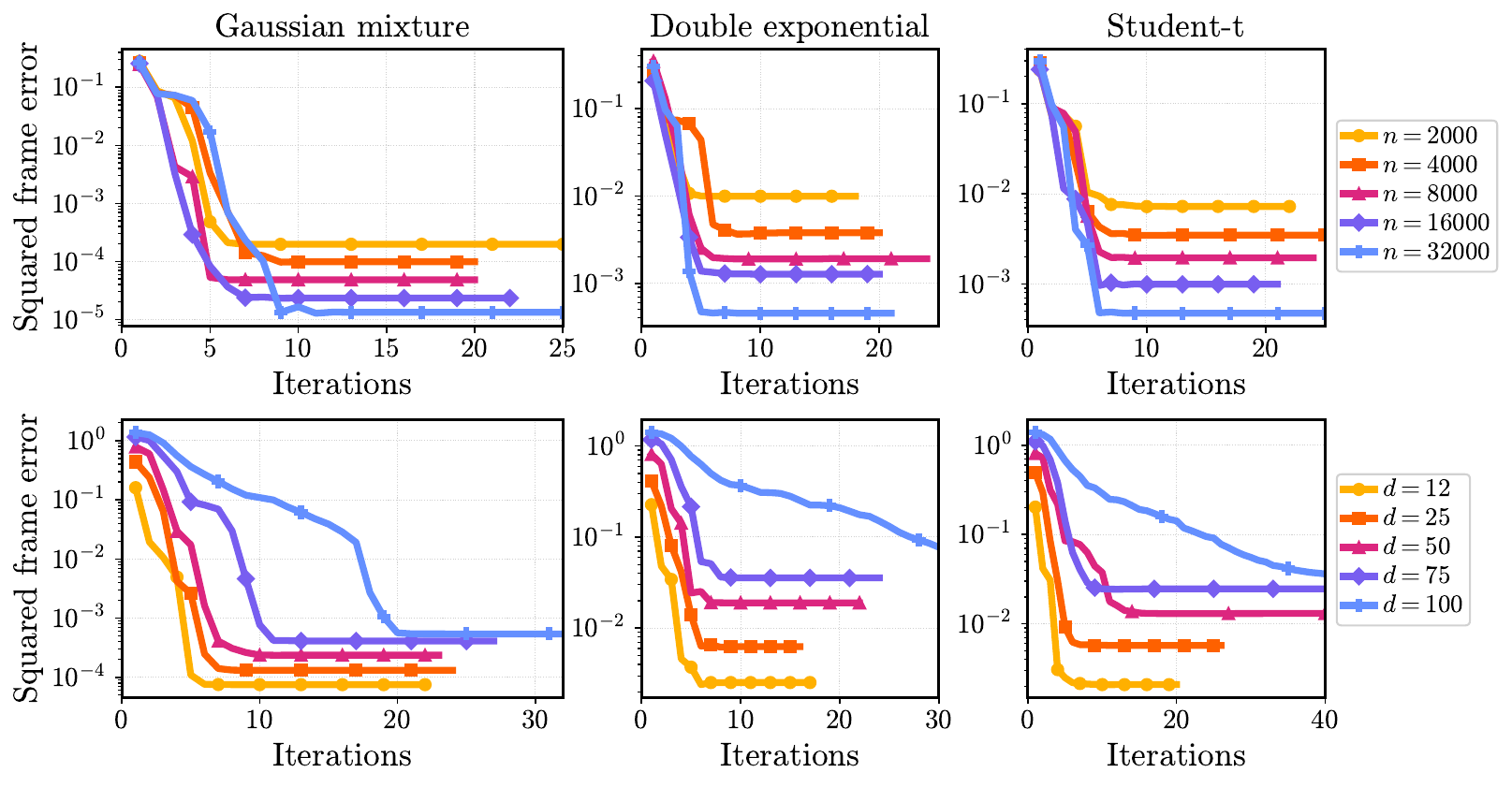}
    \caption{\small \textbf{Squared frame error versus iteration count.} The y-axis shows squared frame error on a logarithmic scale. \textit{Top row:} Fixed $d = 15$, varying $n \in \{2000, 4000, 8000, 16000, 32000\}$; curves plateau at a floor decreasing with $n$. \textit{Bottom row:} Fixed $n = 5000$, varying $d \in \{12, 25, 50, 75, 100\}$; the floor rises with $d$. Columns correspond to Gaussian mixture, double-exponential, and Student-$t$ ($\nu = 3$) data.}
    \label{fig:log_iterations}
\end{figure}

Theorem~\ref{thm:SGD_convergence_thm_prob} analyzes online SGD with a fresh independent batch at every iteration. For simplicity, Figure~\ref{fig:log_iterations} instead reuses one fixed empirical sample throughout, applying ordinary steepest descent to its empirical objective rather than redrawing batches at every step, initialized at $\exp(0.1\,A)\Ustar$ with $A$ skew-symmetric. The top row fixes $d=15$ and varies $n$; the bottom row fixes $n=5000$ and varies $d$; columns use the three distributions from Section~\ref{sec:exp:rate}.

The top-row curves show rapid initial decay followed by an empirical floor that decreases with $n$. This resembles the theorem's contraction-plus-noise-floor mechanism, but the protocol does not identify its constants or verify the predicted online-SGD squared-error rate $\bigO(T/B)$ for fixed problem parameters, equivalently $\bigO((\log n)^2/n)$ under the prescribed sample-budget allocation. The Gaussian mixture has lower floors than the two misspecified distributions.

The bottom row analyzes dimension at fixed sample size. For the Gaussian mixture, the two-phase structure persists and the noise floor increases monotonically with $d$, consistent with a dimension-dependent constant. For the double-exponential and Student-$t$ distributions, convergence remains qualitatively similar at lower dimensions but deteriorates for $d\geq75$; at $d=100$, a stable plateau is barely reached. Hence, the method remains effective at moderate dimensions, while the noise floor can grow faster under distributional misspecification.

\subsection{Comparison with competing methods}
\label{sec:exp:compare}

We compare Riemannian gradient descent with standard alternating-minimization~\cite{lesage2005,rusu2016fast}, which seeks the same estimator, and two further methods based on different objectives, PCA+Varimax~\cite{kaiser1958,rohe2023} and Robust Subspace Clustering (RSC)~\cite{SoltanolkotabiElhamifar2014}. Since in the first comparison we solve the same objective by two different numerical strategies, matching accuracy serves as a sanity check on the two implementations; the comparison with PCA+Varimax and RSC is where the estimator's accuracy, robustness, and runtime are tested against genuinely different approaches. The error panels report the squared discrepancy $\varrho(U,V)^2$, where $\varrho(U,V) := \min_{R \in \Spm}\|U-VR\|_{\Fr}$. Since $\varrho(U,V)^2 \leq \dframe(U,V)^2$ (Remark~\ref{rem:frob-geod}), Theorem~\ref{thm:convergence} also controls this metric after squaring its bound. All experiments use Gaussian mixture data with $n = 1000$.

\begin{figure}[htpb!]
    \centering
  \includegraphics[width=\linewidth]{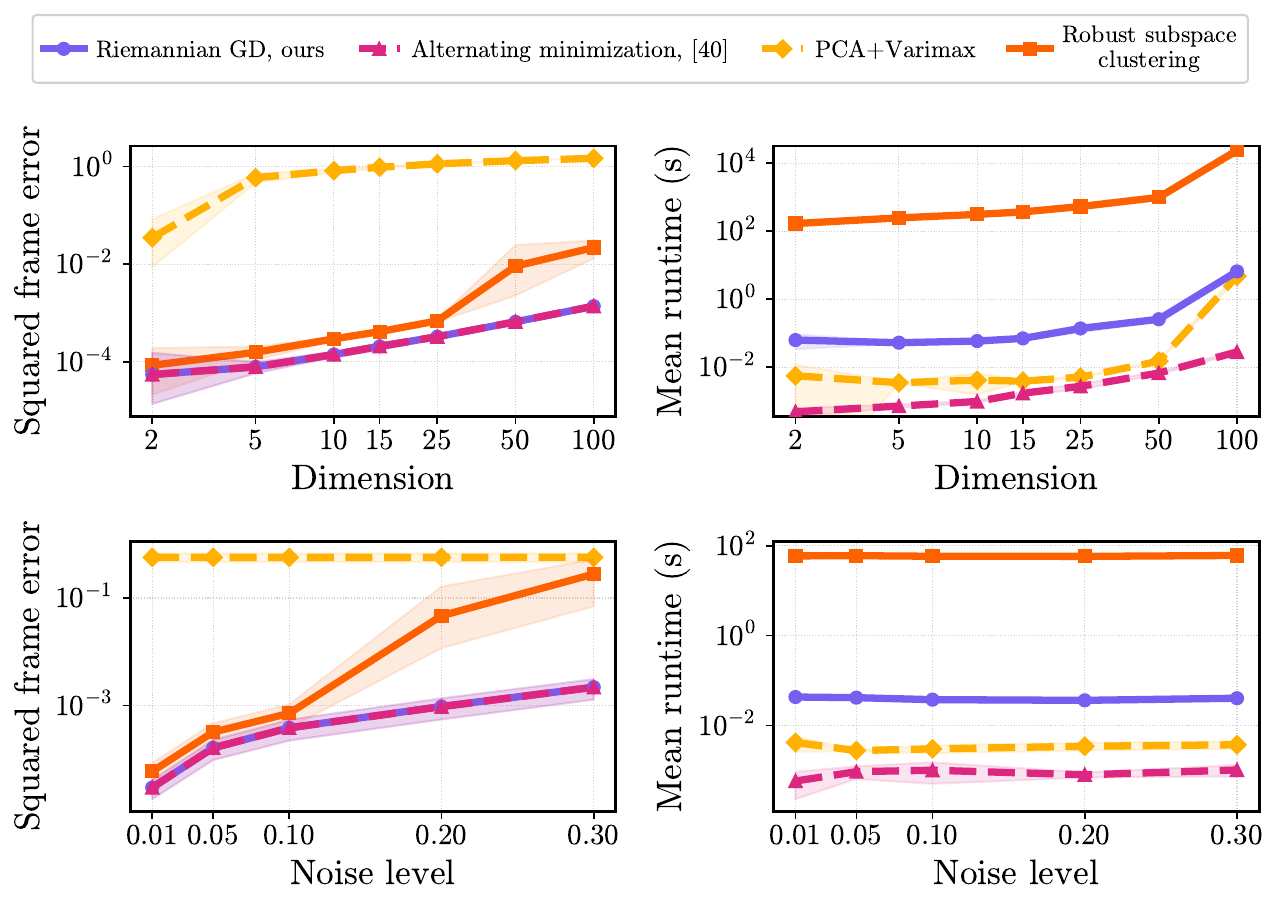}
    \caption{\small \textbf{Comparison with competing methods.}
    Shaded bands denote $\pm 1$ standard deviation over trials. \textit{Top-left:} Mean squared frame error versus $d \in \{2, 5, 10, 15, 25, 50, 100\}$. \textit{Top-right:} Mean runtime versus dimension. \textit{Bottom-left:} Mean squared frame error versus $\tau \in \{0.01, 0.05, 0.10, 0.20, 0.30\}$ at fixed $d = 15$. \textit{Bottom-right:} Mean runtime versus noise level.}
    \label{fig:comparison}
\end{figure}

\paragraph{Frame recovery accuracy}
In the plotted trials, $\name$ and its alternating-minimization counterpart overlap almost exactly across all dimensions, as expected from their shared objective at $s=1$; the largest absolute difference between the two across individual trials is approximately $4.4\times10^{-7}$. Both have substantially smaller squared frame error than RSC and PCA+Varimax. PCA+Varimax remains at errors of order one across all dimensions, while RSC's error grows steadily with $d$ and its band broadens at larger $d$. These observations describe the plotted implementation and do not certify global empirical minimizers or identify the cause of a method's error.

\paragraph{Noise robustness}
Across the noise sweep, $\name$ and its counterpart under alternating minimization track each other closely, with squared error rising with $\tau$ to approximately $2\times10^{-3}$ at $\tau=0.30$ (largest pairwise difference approximately $7.7\times10^{-5}$). RSC is comparable through $\tau=0.10$ and deteriorates markedly at $\tau=0.20$ and $0.30$, with a wider band. PCA+Varimax remains at errors of order one over the sweep.

\paragraph{Computational cost}
Despite matching accuracy, alternating minimization is substantially faster than Riemannian GD: at $d=100$, mean runtimes are approximately $0.028$ and $6.55$ seconds, respectively, and at noise level $\tau=0.30$ in $d=5$, approximately $0.0010$ and $0.0403$ seconds. This gap reflects a difference in two numerical strategies for the same target, not a difference in what is being estimated: alternating minimization exploits closed-form top-one hard-thresholding and an SVD-based Procrustes update, a solver specific to this fixed-sparsity, batch setting, while Riemannian GD uses general-purpose gradient and line-search computations that extend directly to the online, streaming-data analysis of Section~\ref{sec:analysis}. Both algorithms remain far faster than RSC, whose mean runtime exceeds $10^3$ seconds at $d=100$ and $60$ seconds across the noise sweep; PCA+Varimax is comparable in speed to alternating minimization at small $d$ but grows to match Riemannian GD's runtime by $d=100$. These measurements are hardware- and implementation-specific, not complexity guarantees.

\section*{AI disclaimer}
We prepared an initial draft of this manuscript before summer of 2026, developing the arguments without AI assistance. That draft contained all the results presented here, including the $\widetilde{\bigO}_d(n^{-1/2})$ rates for both Theorem~\ref{thm:convergence} and Theorem~\ref{thm:SGD_convergence_thm_prob}. Following the release of ChatGPT 5.6 Sol, we used this model to adversarially review the manuscript for technical gaps. Apart from minor corrections, the main theoretical change arising from this process was a refinement of a covering argument in the proof of Theorem~\ref{thm:convergence}, which improved the dimensional dependency in the rate from exponential to polynomial. We have also used ChatGPT to assist with coding, polish the writing, and conduct literature searches. We take full responsibility for the content of this paper.
\par

\bibliographystyle{abbrvnat}
\bibliography{biblio}
\appendix
\crefalias{section}{appendix}
\crefalias{subsection}{subappendix}
\crefalias{subsubsection}{subsubappendix}
\section{Auxiliary results common to the main theorems}\label{app:common-auxiliary}

This section collects four auxiliary results used throughout the proofs that follow.

\subsection{Wasserstein projection onto frame axes}\label{app:wasserstein-reduction}

The following proposition justifies the reduction, used already in the introduction, from the measure-approximation problem to the finite-dimensional objective $\Ec$.

\begin{proposition}[Nearest-axis Wasserstein reduction]
\label{prop:nearest-axis-wasserstein}
Fix $\mu\in\Pc_2(\RR^d)$ and $U\in\orth(d)$, and set $S_U:=\bigcup_{j=1}^d\Span\{u_j\}$. There exists a Borel nearest-axis projection $P_U\colon\RR^d\to S_U$, and every Borel nearest-axis projection satisfies
\begin{equation}\label{eq:fixed-frame-wasserstein-reduction}
W_2^2(\mu,P_U\#\mu)=\inf_{\nu\in\F_U}W_2^2(\mu,\nu)=\int_{\RR^d}\dist(x,S_U)^2\,d\mu(x)=\int_{\RR^d}\ell(x,U)\,d\mu(x).
\end{equation}
Consequently, $\inf_{\nu\in\F}W_2^2(\mu,\nu)=\min_{U\in\orth(d)}\Ec(U)$. For $\mu_n=n^{-1}\sum_{i=1}^n\delta_{x_i}$, the same identity gives $\inf_{\nu\in\F_U}W_2^2(\mu_n,\nu)=\Ec_n(U)$.
\end{proposition}

\begin{proof}
Let $s_j(x,U):=|u_j^\top x|^2$ and choose $j_U(x):=\min\argmax_j s_j(x,U)$. The selection is Borel because
\[
\{(x,U):j_U(x)=j\} =\bigcap_{k<j}\{s_j>s_k\}\cap\bigcap_{k>j}\{s_j\geq s_k\}.
\]
Hence $P_U(x):=(u_{j_U(x)}^\top x)u_{j_U(x)}$ is a Borel metric projection onto the closed set $S_U$. Applying \cite[Lemma~A.1 and (35)]{cutler2026radius} with $p=2$ and $K=S_U$ gives the first three terms of~\eqref{eq:fixed-frame-wasserstein-reduction}, for this or any other Borel nearest-axis selection. The last term follows from
\[
\dist(x,S_U)^2 =\|x\|^2-\max_j|u_j^\top x|^2 =\ell(x,U),
\]
proving~\eqref{eq:fixed-frame-wasserstein-reduction}. Since $\F=\bigcup_U\F_U$, taking the infimum over $U$ gives the global identity; moreover, $\ell(x,\cdot)$ is continuous and bounded by $\|x\|^2$, so dominated convergence and compactness of $\orth(d)$ give attainment. Applying the fixed-frame identity to $\mu_n$ yields the empirical claim. Finally, $S_{UR}=S_U$ for every signed permutation $R$, so the value is permutation-invariant.
\end{proof}

\subsection{Hilbert-space concentration}\label{app:hilbert-concentration}

The following concentration inequality for sums of independent Hilbert-space-valued random variables is a common ingredient in the empirical-minimizer and online-SGD proofs in Appendices~\ref{app:erm-proofs} and~\ref{app:SGD}.

\begin{lemma}[Bennett-type concentration for empirical subgradients]
\label{lem:pinelis-hilbert-bennett}
Let $Y_1,\ldots,Y_n$ be independent mean-zero random variables taking values in the tangent space $(T_U \orth(d), \|\cdot\|_{\Fr})$ at some fixed $U \in \orth(d)$. Suppose
there exist $\Gamma, \varsigma > 0$ such that, for all integers $m \geq 2$,
\begin{equation}\label{eq:pinelis-moment-growth}
  \sum_{k=1}^n\EE[\|Y_k\|_{\Fr}^m] \leq m!\,\Gamma^{m-2}\,\frac{\varsigma^2}{2}.
\end{equation}
Then, for all $r \geq 0$,
\begin{equation}\label{eq:pinelis-bennett-hilbert}
  \PP\!\left(\left\|\sum_{k=1}^n Y_k\right\|_{\Fr} \geq r\right) \leq 2\exp\!\left(-\frac{r^2}{\varsigma^2 + \Gamma r + \varsigma\sqrt{\varsigma^2+2\Gamma r}}\right).
\end{equation}
In particular, the right-hand side is at most $2\exp\bigl[-\tfrac14\min\{r^2/\varsigma^2,r/\Gamma\}\bigr]$.
\end{lemma}

\begin{proof}
The tangent space $T_U\orth(d)$ is a finite-dimensional Hilbert space with inner product $\langle v,w\rangle_\Fr = \Tr(v^\top w)$; every Hilbert space is $(2,1)$-smooth in the sense of Pinelis~\cite{pinelis1994optimum} with constant $D=1$, as the parallelogram identity is equivalent to $(2,1)$-smoothness. Set $\F_k:=\sigma(Y_1,\ldots,Y_k)$ and $f_k:=\sum_{i=1}^kY_i$. Independence and centering make $(f_k,\F_k)$ a martingale with differences $Y_k$. In the independent setting, $\EE_{k-1}[\|Y_k\|_\Fr^m] = \EE[\|Y_k\|_\Fr^m]$ almost surely, so the conditional moment sum in Pinelis' theorem coincides with the unconditional sum in~\eqref{eq:pinelis-moment-growth}. The stated bound~\eqref{eq:pinelis-bennett-hilbert} then follows from \cite[Theorem~3.3]{pinelis1994optimum}, with the correction in \cite{pinelis1999correction}. Moreover, $\sqrt{\varsigma^2+2\Gamma r}\leq\varsigma+\Gamma r/\varsigma$, so the denominator in~\eqref{eq:pinelis-bennett-hilbert} is at most $2\varsigma^2+2\Gamma r\leq4\max\{\varsigma^2,\Gamma r\}$, and the exponent is at least $\tfrac14\min\{r^2/\varsigma^2,r/\Gamma\}$.
\end{proof}

\subsection{Metric and exponential-map comparisons}\label{app:metric-comparisons}

The experiments in Section~\ref{sec:experiments} measure frame recovery using the squared Frobenius discrepancy $\varrho^2$, whereas our theoretical guarantees are expressed in terms of $\dframe$. The following comparison relates these two error measures.

\begin{remark}[Frobenius and geodesic distance]\label{rem:frob-geod}
For $U,V\in\orth(d)$ in the same connected component, $\|U-V\|_\Fr\leq d_g(U,V)$ \cite[Ch.~6]{lee2018riemannian}. Taking the minimum over signed permutations gives $\varrho(U,V)\leq\dframe(U,V)$.
\end{remark}

The online-SGD analysis also uses the following nonexpansiveness property of the exponential map on $\orth(d)$.

\begin{lemma}[Global exponential-map comparison]
\label{lem:orth-exp-nonexpansive}
For every $U\in\orth(d)$ and $v,w\in T_U\orth(d)$,
\[
d_g\!\left(\Exp_U(v),\Exp_U(w)\right)\leq\|v-w\|_{\Fr}.
\]
\end{lemma}

\begin{proof}
Write $v=U\Omega_1$ and $w=U\Omega_2$ with $\Omega_1,\Omega_2$ skew-symmetric, and set $X_s=\Omega_1+s(\Omega_2-\Omega_1)$ and $c(s)=Ue^{X_s}$.  The left-trivialized velocity is
\[
c(s)^\top\dot c(s)=\int_0^1e^{-tX_s}(\Omega_2-\Omega_1)e^{tX_s}\,dt.
\]
Orthogonal conjugation preserves the Frobenius norm, so $\|\dot c(s)\|_{\Fr}\leq\|\Omega_2-\Omega_1\|_{\Fr}$.  Hence the length of $c$ is at most $\|\Omega_2-\Omega_1\|_{\Fr}=\|v-w\|_{\Fr}$, which bounds the distance between its endpoints.
\end{proof}

For tangent vectors of norm at most the injectivity radius, the estimate in Lemma~\ref{lem:orth-exp-nonexpansive} also follows from Toponogov's comparison theorem, since $\orth(d)$ has nonnegative sectional curvature. The proof above gives it for all $v,w\in T_U\orth(d)$, however large. This matters in the proof of Lemma~\ref{lem:sgd-sharp-contraction}.

\subsection{Active-tie nullity from differentiability}\label{app:diff-implies-tie-null}

We use the following lemma in the pointwise concentration of Appendix~\ref{app:pointwise}, in the analysis of empirical minimizers in Appendix~\ref{app:subdiff}, and in the online-SGD analysis of Appendix~\ref{app:SGD}. 

\begin{lemma}[Differentiability implies active tie nullity]\label{lem:diff-implies-tie-null}
    Let $\mu\in\Pc_2(\RR^d)$, and suppose that $\Ec$ is differentiable at $U\in\orth(d).$ Then
    \[
    \PP_{x\sim\mu}\bigl (x\neq0\text{ and }|J(x;U)|\geq2 \bigr)=0.
    \]
    Consequently, every measurable selection of the active gradients $\gamma_j(x,U)$, $j\in J(x;U)$, agrees with $g_{\mathrm{can}}(x,U)$ almost surely and has expectation $\Rgrad\Ec(U)$. Moreover, the ambient extension $\bar\Ec$ is differentiable at $U$ with 
    \[
     \nabla\bar\Ec(U)=\EE_{x\sim\mu}[h(x,U)].
    \]
\end{lemma}

Under Assumptions~\ref{ass:concentration} and~\ref{ass:strong-convexity}, the lemma applies at every $U\in\Nc$. The null set depends on $U$, which suffices for fixed-frame concentration and, conditionally on $\F_t$, for the fresh-batch SGD analysis.

\begin{proof}
    \emph{Directional derivatives.} Fix $v\in T_U\orth(d)$. By~\eqref{eq:riemannian-gradient-fi}, the smooth branch $\|x\|_2^2-(u_j^{\top}x)^2$ has a Riemannian gradient $\gamma_j(x,U)$, so its derivative along $t\mapsto\Exp_U(tv)$ at $t=0$ is $\langle\gamma_j(x,U), v\rangle_{\Fr}$. Branches inactive at $U$, if any, remain inactive for small $t$. Hence, as $t\downarrow0$ we have 
    \[
    \ell(x,\Exp_U(tv))-\ell(x,U) = t\min_{j\in J(x;U)} \langle\gamma_j(x,U), v\rangle_{\Fr}+\littleo(t).
    \]
    Each score $(u_j^\top x)^2$ is $2\|x\|_2^2$-Lipschitz in the unit vector $u_j$, and a finite minimum preserves a common Lipschitz constant. Together with $\|\Exp_U(tv)-U\|_{\Fr}\leq t\|v\|_\Fr$ from Remark~\ref{rem:frob-geod}, this shows that the difference quotients $t^{-1}\left (\ell(x,\Exp_U(tv))-\ell(x,U)\right )$ are bounded by the integrable envelope $2\|x\|_2^2\|v\|_\Fr$. Dominated convergence and differentiability of $\Ec$ at $U$ therefore give 
    \begin{equation}\label{eq:tie-directional}
        \langle\Rgrad\Ec(U),v\rangle_\Fr = \EE\left [\min_{j\in J(x;U)}\langle\gamma_j(x,U),v \rangle_\Fr\right].
    \end{equation}
    Adding~\eqref{eq:tie-directional} for $v$ and $-v$, and using $\min_j(-a_j)=-\max_ja_j$, yields
    \[
    \EE\left[ \max_{j\in J(x;U)}\langle\gamma_j(x,U),v\rangle_{\Fr} - \min_{j\in J(x;U)}\langle\gamma_j(x,U),v\rangle_{\Fr}\right]=0,
    \]
    so the nonnegative integrand vanishes almost surely. 

    \emph{Tie nullity.} Fix $i<j$ and take $v:=u_je_i^\top-u_ie_j^\top$, which belongs to $T_U\orth(d)$ because $U^{\top}v=e_je_i^\top-e_ie_j^{\top}$ is skew-symmetric. Moreover, $ve_i=u_j$ and $ve_j=-u_i$. Since $\langle\gamma_k(x,U),v\rangle_{\Fr}=-2(u_k^\top x)(x^\top ve_k)$ by~\eqref{eq:riemannian-gradient-fi}, we obtain 
    \begin{equation}\label{eq:tie-rotation}
        \bigl\langle\gamma_i(x,U)-\gamma_j(x,U),\,u_je_i^\top-u_ie_j^\top\bigr\rangle_{\Fr}=-4(u_i^\top x)(u_j^\top x).
    \end{equation}
    If $x\neq0$ and $i,j\in J(x;U)$, then $|u_i^\top x|=|u_j^\top x|=\max_k|u_k^\top x|>0$, so, for the direction $v$ above, 
    \[
    \begin{aligned}
        &\max_{k\in J(x;U)}\langle \gamma_k(x,U),v\rangle_\Fr - \min_{k\in J(x;U)}\langle \gamma_k(x,U),v\rangle_\Fr\\
        &\qquad \geq \left |\langle \gamma_i(x,U)-\gamma_j(x,U),v\rangle_\Fr\right | = 4|(u_i^\top x)(u_j^\top x)|=4(u_i^\top x)^2>0.
    \end{aligned}
    \]
    Since this integrand vanishes almost surely by the first step, we conclude that
    $\PP(x\neq0,\ i,j\in J(x;U))=0$, and a union over the $d(d-1)/2$ pairs proves the first claim.

    \emph{Consequence.} Almost surely, either $J(x;U)$ is a singleton or $x=0$, in which case every $\gamma_j(x,U)$ vanishes. Hence every active gradient selection $\gamma_j(x,U)$ with $j\in J(x;U)$ agrees with $g_{\mathrm{can}}(x,U)$ almost surely. The integrand in~\eqref{eq:tie-directional} therefore equals $\langle g_{\mathrm{can}}(x,U),v\rangle_\Fr$ almost surely. Since both $\Rgrad\Ec(U)$ and $\EE[g_{\mathrm{can}}(x,U)]$ lie in $T_U\orth(d)$, they coincide, which proves our claims on selections. Likewise, for almost every $x\neq0$, the scores $A\mapsto(a_j^\top x)^2$ keep a unique maximizer on a neighborhood of $U$ in $\RR^{d\times d}$, while at $x=0$ the integrand of $\bar\Ec$ is identically zero and $h(0,U)=0$. Thus, the integrand of $\bar\Ec$ is differentiable at $U$ almost surely with gradient $h(x,U)$. Near $U$, it is $3\|x\|^2_2$-Lipschitz in the Frobenius norm, so dominated differentiation gives $\nabla\bar\Ec(U)=\EE[h(x,U)]$.
\end{proof}

\section{Verification of the Gaussian convolution model}
\label{appendix:gaussian-proofs}

This appendix verifies the assumptions of our analysis for the Gaussian convolution model of Example~\ref{ex:gaussian-convolution}. We first establish Proposition~\ref{prop:gaussian-convolution-assumptions} for general axial laws at each fixed dimension, under sufficiently small noise. For balanced Gaussian mixtures, the subsections on dimension-uniform curvature and the global gap develop sharper estimates that accommodate the noise scaling $\tau^2\leq b_+/d$.

For $\lambda\in\{\mu,\nu\}$, write $\Ec_\lambda$ for the population risk associated with $\lambda$. For each $r>0$, let $\Nc=B_r(\Ustar)$ and consider its orbit $\Norbit$ as defined in Section~\ref{sec:prob_setting}. Write $\Delta_\nu:=\inf_{U\notin\Norbit}(\Ec_\nu(U)-\Ec_\nu(\Ustar))$ and $\bar a:=\sum_j a_j$.

\begin{proposition}[Gaussian convolution model satisfies the assumptions]
\label{prop:gaussian-convolution-assumptions}
Fix $d\geq2$. There exists $r_0\in(0,1/8]$ such that, for every $0<r\leq r_0$, we have $\Delta_\nu>0$ and $c_{\Nc}>1/\sqrt2$. Choose the constants
\[
K\leq2\max_j\|S_j^2\|_{\psi_1}+2\tau^2\bigl\|\|Z\|^2\bigr\|_{\psi_1},
\qquad
m=\frac{a_{\min}}{2},
\qquad
L=4a_{\max},
\qquad\text{and}\qquad
\rho=\frac{\Delta_\nu}{4},
\]
and the noise threshold
\[
\tau_0:=\min\{\tau_1,\tau_2\},
\quad\text{with}\quad
\tau_1:=\Bigl(\frac m C\Bigr)^{1/3}\quad\text{and}
\quad
\tau_2:=\min\Bigl\{1,\frac{\Delta_\nu}{8\sqrt{2d(2\bar a+d)}}\Bigr\},
\]
for $C=C(d,M)>0$. Then, for every $0<\tau<\tau_0$, there exists $\eta_0\in(0,1/2]$ satisfying $1/\sqrt2+r+\eta_0<1$ when $d\geq4$, such that Assumptions~\ref{ass:concentration}--\ref{ass:switching-margin} hold on $\Nc$ with the constants above and some $\Cmar>0$, and
\[
\arg\min_{U\in\orth(d)}\Ec(U)=\Ustar\Spm.
\]
\end{proposition}

The proof is given in the final subsection of this appendix, together with that of Proposition~\ref{prop:balanced-gaussian-mixture-uniform}, after establishing the necessary curvature, sublevel-set containment, and switching-margin estimates.

\subsection{Variational formulas and noiseless curvature}

For $\lambda\in\{\mu,\nu\}$, signed-permutation invariance gives, for every $R\in\Spm$,
\[
\Ec_\lambda(UR)=\Ec_\lambda(U),\quad
\Rgrad\Ec_\lambda(UR)=\Rgrad\Ec_\lambda(U)R,
\quad\text{and}\quad
\Hess\,\Ec_\lambda(UR)[\xi R]
=\Hess\,\Ec_\lambda(U)[\xi]R,
\]
where the derivative identities hold wherever the corresponding derivatives exist, and 
\[
\Hess\,\Ec_\lambda(U)[\xi,\xi]:=\langle\Hess\,\Ec_\lambda(U)[\xi],\xi\rangle. 
\]
Consequently, local curvature estimates near $\Ustar$ extend to every signed-permutation representative.

We first derive the variations of the population risk along general curves in $\orth(d)$. For a curve $\gamma(t)=[u_1(t),\ldots,u_d(t)]$, Proposition~\ref{prop:nearest-axis-wasserstein} gives
\begin{equation}\label{WFA1-energy}
\Ec(\gamma(t))
=\int_{\RR^d}\min_{1\leq j\leq d}
\left(\|x\|_2^2-\langle u_j(t),x\rangle^2\right)\,d\mu(x).
\end{equation}
The next two lemmas differentiate this finite minimum, first within the active cells and then accounting for their moving boundaries. We write
\[
\ell_i(t,x):=\|x\|_2^2-\langle u_i(t),x\rangle^2
\qquad\text{and}\qquad
S_i(t):=\{x:i=\min\argmin_j\ell_j(t,x)\}.
\]
The sets $S_i(t)$ form the nearest-axis partition under the lexicographic tie-breaking rule.

\begin{lemma}[First variation of the population risk]
\label{lemma:first-variation-WFA1}
Let $\gamma$ be a $C^1$ curve in $\orth(d)$ and let $\mu\in\Pc_2(\RR^d)$. Suppose that nonzero active ties have $\mu$-measure zero at $\gamma(t)$ for every $t$. Then $t\mapsto\Ec(\gamma(t))$ is differentiable, with
\begin{equation}\label{eq:first-variation-exact}
\frac{d}{dt}\Ec(\gamma(t))
=-2\sum_{i=1}^d\int_{S_i(t)}
\langle u_i(t),x\rangle\langle\dot u_i(t),x\rangle\,d\mu(x).
\end{equation}
\end{lemma}

\begin{proof}
By hypothesis, the minimizing branch is unique for $\mu$-almost every nonzero $x$ at each fixed $t$. At $x=0$, all branches and their derivatives vanish. On a compact time interval, the branch difference quotients are bounded by $C\|x\|_2^2$, so dominated convergence differentiates the finite minimum and gives~\eqref{eq:first-variation-exact}.
\end{proof}

\begin{lemma}[Second variation across active walls]
\label{lemma:second-variation-WFA1}
Let $\gamma$ be $C^2$, let $\mu\in\Pc_2(\RR^d)$ have density $d\mu(x)=p(x)\,dx$ with $p\in C^1(\RR^d)$, and assume
\[
M_{\partial,4}(p):=
\sup_{\|n\|_2=1}\int_{n^\perp}(1+\|x\|_2^4)p(x)\,d\mathcal H^{d-1}(x)<\infty.
\]
For $i<j$ and $\sigma\in\{-1,1\}$, define the regular active wall
\[
\begin{aligned}
\Sigma_{ij}^{\sigma}(t):=\bigl\{x:\;&
\langle u_i(t),x\rangle=\sigma\langle u_j(t),x\rangle,
\quad |\langle u_j(t),x\rangle|>0,\\
&|\langle u_i(t),x\rangle|>|\langle u_k(t),x\rangle|
\ \text{for }k\notin\{i,j\}\bigr\}.
\end{aligned}
\]
Then $U\mapsto\Ec(U)$ belongs to $C^2(\orth(d))$. In particular, $t\mapsto\Ec(\gamma(t))$ is $C^2$ and
\begin{align}
\frac{d^2}{dt^2}\Ec(\gamma(t))
&=-\sqrt2\sum_{i<j}\sum_{\sigma\in\{-1,1\}}
\int_{\Sigma_{ij}^{\sigma}(t)}
|\langle u_j(t),x\rangle|
\langle\dot u_j(t)-\sigma\dot u_i(t),x\rangle^2
p(x)\,d\mathcal H^{d-1}(x)\label{2nd-var-1}\\
&\quad-2\sum_{i=1}^d\int_{S_i(t)}
\left(\langle\dot u_i(t),x\rangle^2
+\langle u_i(t),x\rangle\langle\ddot u_i(t),x\rangle\right)p(x)\,dx.\label{2nd-var-2}
\end{align}
\end{lemma}

\begin{proof}
First multiply the integrand in~\eqref{eq:first-variation-exact} by a smooth radial cutoff and differentiate on each cell. Away from intersections of tie hyperplanes, the implicit-function theorem provides a smooth local flow of the active walls. Such intersections, including triple-active points, have codimension at least two and hence zero $\mathcal H^{d-1}$ measure. For $g_{ij}:=\ell_i-\ell_j$, the outward unit normal to $S_i$ along a regular wall is $\nabla_xg_{ij}/\|\nabla_xg_{ij}\|$, and the corresponding normal velocity is $-\partial_tg_{ij}/\|\nabla_xg_{ij}\|$. Since the two active branches agree on their common wall, their first-order boundary contributions cancel. Pairing the boundary contributions when differentiating the
first-variation formula instead gives
\[
-\int_{\{g_{ij}=0\}_{\mathrm{active}}} \frac{(\partial_tg_{ij})^2}{\|\nabla_xg_{ij}\|_2} p\,d\mathcal H^{d-1}.
\]
On the signed tie hyperplane $\langle u_i,x\rangle=\sigma\langle u_j,x\rangle$, we have $\|\nabla_xg_{ij}\|_2 =2\sqrt2|\langle u_j,x\rangle|$ and $\partial_tg_{ij} =2\langle u_j,x\rangle \langle\dot u_j-\sigma\dot u_i,x\rangle$. Substitution gives the wall term in~\eqref{2nd-var-1}; the quotient extends continuously by zero where both projections vanish, a codimension-two set. Differentiating the integrand within each cell gives the bulk term~\eqref{2nd-var-2}.

To remove the cutoff, choose smooth local orthogonal parametrizations of the finitely many moving tie hyperplanes. Under these parametrizations, the active-wall indicators converge almost everywhere, since their discontinuities occur only on the codimension-two strata. Dominated convergence therefore gives continuity of the truncated wall integrals, while the bulk integrals are continuous by the same theorem, using the integrable bound $C\|x\|_2^2p(x)$.

The truncation can be removed uniformly in local coordinates. Indeed, the bulk tails vanish by the finite second moment of $\mu$, whereas the wall tails are bounded uniformly by $CM_{\partial,4}(p)/R$, since $\|x\|_2^3\one_{\{\|x\|_2>R\}} \leq\|x\|_2^4/R$. The truncated objectives and their first and second derivatives thus converge locally uniformly to the corresponding expressions without the cutoff. Since the truncated objectives are $C^2$, their limit is $C^2$ as well, with second variation given by~\eqref{2nd-var-1}--\eqref{2nd-var-2}.
\end{proof}

\begin{corollary}[Gaussian-convolution applicability]
\label{cor:gaussian-variation-applicability}
If $X=Y+\tau Z$, where $Z\sim\Nc(0,I_d)$ is independent of $Y$, $\tau>0$, and $\EE\|Y\|_2^4<\infty$, then the density $p_\tau$ of $X$ satisfies the hypotheses of Lemmas~\ref{lemma:first-variation-WFA1} and~\ref{lemma:second-variation-WFA1}.
\end{corollary}

\begin{proof}
Gaussian convolution makes $p_\tau$ smooth. For $H=n^\perp$, let $G_H$ be a standard Gaussian vector in $H$. Conditioning on $Y=y$ gives
\[
\int_H(1+\|x\|_2^4)\phi_\tau(x-y)\,d\mathcal H^{d-1}(x)
=\phi_{\tau,1}(n^\top y)\,\EE\!\left[1+\|P_Hy+\tau G_H\|_2^4\right]
\leq C_{d,\tau}(1+\|y\|_2^4).
\]
The bound is uniform in $n$. Averaging over $Y$ yields $M_{\partial,4}(p_\tau)<\infty$, while the finite second moment and absolute continuity are immediate.
\end{proof}

The second-variation formula separates the contributions from the active cells and their moving boundaries. In the noiseless model, the nearest-axis assignment is fixed near $\Ustar$, allowing an exact curvature calculation that will serve as the baseline for the positive-noise analysis.

\begin{lemma}[Exact noiseless curvature]\label{axis-case-exact}
For the signal law $\nu$ of Example~\ref{ex:gaussian-convolution}, the Hessian at $\Ustar$ satisfies, for every $\xi=\Ustar\Omega\in T_{\Ustar}\orth(d)$ with $\Omega^\top=-\Omega$,
\[
2a_{\min}\|\xi\|_{\Fr}^2
\leq\Hess\,\Ec_\nu(\Ustar)[\xi,\xi]
=2\sum_{j<k}(a_j+a_k)\Omega_{jk}^2
\leq2a_{\max}\|\xi\|_{\Fr}^2.
\]
Moreover, there exists $r_0\in(0,1/8]$, depending only on $d$ and $\nu$, such that, for every $0<r\leq r_0$, the ball $\Nc=B_r(\Ustar)$ satisfies the geometric conditions of Assumption~\ref{ass:strong-convexity}, $\dframe(U,V)=d_g(U,V)$ for all $U,V\in\Nc$, $c_{\Nc}>1/\sqrt2$, and the translates $\Nc R$, $R\in\Spm$, are pairwise disjoint. The noiseless risk is smooth on a neighborhood of $\overline{\Nc}$ and satisfies
\[
a_{\min}I\preceq\Hess\,\Ec_\nu(U)\preceq3a_{\max}I,
\qquad \text{for }\,U\in\Nc.
\]
Consequently, $\Ec_\nu$ is geodesically $a_{\min}$-strongly convex on each $\Nc R$.
\end{lemma}

\begin{proof}
For $U$ sufficiently close to $\Ustar$, we have $|\langle u_j,\Ustar_j\rangle|>|\langle u_k,\Ustar_j\rangle|$ whenever $k\neq j$. Thus the nearest candidate axis to every nonzero point of $\Span\{\Ustar_j\}$ is $\Span\{u_j\}$, giving
\[
\Ec_\nu(U)=\sum_{j=1}^d a_j(1-\langle u_j,\Ustar_j\rangle^2).
\]
Along $\gamma(t)=\Ustar\exp(t\Omega)$, we have
\[
\bigl(\exp(t\Omega)\bigr)_{jj}
=1-\frac{t^2}{2}\sum_{k\neq j}\Omega_{jk}^2+o(t^2).
\]
Twice differentiating the risk identity yields the stated Hessian formula, with the two-sided bounds following from $\|\xi\|_{\Fr}^2=2\sum_{j<k}\Omega_{jk}^2$ and $2a_{\min}\leq a_j+a_k\leq2a_{\max}$.

We now extend these curvature bounds to a neighborhood of $\Ustar$. The fixed-branch expression is smooth, so continuity of its Hessian gives a sufficiently small radius $r_0\in(0,1/8]$ such that every ball $\Nc=B_r(\Ustar)$, $0<r\leq r_0$, has its closure within the fixed-branch region and satisfies $a_{\min}I\preceq\Hess\,\Ec_\nu\preceq3a_{\max}I$. Since $\orth(d)$ has positive convexity radius, we may further decrease $r_0$ so that these balls are geodesically convex, with unique minimizing geodesics and smoothly varying Riemannian logarithms.

It remains to verify that the same neighborhood satisfies the frame-alignment and switching margin conditions. For $U,V\in\Nc$ and $R\in\Spm\setminus\{I\}$, the separation of distinct signed permutations gives
\[
\|UR-V\|_{\Fr}
\geq\|R-I\|_{\Fr}-\|U-V\|_{\Fr}
\geq2-2r>2r>d_g(U,V).
\]
Thus no nontrivial signed permutation decreases the geodesic distance, proving $\dframe(U,V)=d_g(U,V)$ and the pairwise disjointness of the translates $\Nc R$. 
The restriction $r_0\leq1/8$ is chosen to ensure the geometric margin condition in Assumption~\ref{ass:switching-margin}. Indeed, since $r_{\Nc}=r\leq1/8$ and $\cos(2r)-2r$ is decreasing in $r$,
\[
c_{\Nc}=\cos(2r)-2r
\geq\cos(1/4)-1/4
>\frac1{\sqrt2}.
\]
Together with geodesic convexity, the uniform Hessian lower bound gives $a_{\min}$-strong convexity on $\Nc$, which extends to every $\Nc R$ by signed-permutation invariance.
\end{proof}

\subsection{Hessian stability and fixed-dimensional local regularity}

The preceding subsection established local strong convexity of the noiseless risk $\Ec_\nu$. Gaussian convolution allows samples to switch their nearest axis, so the corresponding curvature bounds do not follow directly. Throughout this subsection, we use the axial-mixture construction of Example~\ref{ex:gaussian-convolution}, requiring only finite fourth moments of the component laws in place of their sub-exponential moment condition. The next theorem is proved after Lemma~\ref{lem:hessian-stability} below.

\begin{theorem}\label{thm:strong-cvx}
    Fix $0<r\leq r_0$ as in Lemma~\ref{axis-case-exact} and set $\Nc=B_r(\Ustar)$. Then, for every $0<\tau<\tau_1$, the population risk $\Ec$ belongs to $C^2(\Norbit)$ and satisfies
    \[
    mI\preceq\Hess\,\Ec(U)\preceq LI,
    \qquad \text{for}\quad U\in\Norbit,
    \]
    with $m=a_{\min}/2$ and $L=4a_{\max}$. Consequently, $\Ec$ is geodesically $m$-strongly convex on each component $\Nc R$, $R\in\Spm$, and its Riemannian gradient is $L$-Lipschitz in the sense of Assumption~\ref{ass:strong-convexity}.
\end{theorem}

\begin{lemma}[Hessian stability under convolution]
\label{lem:hessian-stability}
Assume the components of the axial-mixture law $\nu$ have densities bounded by a common constant $M<\infty$ and finite fourth moments. For each $\tau>0$, let $\mu$ be the law of $Y+\tau Z$. Fix an open ball $\Nc_0=B_r(\Ustar)$ with $0<r\leq1/2$ whose closure lies in the aligned fixed-branch neighborhood. There exists $C=C(d,M)<\infty$, independent of $\tau$, such that, for every $\tau>0$,
\[
\sup_{U\in\Nc_0}
\bigl\|\Hess\,\Ec_\mu(U)-\Hess\,\Ec_\nu(U)\bigr\|_{\opp}
\leq C\tau^3.
\]
Consequently, if $\Hess\,\Ec_\nu\succeq2mI$ on $\Nc_0$ for some $m>0$ and $\tau\leq(m/C)^{1/3}$, then
\[
mI\preceq\Hess\,\Ec_\mu(U)
\preceq\Hess\,\Ec_\nu(U)+mI
\qquad\text{for every}\quad U\in\Nc_0.
\]
\end{lemma}

\begin{proof}[Proof of Lemma~\ref{lem:hessian-stability}]
We first establish a uniform margin on the true axes and use it to localize changes in nearest-axis assignments. We then reduce the Hessian comparison to two contributions, corresponding to the cellwise second derivatives and the moving cell boundaries, and bound each uniformly by $C\tau^3$. The curvature conclusion will follow from the resulting operator-norm estimate.

\emph{Uniform margin and switching localization.} Fix $U=[u_1,\ldots,u_d]\in\Nc_0$. Since $\|u_j-\ustar_j\|\leq\|U-\Ustar\|_{\Fr}\leq d_g(U,\Ustar)<r$, we have
\[
u_j^\top\ustar_j=1-\frac{\|u_j-\ustar_j\|^2}{2}\geq1-\frac{r^2}{2}.
\]
Moreover, orthogonality of the columns of $\Ustar$ gives $|u_\ell^\top\ustar_j|=|(u_\ell-\ustar_\ell)^\top\ustar_j|\leq r$ for $\ell\neq j$. Consequently,
\[
(u_j^\top\ustar_j)^2-(u_\ell^\top\ustar_j)^2\geq\Bigl(1-\frac{r^2}{2}\Bigr)^2-r^2\geq\frac12,
\]
where the last inequality holds for $r\leq1/2$.
For $y=s\ustar_j$, the identity
\[
\dist^2(y,\Span(u))=\|y\|^2-(u^\top y)^2=s^2\bigl(1-(u^\top\ustar_j)^2\bigr)
\]
therefore yields
\[
\dist^2(y,\Span(u_\ell))-\dist^2(y,\Span(u_j))\geq\frac{s^2}{2},
\]
uniformly over $U\in\Nc_0$, $j$, and $\ell\neq j$.

We next use this margin to show that switching can occur only near the origin. Condition on the mixture index $J=j$, write $Y=s\ustar_j$ and $X=Y+\tau Z$, and let $i_U(x)$ denote the lexicographic cell index. Since the axial laws have densities, $s=0$ is a null event, and the preceding margin gives $i_U(Y)=j$ almost surely. Define $A_U:=\{i_U(X)\neq i_U(Y)\}$ and, for $\ell\neq j$,
\[
\Delta_{j\ell,U}(x):=\dist^2(x,\Span(u_\ell))-\dist^2(x,\Span(u_j))=x^\top A_{j\ell,U}x,
\]
where $A_{j\ell,U}:=u_ju_j^\top-u_\ell u_\ell^\top$. Since $\|A_{j\ell,U}\|_{\opp}\leq2$,
\[
|\Delta_{j\ell,U}(x)-\Delta_{j\ell,U}(y)|\leq2(\|x\|+\|y\|)\|x-y\|.
\]
On $A_U$, some $\ell\neq j$ satisfies $\Delta_{j\ell,U}(X)\leq0$, whereas the margin gives $\Delta_{j\ell,U}(Y)\geq s^2/2$. Combining these inequalities with $\|X-Y\|=\tau\|Z\|$ and $\|X\|\leq|s|+\tau\|Z\|$ yields
\[
\frac{s^2}{2}\leq2(\|X\|+\|Y\|)\|X-Y\|\leq2(2|s|+\tau\|Z\|)\tau\|Z\|\leq4(|s|+\tau\|Z\|)\tau\|Z\|.
\]
Choose $\Csw:=10$, for which $4(\Csw^{-1}+\Csw^{-2})=0.44<\tfrac12$. If $|s|>\Csw\tau\|Z\|$, the preceding bound would imply
\[
\frac{s^2}{2}\leq4(|s|+\tau\|Z\|)\tau\|Z\|<4(\Csw^{-1}+\Csw^{-2})s^2<\frac{s^2}{2},
\]
a contradiction. Thus, uniformly over $U\in\Nc_0$ and the mixture components,
\(
A_U\subseteq\{|s|\leq\Csw\tau\|Z\|\}.
\)

\emph{Reduction to quadratic forms.} Fix $\xi\in T_U\orth(d)$ with $\|\xi\|_{\Fr}=1$ and let $\gamma(t)=U\exp(tU^\top\xi)$ be the geodesic starting at $U$ with velocity $\xi$. Its second-order expansion gives
\[
\gamma(t)=U+t\xi-\frac{t^2}{2}U\xi^\top\xi+\bigO(t^3),
\qquad\text{so}\quad
\ddot\gamma(0)=-U\xi^\top\xi.
\]
Write
\[
Q_\mu(U,\xi):=\left.\frac{d^2}{dt^2} \Ec_\mu(\gamma(t))\right|_{t=0}
\qquad\text{and}\qquad
Q_\nu(U,\xi):=\left.\frac{d^2}{dt^2}\Ec_\nu(\gamma(t))\right|_{t=0}.
\]
These are the quadratic forms of the corresponding Riemannian Hessians. Since their difference is self-adjoint, it suffices to establish
\[
\sup_{U\in\Nc_0}\sup_{\|\xi\|_{\Fr}=1}|Q_\mu(U,\xi)-Q_\nu(U,\xi)|\leq C\tau^3.
\]
To this end, let $\Sigma_{ij}^{\sigma}(U):=\Sigma_{ij}^{\sigma}(0)$ and decompose the second-variation formula of Lemma~\ref{lemma:second-variation-WFA1}, which applies to $p_\mu$ by Corollary~\ref{cor:gaussian-variation-applicability}, into the wall term
\[
\mathcal B_\mu(U,\xi):=-\sqrt2\sum_{i<j}\sum_{\sigma\in\{-1,1\}}\int_{\Sigma_{ij}^{\sigma}(U)}|\langle u_j,x\rangle|\langle\xi_j-\sigma\xi_i,x\rangle^2p_\mu(x)\,d\mathcal H^{d-1}(x)
\]
and the bulk term
\[
\mathcal R_\mu(U,\xi):=-2\sum_{i=1}^d\int_{S_i(U)}
\left(\langle\xi_i,x\rangle^2+\langle u_i,x\rangle\langle\ddot\gamma(0)_i,x\rangle\right)d\mu(x).
\]
Thus $Q_\mu=\mathcal B_\mu+\mathcal R_\mu$. For the singular noiseless law, however, Lemma~\ref{axis-case-exact} ensures that the active branch is fixed on each signal axis, so $Q_\nu=\mathcal R_\nu$. It remains to control the bulk difference $\mathcal R_\mu-\mathcal R_\nu$ and the noisy wall term $\mathcal B_\mu$.

\emph{Bulk contribution.}
Set $a:=\ddot\gamma(0)$. The preceding expansion gives $\|a\|_{\Fr}\leq\|\xi\|_{\Fr}^2=1$. For each $i$, define
\[
q_{U,\xi,i}(x):=-2\left(\langle\xi_i,x\rangle^2 +\langle u_i,x\rangle\langle a_i,x\rangle \right),
\]
and set
$q_{U,\xi}(x) :=\sum_{i=1}^d q_{U,\xi,i}(x) \one_{\{x\in S_i(U)\}}$.
Then
\[
\mathcal R_\mu(U,\xi)=\EE[q_{U,\xi}(X)]\qquad\text{and}\qquad \mathcal R_\nu(U,\xi)=\EE[q_{U,\xi}(Y)].
\]

We first bound the weighted contribution of the switching event. By the localization established above and the density bound $M$ on each axial law, conditioning on $Z$ and $J=j$ gives
\[
\PP(A_U\mid Z,J=j)\leq2M\Csw\tau\|Z\|
\]
and
\[
\EE[s^2\one_{A_U}\mid Z,J=j] \leq\frac{2M}{3}\Csw^3\tau^3\|Z\|^3.
\]
Since $\|X\|^2\leq2s^2+2\tau^2\|Z\|^2$, averaging over $Z$ and $J$ yields
\begin{equation}\label{eq:switching-weighted-moment}
\sup_{U\in\Nc_0}\EE\!\left[\|X\|^2\one_{A_U}\right]\leq 4M\left(\frac{\Csw^3}{3}+\Csw\right) \EE\|Z\|^3\,\tau^3.
\end{equation}

To exploit this bound, we first cancel the contribution obtained by evaluating the same branch at $X$ and $Y$. Indeed, differentiating the unit-column constraint twice gives $\langle u_J,a_J\rangle=-\|\xi_J\|^2$. Gaussian isotropy therefore implies
\[
\EE\!\left[ q_{U,\xi,J}(X)-q_{U,\xi,J}(Y) \mid Y,J \right] =-2\tau^2 \bigl(\|\xi_J\|^2+\langle u_J,a_J\rangle\bigr) =0.
\]
Adding and subtracting this fixed-branch contribution consequently gives
\[
\mathcal R_\mu(U,\xi)-\mathcal R_\nu(U,\xi) = \EE\!\left[
\bigl( q_{U,\xi,i_U(X)}(X)-q_{U,\xi,J}(X) \bigr)\one_{A_U} \right].
\]
Moreover, $\|\xi_i\|\leq 1$ and $\|a_i\|\leq 1$ imply $|q_{U,\xi,i}(x)|\leq4\|x\|^2$, so the integrand is bounded in absolute value by $8\|X\|^2\one_{A_U}$. Applying \eqref{eq:switching-weighted-moment} yields
\begin{equation}\label{eq:bulk-hessian-stability}
\sup_{U\in\Nc_0} \sup_{\|\xi\|_{\Fr}=1} |\mathcal R_\mu(U,\xi)-\mathcal R_\nu(U,\xi)| \leq C_{\mathrm{bulk}}\tau^3,
\end{equation}
where
\[ 
C_{\mathrm{bulk}} :=32M\left(\frac{\Csw^3}{3}+\Csw\right) \EE\|Z\|^3.
\]

\emph{Wall contribution.}
The signed-wall formula, together with $(a+b)^2\leq 2(a^2+b^2)$, gives
\[
|\mathcal B_\mu(U,\xi)| \leq4\sqrt2(d-1)\|\xi\|_{\Fr}^2 \max_{i<j,\,\sigma} \int_{\Sigma_{ij}^{\sigma}(U)} \|x\|_2^3p_\mu(x)\,d\mathcal H^{d-1}(x).
\]
Thus it suffices to show that
\[
\sup_{U\in\Nc_0}\max_{i<j,\,\sigma} \int_{\Sigma_{ij}^{\sigma}(U)} \|x\|_2^3p_\mu(x)\,d\mathcal H^{d-1}(x) \leq C_6\tau^3,
\]
where $C_6$ depends only on $d$ and $M$.

We first show that every nonzero point on a true axis has a neighborhood free of cell
boundaries. Fix $U\in\Nc_0$, set $\vartheta:=1/20$, and let $y=s\ustar_j$ with $s\neq0$. For $\ell\neq j$, define
\[
G_{j\ell,U}(x) :=\dist^2(x,\Span(u_\ell)) -\dist^2(x,\Span(u_j)).
\]
The margin estimate gives $G_{j\ell,U}(y)\geq s^2/2$, whereas
\[
|G_{j\ell,U}(x)-G_{j\ell,U}(y)| \leq2(\|x\|+\|y\|)\|x-y\|.
\]
Whenever $\|x-y\|\leq\vartheta|s|$, we have $\|x\|\leq(1+\vartheta)|s|$, and hence
\[
|G_{j\ell,U}(x)-G_{j\ell,U}(y)| \leq2(2+\vartheta)\vartheta s^2 =0.205s^2 \leq\frac{s^2}{4}.
\]
Consequently, $G_{j\ell,U}(x)\geq s^2/4>0$ throughout this ball, so
\[
B(y,\vartheta|s|) \cap\bigcup_{i=1}^d\partial S_i(U) =\varnothing.
\]

We now use this separation to control the Gaussian mass on the cell boundaries. The convolution density satisfies
\[
p_\mu(x) =\sum_{j=1}^d p_j \int_{\RR}\phi_\tau(x-s\ustar_j)\, d\lambda_j(s).
\]
For each $i,j$, define
\[
K_{ij}(U,s,\tau) :=\int_{\partial S_i(U)} \|x\|^3\phi_\tau(x-s\ustar_j)\,d\sigma(x),
\]
so that
\[
\int_{\partial S_i(U)} \|x\|^3p_\mu(x)\,d\sigma(x) =\sum_{j=1}^d p_j \int_{\RR}K_{ij}(U,s,\tau)\,d\lambda_j(s).
\]
The separation established above gives
\[
\partial S_i(U) \subseteq \{x:\|x-s\ustar_j\|\geq\vartheta|s|\},
\]
which also holds trivially for $s=0$. Moreover, each cell boundary lies in a finite union of signed bisector hyperplanes
\[
H_{i\ell}^{\pm}(U) :=\{x:\langle x,u_i\pm u_\ell\rangle=0\}, \quad\text{for}\quad \ell\neq i.
\]
It therefore suffices to bound integrals of the form
\[
I_H(y,\tau,r) :=\int_{H\cap\{\|x-y\|\geq r\}} \|x\|^3\phi_\tau(x-y)\,d\sigma(x)
\]
for an origin hyperplane $H$, with $y=s\ustar_j$ and $r=\vartheta|s|$.

We claim that
\[
I_H(y,\tau,r) \leq C_8\tau^{-1} (\|y\|^3+\tau^3)e^{-r^2/(4\tau^2)},
\]
uniformly over such hyperplanes. Indeed, on $\{\|x-y\|\geq r\}$,
\[
e^{-\|x-y\|^2/(2\tau^2)} \leq e^{-r^2/(4\tau^2)} e^{-\|x-y\|^2/(4\tau^2)}.
\]
Using the Gaussian density formula gives
\[
I_H(y,\tau,r) \leq C\tau^{-d}e^{-r^2/(4\tau^2)} \int_H\|x\|^3 e^{-\|x-y\|^2/(4\tau^2)}\,d\sigma(x).
\]
To estimate the remaining integral, write $y=y_H+y_\perp$, with $y_H\in H$ and $y_\perp\perp H$, and parametrize $H$ by $x=y_H+2\tau z$, with $z\in H$.
Then
\[
\|x-y\|^2 =4\tau^2\|z\|^2+\|y_\perp\|^2\qquad\text{and} \qquad d\sigma(x)=(2\tau)^{d-1}dz.
\]
Since
$\|x\|^3\leq C(\|y\|^3+\tau^3\|z\|^3)$, the Gaussian moment integral yields
\[
\int_H\|x\|^3 e^{-\|x-y\|^2/(4\tau^2)}\,d\sigma(x)\leq C\tau^{d-1}(\|y\|^3+\tau^3) \int_H(1+\|z\|^3)e^{-\|z\|^2}\,dz\leq C\tau^{d-1}(\|y\|^3+\tau^3),
\]
proving the claimed hyperplane estimate.

Applying this bound to the finitely many hyperplanes covering $\partial S_i(U)$, with $y=s\ustar_j$ and $r=\vartheta|s|$, gives
\[
K_{ij}(U,s,\tau) \leq C_9\tau^{-1} (|s|^3+\tau^3)e^{-cs^2/\tau^2}, \qquad\text{with}\quad c:=\frac{\vartheta^2}{4},
\]
where $C_9$ depends only on $d$. Consequently,
\[
\int_{\partial S_i(U)} \|x\|^3p_\mu(x)\,d\sigma(x)\leq C_9\sum_{j=1}^d p_j \int_{\RR} \tau^{-1}(|s|^3+\tau^3) e^{-cs^2/\tau^2}\,d\lambda_j(s).
\]
Since each $\lambda_j$ has density bounded by $M$, the change of variables $s=\tau r$
gives
\[
\int_{\RR} \tau^{-1}(|s|^3+\tau^3) e^{-cs^2/\tau^2}\,d\lambda_j(s)\leq M\tau^3\int_{\RR} (|r|^3+1)e^{-cr^2}\,dr \leq C\tau^3.
\]
Using $\sum_jp_j=1$, we conclude that
\[
\sup_{U\in\Nc_0} \int_{\partial S_i(U)} \|x\|^3p_\mu(x)\,d\sigma(x) \leq C_6\tau^3.
\]
Every active wall $\Sigma_{ij}^{\sigma}(U)$ is contained in $\partial S_i(U)$ up to an $\mathcal H^{d-1}$-null stratum, so this establishes the required wall estimate. In particular,
\[
\sup_{U\in\Nc_0,\,\|\xi\|_{\Fr}=1} |\mathcal B_\mu(U,\xi)| \leq C_{\mathrm{wall}}\tau^3, \quad\text{where}\quad C_{\mathrm{wall}}:=4\sqrt2(d-1)C_6.
\]
Combining this estimate with \eqref{eq:bulk-hessian-stability} yields
\[
\sup_{U\in\Nc_0} \sup_{\|\xi\|_{\Fr}=1} |Q_\mu(U,\xi)-Q_\nu(U,\xi)| \leq C\tau^3 \quad\text{with}\quad C:=C_{\mathrm{bulk}}+C_{\mathrm{wall}}.
\]
Since the Hessian difference is self-adjoint, this proves
\[
\sup_{U\in\Nc_0} \bigl\| \Hess\,\Ec_\mu(U)-\Hess\,\Ec_\nu(U) \bigr\|_{\opp} \leq C\tau^3,
\]
where $C$ depends only on $d$ and $M$, and the estimate holds for every $\tau>0$. Finally, if $\Hess\,\Ec_\nu\succeq2mI$ on $\Nc_0$ and $\tau\leq(m/C)^{1/3}$, then the operator-norm bound gives
\[
\Hess\,\Ec_\nu(U)-mI \preceq\Hess\,\Ec_\mu(U) \preceq\Hess\,\Ec_\nu(U)+mI \qquad\text{and}\qquad \Hess\,\Ec_\mu(U)\succeq mI,
\]
for every $U\in\Nc_0$, completing the proof.
\end{proof}

\begin{proof}[Proof of Theorem~\ref{thm:strong-cvx}]
Fix $r\leq r_0$ and set $\Nc=B_r(\Ustar)$, so that $r\leq1/8$. By Lemma~\ref{axis-case-exact}, $\Hess\,\Ec_\nu\succeq a_{\min}I=2mI$ on $\Nc$, and $\Nc$ satisfies the neighborhood and alignment conditions in Assumption~\ref{ass:strong-convexity}. Lemma~\ref{lem:hessian-stability} therefore gives, for $\tau<\tau_1$ and every $U\in\Nc$,
\[
mI\preceq\Hess\,\Ec_\mu(U)\preceq\Hess\,\Ec_\nu(U)+mI\preceq(3a_{\max}+m)I\preceq4a_{\max}I,
\]
where the last inequality uses $m=a_{\min}/2\leq a_{\max}/2$. Moreover, Corollary~\ref{cor:gaussian-variation-applicability} and Lemma~\ref{lemma:second-variation-WFA1} give $\Ec_\mu\in C^2(\orth(d))$. The lower bound yields geodesic $m$-strong convexity on the geodesically convex ball $\Nc$. For the gradient, differentiating $t\mapsto\mathsf P_{\gamma(t)\to U}\Rgrad\Ec(\gamma(t))$ along the minimizing geodesic $\gamma$ from $V$ to $U$ gives
\[
\bigl\|\Rgrad\Ec(U)-\mathsf P_{V\to U}\Rgrad\Ec(V)\bigr\|_{\Fr}
\leq\sup_{\Nc}\|\Hess\,\Ec\|_{\opp}\,d_g(U,V)
\leq4a_{\max}\,d_g(U,V).
\]
Equivariance transports these properties to every $\Nc R$, which proves the theorem with $m=a_{\min}/2$ and $L=4a_{\max}$.
\end{proof}

Local strong convexity alone does not identify $\Ustar$ as a minimizer. The symmetry of the Gaussian convolution model also makes $\Ustar$ a critical point, giving the following result.

\begin{corollary}[Criticality and local minimality]
\label{thm:crit-pt}
For every $\tau>0$, all frames in $\Ustar\Spm$ are critical points of $\Ec$. Moreover, if $\tau<\tau_1$, each $\Ustar R$ is the unique minimizer of $\Ec$ on $\Nc R$, and hence a strict local minimizer on $\orth(d)$.
\end{corollary}

\begin{proof}
    It suffices to prove criticality at $\Ustar$. Rotate coordinates so that $\Ustar=I$, fix $\tau>0$ and $\xi\in T_I\orth(d)$, and let $i^\star(x)$ denote the lexicographic nearest-axis index. The first-variation formula gives
    \[
        D\Ec(I)[\xi]
        =-2\sum_{i=1}^d\EE\!\left[
        \one\{i^\star(X)=i\}X_i\langle\xi_i,X\rangle\right].
    \]
    Condition on the latent axial index $J=j$ and amplitude, so $X=s e_j+\tau Z$. If $i=j$, reflect $Z$ by fixing coordinate $j$ and reversing all other coordinates. This preserves the Gaussian law, every squared coordinate score, and hence $i^\star(X)$; it fixes $X_j$ and reverses $\langle\xi_j,X\rangle$ because $\xi_j\perp e_j$. If $i\neq j$, reverse coordinate $i$ only. This again preserves the law and the cell index, reverses $X_i$, and fixes $\langle\xi_i,X\rangle$ because $\xi_i\perp e_i$. Therefore, every conditional expectation in the display is zero, and, $\Ustar$ is critical for $\tau>0$. Signed-permutation equivariance gives $\Rgrad\Ec(\Ustar R)=0$ for every $R\in\Spm$. When $\tau<\tau_1$, Theorem~\ref{thm:strong-cvx} gives geodesic strong convexity on each $\Nc R$. Consequently, $\Ustar R$ is its unique minimizer and a strict local minimizer on $\orth(d)$.
\end{proof}

\subsection{Dimension-uniform curvature for balanced Gaussian mixtures}
The preceding subsection establishes local regularity for general axial laws at each fixed dimension. However, the constant $C(d,M)$ in Lemma~\ref{lem:hessian-stability} is not controlled uniformly in $d$. The following proposition obtains the sharper curvature estimates needed for balanced Gaussian mixtures by controlling the switching contribution directly. The corresponding global gap is established in the next subsection.

\begin{proposition}[Dimension-uniform Gaussian-mixture curvature]
\label{prop:balanced-gaussian-curvature}
Fix $0<\sigma_-^2\leq\sigma_+^2<\infty$ and $b_+>0$. Specialize Example~\ref{ex:gaussian-convolution} to $p_j=1/d$ and $S_j\sim\Nc(0,\sigma_j^2)$, where $\sigma_-^2\leq\sigma_j^2\leq\sigma_+^2$. There exist $r_0>0$, depending only on $\sigma_-^2$ and $\sigma_+^2$, and $d_0\geq4$, depending additionally on $b_+$, such that, for every $d\geq d_0$ and $0<\tau^2\leq b_+/d$, the following holds. For every $U\in B_{r_0}(\Ustar)$ and $\xi\in T_U\orth(d)$,
\begin{equation}\label{eq:dimension-uniform-hessian}
\frac{\sigma_-^2}{2d}\|\xi\|_{\Fr}^2
\leq\Hess\,\Ec(U)[\xi,\xi]
\leq\frac{4\sigma_+^2}{d}\|\xi\|_{\Fr}^2.
\end{equation}
Moreover,
\[
\frac{\sigma_-^2}{\log2}
\leq K:=\bigl\|\,\|X\|^2\bigr\|_{\psi_1}
\leq\frac{16}{3}\sigma_+^2+4b_+.
\]

The ball $B_{r_0}(\Ustar)$ is geodesically convex, has pairwise disjoint signed-permutation translates, and satisfies $\dframe=d_g$. Consequently, Assumptions~\ref{ass:subexp-norms} and~\ref{ass:strong-convexity} hold on this ball with $K$ as above, $m=\sigma_-^2/(2d)$, and $L=4\sigma_+^2/d$.
\end{proposition}

\begin{proof}
By left multiplication with $\UstarT$, we may assume $\Ustar=I$. Set $b:=\tau^2$ and, for $U$ near $I$, let $v_j:=U^\top e_j$ denote the $j$th row of $U$. Conditional on $J=j$, the score vector $q:=U^\top X$ satisfies
\[
q\mid\{J=j\}\sim\Nc(0,bI_d+\sigma_j^2v_jv_j^\top).
\]
Define $\psi_j(q):=\max_iq_i^2-q_j^2\geq0$. Since $\EE[q_j^2\mid J=j]=b+\sigma_j^2U_{jj}^2$, the population risk decomposes as
\begin{equation}\label{eq:oracle-switching-decomposition}
\Ec(U)=C_b-\frac1d\sum_{j=1}^d\sigma_j^2U_{jj}^2-\mathcal S(U),
\qquad\text{where}\quad
\mathcal S(U):=\frac1d\sum_{j=1}^d\EE_{v_j}[\psi_j(q)],
\end{equation}
and $C_b$ is independent of $U$. The middle term corresponds to fixing each observation to its true axis, whereas $\mathcal S$ accounts for the gain from switching to another coordinate. We show that the former has curvature of order $d^{-1}$, while the curvature of the switching correction is $\littleo(d^{-1})$.

\emph{Oracle term.}
Along a geodesic $U(t)=U\exp(t\Omega)$ with $\Omega^\top=-\Omega$, the second derivative of the oracle risk is
\begin{equation}\label{eq:oracle-hessian-near-frame}
-\frac2d\sum_{j=1}^d\sigma_j^2 \left((U\Omega)_{jj}^2+U_{jj}(U\Omega^2)_{jj}\right).
\end{equation}
At $U=I$, this equals $\frac2d\sum_{j<k}(\sigma_j^2+\sigma_k^2)\Omega_{jk}^2$, which lies between $(2\sigma_-^2/d)\|\Omega\|_{\Fr}^2$ and $(2\sigma_+^2/d)\|\Omega\|_{\Fr}^2$.
To extend these bounds to a neighborhood of $I$, suppose $\|U-I\|_{\Fr}\leq r\leq1$. Since $\operatorname{diag}(U\Omega) =\operatorname{diag}((U-I)\Omega)$, we have $\|\operatorname{diag}(U\Omega)\|_2 \leq r\|\Omega\|_{\Fr}$.

For the remaining term, write $U_{jj}(U\Omega^2)_{jj}-(\Omega^2)_{jj} =(U_{jj}-1)(U\Omega^2)_{jj} +((U-I)\Omega^2)_{jj}$. Using $\sum_j|a_{jj}|\leq\|a\|_{S^1}$ and $\|\Omega^2\|_{S^1}=\|\Omega\|_{\Fr}^2$ gives
\[
\left|\sum_j\sigma_j^2
\left\{U_{jj}(U\Omega^2)_{jj}-(\Omega^2)_{jj}\right\}\right|
\leq2\sigma_+^2r\|\Omega\|_{\Fr}^2.
\]
The nuclear-norm estimate avoids the factor $\sqrt d$ that would arise from a diagonal Cauchy--Schwarz bound. Together with the bound on $\operatorname{diag}(U\Omega)$, it shows that the oracle Hessian differs from its value at $U=I$ by at most $(6\sigma_+^2r/d)\|\Omega\|_{\Fr}^2$.

Choose $r_0\leq\sigma_-^2/(6\sigma_+^2)$, shrinking it by a dimension-independent amount if necessary so that the ball $\{d_g(U,I)<r_0\}$ is geodesically convex, its signed-permutation translates are pairwise disjoint, and $\dframe=d_g$ on it. The preceding estimate then places the oracle Hessian between $(\sigma_-^2/d)\|\Omega\|_{\Fr}^2$ and $(3\sigma_+^2/d)\|\Omega\|_{\Fr}^2$ throughout this ball.

\emph{Switching term.}
To control the curvature of $\mathcal S$, we differentiate the Gaussian density rather than the nonsmooth function $\psi_j$. For a unit vector $v$, the density of $\Nc(0,bI_d+\sigma_j^2vv^\top)$ is
\[
p_{j,v}(q)=c_{j,b}\exp\!\left( -\frac{\|q\|^2}{2b}+\frac{\delta_j}{2}(v^\top q)^2 \right), \qquad \text{with}\qquad\delta_j:=\frac1b-\frac1{b+\sigma_j^2}\leq\frac1b.
\]
Along $v(t)=\exp(-t\Omega)v$, set $w:=v'(0)=-\Omega v$ and $z:=v''(0)=\Omega^2v$. Differentiating the density twice gives
\begin{equation}\label{eq:rotated-gaussian-density-second-derivative}
\frac{\partial_t^2p_{j,v(t)}(q)|_{t=0}}{p_{j,v}(q)}
=\delta_j\left((w^\top q)^2+(v^\top q)(z^\top q)\right)
+\delta_j^2(v^\top q)^2(w^\top q)^2.
\end{equation}
For fixed $b>0$, the polynomial factor is Gaussian-integrable locally uniformly in $v$, justifying differentiation under the integral.

To estimate the resulting terms, write $q=\alpha v+\varepsilon$, where $\alpha\sim\Nc(0,b+\sigma_j^2)$ and $\varepsilon\sim\Nc(0,b(I_d-vv^\top))$ are independent, and set $\bar\varepsilon:=\max_i|\varepsilon_i|$. For $v=v_j$, the bound $\|U-I\|_{\Fr}\leq r_0$ gives $|(v_j)_j-1|\leq r_0$ and $|(v_j)_i|\leq r_0$ for $i\neq j$. Consequently,
\[
|(v_j)_j|-\max_{i\neq j}|(v_j)_i|\geq1-2r_0\geq\tfrac23,
\]
where the last inequality uses $r_0\leq\sigma_-^2/(6\sigma_+^2)\leq1/6$. If $\psi_j(q)>0$, some $i\neq j$ satisfies $|q_i|>|q_j|$. Comparing these coordinates using $q=\alpha v+\varepsilon$ gives $\tfrac23|\alpha|<2\bar\varepsilon$, and then $\max_i|q_i|\leq|\alpha|+\bar\varepsilon\leq4\bar\varepsilon$. Hence
\begin{equation}\label{eq:small-signal-switching-event}
|\alpha|\leq3\bar\varepsilon\qquad \text{and}\qquad 0\leq\psi_j(q)\leq16\bar\varepsilon^2.
\end{equation}
Thus switching requires the signal amplitude to be comparable to the largest noise coordinate.

We next record the moment bounds used to control this event. For every fixed integer $k\geq1$, a Gaussian union bound followed by tail integration gives
\begin{equation}\label{eq:gaussian-coordinate-maximum-moments}
\EE[\bar\varepsilon^k]\leq C_k\bigl(b\log(ed)\bigr)^{k/2}.
\end{equation}
Moreover, the density of $\alpha$ is bounded by $(2\pi\sigma_-^2)^{-1/2}$, uniformly in $j$, $d$, and $b$. Integrating this bound over $[-t,t]$ gives, for every $t\geq0$,
\begin{equation}\label{eq:truncated-signal-moments}
\PP(|\alpha|\leq t)\leq Ct,
\qquad
\EE[|\alpha|\one_{\{|\alpha|\leq t\}}]\leq Ct^2,
\quad\text{and}\quad
\EE[\alpha^2\one_{\{|\alpha|\leq t\}}]\leq Ct^3.
\end{equation}

Since $w\perp v$, we have $w^\top q=w^\top\varepsilon$. Similarly, $v^\top z=-\|w\|^2$ implies $z^\top q=-\alpha\|w\|^2+z^\top\varepsilon$. We can therefore apply
\eqref{eq:small-signal-switching-event}--\eqref{eq:truncated-signal-moments} to the four terms in \eqref{eq:rotated-gaussian-density-second-derivative}. Conditioning first on $\varepsilon$ and then using Cauchy--Schwarz, together with
$\delta_j\leq1/b$, $\EE(w^\top\varepsilon)^4=3b^2\|w\|^4$, and $\EE(z^\top\varepsilon)^2\leq b\|z\|^2$, gives
\begin{align*}
\delta_j\EE[\psi_j(w^\top\varepsilon)^2]
&\leq C\delta_j\EE[\bar\varepsilon^3(w^\top\varepsilon)^2]
\leq Cb^{3/2}\log(ed)^{3/2}\|w\|^2,\\
\delta_j\|w\|^2\EE[\psi_j\alpha^2]
&\leq C\delta_j\|w\|^2\EE[\bar\varepsilon^5]
\leq Cb^{3/2}\log(ed)^{5/2}\|w\|^2,\\
\delta_j\EE[\psi_j|\alpha|\,|z^\top\varepsilon|]
&\leq C\delta_j\EE[\bar\varepsilon^4|z^\top\varepsilon|]
\leq Cb^{3/2}\log(ed)^{2}\|z\|,\quad\text{and}\\
\delta_j^2\EE[\psi_j\alpha^2(w^\top\varepsilon)^2]
&\leq C\delta_j^2\EE[\bar\varepsilon^5(w^\top\varepsilon)^2]
\leq Cb^{3/2}\log(ed)^{5/2}\|w\|^2.
\end{align*}
For example, the last bound uses $\EE[\alpha^2\one_{\{|\alpha|\leq 3\bar\varepsilon\}}\mid\varepsilon]\leq C\bar\varepsilon^3$, followed by
\[
\EE[\bar\varepsilon^5(w^\top\varepsilon)^2]\leq(\EE\bar\varepsilon^{10})^{1/2}(\EE(w^\top\varepsilon)^4)^{1/2}\leq Cb^{7/2}\log(ed)^{5/2}\|w\|^2.
\]
In the third line, the corresponding conditional bound is
$\EE[|\alpha|\one_{\{|\alpha| \leq3\bar\varepsilon\}}\mid\varepsilon] \leq C\bar\varepsilon^2$. Unlike the other terms, this contribution depends linearly on $\|z\|$, which accounts for the factor $\sqrt d$ in the subsequent summation.

Consequently, writing $s_j(U):=\EE_{v_j}\psi_j(q)$, the four bounds give
\begin{equation}\label{eq:one-component-switching-hessian}
\left|\frac{d^2}{dt^2}s_j(U\exp(t\Omega))\bigg|_{t=0}\right|
\leq Cb^{3/2}\log(ed)^{5/2}
\left(\|\Omega v_j\|^2+\|\Omega^2v_j\|\right).
\end{equation}
Since $v_1,\ldots,v_d$ form an orthonormal basis, we have
\[
\sum_j\|\Omega v_j\|^2=\|\Omega\|_{\Fr}^2 \qquad \text{and}\qquad\sum_j\|\Omega^2v_j\|\leq\sqrt d\,\|\Omega^2\|_{\Fr} \leq\sqrt d\,\|\Omega\|_{\Fr}^2.
\]
Averaging \eqref{eq:one-component-switching-hessian} over the mixture components and using $b\leq b_+/d$ therefore yields
\begin{equation}\label{eq:uniform-switching-hessian}
|\Hess\,\mathcal S(U)[U\Omega,U\Omega]| \leq
C\frac{\log(ed)^{5/2}}{d^2}\|\Omega\|_{\Fr}^2.
\end{equation}

The constant $C$ in this estimate depends only on $\sigma_-^2,\sigma_+^2$, and $b_+$, and is independent of $d$, $\tau$, $U$, and $\Omega$. The correction is $\littleo(d^{-1})\|\Omega\|_{\Fr}^2$, whereas the oracle curvature is of order $d^{-1}\|\Omega\|_{\Fr}^2$. We may therefore choose $d_0$ so that the right-hand side of~\eqref{eq:uniform-switching-hessian} is at most $(\sigma_-^2/2d)\|\Omega\|_{\Fr}^2$ for all $d\geq d_0$. Combining this estimate with the oracle bounds gives
\[
\frac{\sigma_-^2}{2d}\|\Omega\|_{\Fr}^2\leq\Hess\,\Ec(U)[U\Omega,U\Omega]\leq\frac{4\sigma_+^2}{d}\|\Omega\|_{\Fr}^2,
\]
which proves~\eqref{eq:dimension-uniform-hessian}. On the chosen geodesically convex ball, the lower Hessian bound permits $m=\sigma_-^2/(2d)$, while the upper bound gives the intrinsic gradient-Lipschitz constant $L=4\sigma_+^2/d$. Reflection symmetry of the centered Gaussian components also makes $\Ustar$ critical, as in the proof of Corollary~\ref{thm:crit-pt}.

\emph{Norm of $X$.}
The inequality $\|X\|^2\leq2S_J^2+2\tau^2\|Z\|^2$ reduces the required $\psi_1$ bound to the signal and noise terms. The Gaussian moment-generating function gives $\|S_j^2\|_{\psi_1}=\tfrac83\sigma_j^2$, and convexity in the mixture weights yields $\|S_J^2\|_{\psi_1}\leq\tfrac83\sigma_+^2$. For the noise term, the same calculation gives
\[
\|\tau^2\|Z\|^2\|_{\psi_1}=\frac{2\tau^2}{1-2^{-2/d}}.
\]
To bound this expression uniformly in $d$, write
$h(u):=(1-e^{-u})/u=\int_0^1e^{-us}\,ds$. Since $h$ is decreasing and $u=(2/d)\log2$ decreases with $d$, the function $f(d):=d(1-2^{-2/d}) =2\log2\cdot h((2/d)\log2)$ is increasing, with $f(2)=1$. Consequently, $f(d)\geq1$ for $d\geq2$, and
\[
\|\tau^2\|Z\|^2\|_{\psi_1} =\frac{2\tau^2d}{f(d)} \leq2\tau^2d\leq2b_+.
\]
Subadditivity of the $\psi_1$-gauge therefore yields
\[
\bigl\|\,\|X\|^2\bigr\|_{\psi_1} \leq2\|S_J^2\|_{\psi_1}+2\|\tau^2\|Z\|^2\|_{\psi_1} \leq\frac{16}{3}\sigma_+^2+4b_+,
\]
establishing the upper bound for $K$. Conversely, $\EE\|X\|^2\geq\sigma_-^2$, while Jensen's inequality gives
$\EE\|X\|^2 \leq\log2\cdot\bigl\|\,\|X\|^2\bigr\|_{\psi_1}$. Together these inequalities yield $K\geq\sigma_-^2/\log2$, completing the proof.
\end{proof}

\subsection{Dimension-uniform global gap for balanced Gaussian mixtures}

The local curvature bounds of Proposition~\ref{prop:balanced-gaussian-curvature} do not exclude minimizers outside the aligned neighborhoods. To identify $\Ustar\Spm$ as the complete minimizer orbit and verify Assumption~\ref{ass:sublevel-containment}, we need a positive energy gap outside these neighborhoods. The uniform comparison between $\Ec$ and $\Ec_\nu$ used for general axial laws is too coarse under the noise scaling $\tau^2\leq b_+/d$. Instead, we establish the gap by comparing the captured energies row by row, controlling the switching contribution relative to the noiseless gap.

For $v\in S^{d-1}$, define its axis defect by $\delta(v):=1-\|v\|_\infty^2 =\min_{1\leq i\leq d}\dist^2(v,\Span(e_i))$.

\begin{lemma}[Rowwise switching gain for a Gaussian spike]
\label{lem:balanced-row-switching}
Fix $0<\sigma_-\leq\sigma_+<\infty$. For $\sigma\in[\sigma_-,\sigma_+]$, $v\in S^{d-1}$, $G\sim\Nc(0,1)$, and $Z\sim\Nc(0,I_d)$ independent, set $q=\sigma Gv+\tau Z$ and let $k$ satisfy $|v_k|=\|v\|_\infty$. There is $C<\infty$, depending only on $\sigma_-$ and $\sigma_+$, such that, for every $\tau>0$,
\begin{equation}\label{eq:row-switching-gain}
\EE\bigl[\max_iq_i^2-q_k^2\bigr]
\leq
C\bigl(\tau\sqrt{\log(ed)}+\tau^2\log(ed)\bigr)\delta(v)
+C\tau^3\log(ed)^{3/2}.
\end{equation}
\end{lemma}

\begin{proof}
Write $M:=\max_i|Z_i|$. A Gaussian union bound followed by tail integration gives $\EE M^p\leq C_p\log(ed)^{p/2}$ for every fixed $p\geq1$.

\emph{Row close to a coordinate axis.} Suppose $\delta(v)\leq1/16$, so $|v_k|-\max_{i\neq k}|v_i|\geq\sqrt{15/16}-\sqrt{1/16}>1/2$. If a switch occurs, then $|q_i|>|q_k|$ for some $i\neq k$, and hence
\[
\tfrac12\sigma|G|\leq\sigma|G|(|v_k|-|v_i|)\leq(|q_k|-|q_i|)+2\tau M<2\tau M,
\]
which implies $|G|<C\tau M$. Moreover, the same switch gives $0<|q_i|-|q_k|\leq2\tau M$ and $|q_i|+|q_k|\leq C\tau M$, so $0<q_i^2-q_k^2\leq C\tau^2M^2$. Since $G$ is independent of $Z$ and has a bounded density, $\PP(|G|\leq C\tau M\mid Z)\leq C\tau M$. Conditioning on $Z$ and using the moment bound above therefore yields
\[
\EE\bigl[\max_iq_i^2-q_k^2\bigr]
\leq C\tau^2\EE\bigl[M^2\PP(|G|\leq C\tau M\mid Z)\bigr]
\leq C\tau^3\EE M^3
\leq C\tau^3\log(ed)^{3/2}.
\]

\emph{Row away from every coordinate axis.} Suppose instead that $\delta(v)>1/16$. Since $v_i^2-v_k^2\leq0$ for every $i$, expanding the difference of squared coordinates gives
\[
q_i^2-q_k^2
\leq2\sigma\tau G(v_iZ_i-v_kZ_k)+\tau^2(Z_i^2-Z_k^2)
\leq4\sigma_+\tau|G|M+2\tau^2M^2.
\]
Taking the maximum over $i$ and then expectations yields $\EE[\max_iq_i^2-q_k^2]\leq C(\tau\sqrt{\log(ed)}+\tau^2\log(ed))$. Since $\delta(v)>1/16$, this is bounded by $16C(\tau\sqrt{\log(ed)}+\tau^2\log(ed))\delta(v)$.

Combining the bounds from the two cases proves~\eqref{eq:row-switching-gain}.
\end{proof}

The row defect of $U\in\orth(d)$ is $D(U):=\sum_{j=1}^d\delta(v_j)$, summing the defects of its $d$ rows $v_j:=U^\top e_j$; $D(U)=0$ exactly when $U$ is a signed permutation. The next lemma bounds $D(U)$ below once $U$ is bounded away from every signed permutation in its own component of $\orth(d)$, and is purely geometric, independent of the noise model.

\begin{lemma}[Row defects control distance from the frame orbit]
\label{lem:row-defect-orbit-distance}
For every $r>0$, if $V\in\orth(d)$ satisfies $d_g(V,R)\geq r$ for every signed permutation $R$ in the same connected component as $V$, then
\begin{equation}\label{eq:row-defect-orbit-distance}
D(V)\geq\min\Bigl\{\frac12,\frac{2r^2}{\pi^2}\Bigr\}.
\end{equation}
\end{lemma}

\begin{proof}
If $D(V)\geq1/2$, the claimed bound holds immediately. Suppose therefore that $D(V)<1/2$. Since every row defect is nonnegative, each row has a unique entry of squared magnitude exceeding $1/2$. Moreover, two different rows cannot have their distinguished entries in the same column, since the corresponding squared entries would sum to more than one, contradicting orthogonality.

For each row $j$, let $k_j$ denote the column of its distinguished entry and set $a_j:=|V_{j,k_j}|$. Define the signed permutation $R$ by taking its unique nonzero entry in row $j$ to be $R_{j,k_j}=\operatorname{sign}(V_{j,k_j})$. Then $\langle V,R\rangle_{\Fr}=\sum_ja_j$, and hence
\[
\|V-R\|_{\Fr}^2=2\sum_j(1-a_j)\leq2\sum_j(1-a_j^2)=2D(V)<1.
\]
Thus $V$ and $R$ belong to the same connected component of $\orth(d)$. Indeed, if they belonged to different components, $R^\top V$ would have an eigenvalue $-1$, implying $\|V-R\|_{\Fr}\geq2$, contrary to the preceding bound.

Writing $\theta_\ell$ for the principal rotation angles of $R^\top V$, the geodesic and Frobenius distances satisfy $d_g(V,R)^2=2\sum_\ell\theta_\ell^2$ and $\|V-R\|_{\Fr}^2=8\sum_\ell\sin^2(\theta_\ell/2)$. Since $\sin(t/2)\geq t/\pi$ for $0\leq t\leq\pi$, these identities give
\[
\|V-R\|_{\Fr}^2\geq\frac4{\pi^2}d_g(V,R)^2.
\]
By hypothesis, $d_g(V,R)\geq r$. Combining this with the preceding inequality and $D(V)\geq\tfrac12\|V-R\|_{\Fr}^2$ yields $D(V)\geq2r^2/\pi^2$, proving~\eqref{eq:row-defect-orbit-distance}.
\end{proof}

\begin{proposition}[Dimension-uniform global gap for balanced Gaussian mixtures]
\label{prop:balanced-gaussian-global-gap}
Under the hypotheses of Proposition~\ref{prop:balanced-gaussian-curvature}, fix $0<r\leq\min\{r_0,1/8\}$, set $\Nc=B_r(\Ustar)$ and $\Norbit=\bigcup_{R\in\Spm}\Nc R$, and write $c_r:=\min\{1/2,2r^2/\pi^2\}$. There exists $d_0\geq4$, depending only on $r,\sigma_-^2,\sigma_+^2$, and $b_+$, such that the following two statements hold for every $d\geq d_0$ and $0<\tau^2\leq b_+/d$,
\begin{equation}\label{eq:balanced-global-gap}
\Ec(U)-\Ec(\Ustar)\geq\frac{\sigma_-^2c_r}{4d}
\qquad\text{for every }U\notin\Norbit,
\end{equation}
and, for the corresponding noiseless model,
\begin{equation}\label{eq:balanced-noiseless-gap}
\Delta_\nu(d,r):=\inf_{U\notin\Norbit}\bigl(\Ec_\nu(U)-\Ec_\nu(\Ustar)\bigr)\asymp d^{-1},
\end{equation}
with implicit constants depending only on $r,\sigma_-^2$, and $\sigma_+^2$. Consequently, $\Ustar\Spm$ is the complete minimizer orbit of $\Ec$, and Assumption~\ref{ass:sublevel-containment} holds with $\rho=\sigma_-^2c_r/(8d)$.
\end{proposition}

\begin{proof}
By an orthogonal change of coordinates, assume $\Ustar=I$. We first establish the noisy global gap using the rowwise switching estimate. We then combine this gap with local strong convexity to identify the minimizer orbit and verify sublevel containment. Finally, we establish the matching bounds for the noiseless gap.

Write $H_j(v):=\EE\left(\max_i|\sigma_jGv_i+\tau Z_i|^2\right)$ for the expected energy captured by row $j$ in direction $v$. Conditional on $J=j$, we have $U^\top X\overset d=\sigma_jGv_j+\tau Z$, so $\Ec(U)$ equals a constant minus $d^{-1}\sum_jH_j(v_j)$. Consequently,
\[
d\bigl(\Ec(U)-\Ec(I)\bigr)=\sum_{j=1}^d\bigl[H_j(e_j)-H_j(v_j)\bigr].
\]
For each $j$, choose $k_j$ such that $|(v_j)_{k_j}|=\|v_j\|_\infty$. Exchangeability of the Gaussian coordinates gives
$H_j(e_j)=H_j(e_{k_j})$. The expected energy captured in direction $v_j$ is $\sigma_j^2(v_j)_{k_j}^2+\tau^2$ plus a nonnegative switching gain. Applying Lemma~\ref{lem:balanced-row-switching} to this gain and discarding the nonnegative switching gain at $e_{k_j}$ therefore yields
\[
H_j(e_j)-H_j(v_j) \geq\sigma_j^2\delta(v_j) -C\bigl(\tau\sqrt{\log(ed)}+\tau^2\log(ed)\bigr)\delta(v_j) -C\tau^3\log(ed)^{3/2}.
\] 
Summing over $j$, using $\sigma_j^2\geq\sigma_-^2$, and recalling $D(U)=\sum_j\delta(v_j)$ gives
\[
\Ec(U)-\Ec(I) \geq\Bigl[\sigma_-^2-C\bigl(\tau\sqrt{\log(ed)}+\tau^2\log(ed)\bigr)\Bigr]\frac{D(U)}d -C\tau^3\log(ed)^{3/2}.
\]

Under $\tau^2\leq b_+/d$, the correction in the bracket is $\littleo(1)$, while the final term is $\littleo(d^{-1})$. We may therefore increase $d_0$, depending only on the parameters in the proposition, so that the bracket is at least $\sigma_-^2/2$ and the final term is at most $\sigma_-^2c_r/(4d)$ for every $d\geq d_0$. For $U\notin\Norbit$, Lemma~\ref{lem:row-defect-orbit-distance} gives $D(U)\geq c_r$. Combining these bounds yields
\[
\Ec(U)-\Ec(I)\geq\frac{\sigma_-^2c_r}{2d}-\frac{\sigma_-^2c_r}{4d} =\frac{\sigma_-^2c_r}{4d},
\]
which proves~\eqref{eq:balanced-global-gap}.

We next combine this global separation with the local curvature established in Proposition~\ref{prop:balanced-gaussian-curvature}. That proposition gives
$\Hess\,\Ec\succeq(\sigma_-^2/2d)I$ on each $\Nc R$, while the reflection argument in its proof makes every $IR$, $R\in\Spm$, a critical point. Geodesic strong convexity therefore yields $\Ec(U)\geq\Ec(IR)+(\sigma_-^2/4d)d_g(U,IR)^2$ for every $U\in\Nc R$. Thus $IR$ is the unique minimizer on its component. Since all orbit representatives have the same risk and~\eqref{eq:balanced-global-gap} excludes minimizers outside $\Norbit$, we obtain $\arg\min_{U\in\orth(d)}\Ec(U)=\Spm$, or $\Ustar\Spm$ before the coordinate reduction. Moreover, choosing
$\rho=\sigma_-^2c_r/(8d)$, strictly below the gap in~\eqref{eq:balanced-global-gap}, gives
\[
\{U:\Ec(U)\leq\Ec(\Ustar)+\rho\} \subseteq\Norbit,
\]
which verifies Assumption~\ref{ass:sublevel-containment}.

It remains to establish~\eqref{eq:balanced-noiseless-gap}. At $\tau=0$, the captured energy is exactly $H_j(v)=\sigma_j^2\|v\|_\infty^2$, so
\[
d\bigl(\Ec_\nu(U)-\Ec_\nu(I)\bigr)=\sum_j\sigma_j^2\delta(v_j).
\]
For $U\notin\Norbit$, Lemma~\ref{lem:row-defect-orbit-distance} therefore gives $\Delta_\nu(d,r)\geq\sigma_-^2c_r/d$. For the reverse bound, rotate coordinates $1$ and $2$ by $\phi=r/\sqrt2$, obtaining $U_\phi$ with $d_g(I,U_\phi)=r$. Since $d_g(I,U_\phi)=r$, we have $U_\phi\notin\Nc$. Moreover, for every $R\in\Spm\setminus\{I\}$,
\[
\|U_\phi-R\|_{\Fr} \geq\|I-R\|_{\Fr}-\|U_\phi-I\|_{\Fr} \geq2-r>r,
\]
so $U_\phi$ lies outside every other translate as well. Hence $U_\phi\notin\Norbit$. Only the first two rows have nonzero defects, both equal to $\sin^2\phi$, giving
\[
\Delta_\nu(d,r)\leq\Ec_\nu(U_\phi)-\Ec_\nu(I)=\frac{\sigma_1^2+\sigma_2^2}d\sin^2\Bigl(\frac r{\sqrt2}\Bigr)\leq\frac{2\sigma_+^2}d\sin^2\Bigl(\frac r{\sqrt2}\Bigr).
\]
Together with the preceding lower bound, this proves~\eqref{eq:balanced-noiseless-gap}.
\end{proof}

\subsection{Sublevel-set containment under Gaussian convolution}
\label{app:sublevel-set-containment}

Unlike the rowwise comparison above, which is specific to balanced mixtures, the argument in this section works for general axial laws, albeit for a fixed dimension. Lemma~\ref{lem:row-defect-orbit-distance} gives a positive gap for the noiseless risk, which we transfer to the Gaussian convolution model through a uniform stability bound. Combined with local strong convexity and criticality, this yields sublevel-set containment and identifies $\Ustar\Spm$ as the complete minimizer orbit for sufficiently small noise.

\begin{lemma}[Uniform continuity of the population objective under convolution]
\label{lem:objective-continuity}
Under the Gaussian convolution model, for every $\tau\geq0$,
\[
\sup_{U\in\orth(d)} \bigl|\Ec_\mu(U)-\Ec_\nu(U)\bigr| \leq \tau\sqrt{2d(2\bar a+\tau^2d)}.
\]
\end{lemma}

\begin{proof}
Let $\gamma$ be an optimal coupling of $\mu$ and $\nu$, with first marginal $\mu$ and second marginal $\nu$. For $U\in\orth(d)$, set $S_U=\bigcup_{j=1}^d\Span(u_j)$ and $f_U(x)=\dist(x,S_U)^2$. Then
\[
\Ec_\mu(U)-\Ec_\nu(U) =\int\bigl(f_U(x)-f_U(y)\bigr)\,d\gamma(x,y).
\]
Since the distance to a closed set is $1$-Lipschitz, the factorization $a^2 - b^2 = (a-b)(a+b)$ gives $|f_U(x) - f_U(y)| \leq \|x-y\|(\|x\| + \|y\|)$. Integrating against $\gamma$ and applying the Cauchy--Schwarz inequality yields
$$
\bigl|\Ec_\mu(U) - \Ec_\nu(U)\bigr|
\;\leq\;
\Bigl(\int \|x-y\|^2\,d\gamma\Bigr)^{1/2}
\Bigl(\int (\|x\|+\|y\|)^2\,d\gamma\Bigr)^{1/2}.
$$
Applying $(a+b)^2 \leq 2(a^2+b^2)$ to the second factor and noting that $\gamma$ has marginals $\mu$ and $\nu$, this becomes
$$
\bigl|\Ec_\mu(U) - \Ec_\nu(U)\bigr|
\;\leq\;
W_2(\mu,\nu)\,\sqrt{2\bigl(\EE_\mu\|X\|^2+\EE_\nu\|Y\|^2\bigr)}.
$$
Since $Z$ is isotropic with covariance $I_{d\times d}$, we have $\EE_\mu\|X\|^2 = \bar a + \tau^2 d$, giving the factor $\sqrt{2(2\bar a + \tau^2 d)}$. This bound is independent of $U$, so taking the supremum over $\orth(d)$ and substituting the bound $W_2(\mu,\nu) \leq \tau\sqrt{d}$ completes the proof.
\end{proof}

\begin{proposition}[Sublevel-set containment under Gaussian convolution]
\label{prop:sublevel-set}
Fix $0<r\leq r_0$ as in Lemma~\ref{axis-case-exact}, and set $\Nc_0=B_r(\Ustar)$ and $\Norbit^0:=\bigsqcup_{R\in\Spm}\Nc_0R$. Define
\[
\Delta_\nu:= \inf_{U\notin\Norbit^0} \bigl(\Ec_\nu(U)-\Ec_\nu(\Ustar)\bigr).
\]
Then $\Delta_\nu>0$. Therefore, if $\tau_1$ is the curvature threshold of Theorem~\ref{thm:strong-cvx}, and we set
\[
\tau_0:=\min\left\{\tau_1,\,1,\,\frac{\Delta_\nu}{8\sqrt{2d(2\bar a+d)}}\right\},
\]
it follows that, for every $0<\tau<\tau_0$,
\[
\Srhomu(\mu):=\left\{U\in\orth(d): \Ec_\mu(U)\leq\Ec_\mu(\Ustar)+\rhomu\right\} \subseteq\Norbit^0, \quad\text{with}\quad \rhomu:=\frac{\Delta_\nu}{4}.
\]
Consequently,
\[
\arg\min_{U\in\orth(d)}\Ec_\mu(U)=\Ustar\Spm.
\]
\end{proposition}

\begin{proof}
For $U\in\orth(d)$, set $V:=\UstarT U$. The noiseless risk satisfies
\[
\Ec_\nu(U)-\Ec_\nu(\Ustar) =\sum_{j=1}^d a_j \left(1-\max_i|V_{ji}|^2\right) \geq a_{\min}D(V).
\]
If $U\notin\Norbit^0$, then $d_g(V,R)\geq r$ for every signed permutation $R$ in the same connected component as $V$. Lemma~\ref{lem:row-defect-orbit-distance} therefore gives
\[
\Delta_\nu\geq a_{\min}\min\left\{\frac12,\frac{2r^2}{\pi^2}\right\}>0.
\]

Let $\etau:=\sup_{U\in\orth(d)}|\Ec_\mu(U)-\Ec_\nu(U)|$. Since $\tau<\tau_0\leq1$, Lemma~\ref{lem:objective-continuity} gives $\etau\leq\tau\sqrt{2d(2\bar a+d)}\leq\Delta_\nu/8$. Hence, for $U\in\Srhomu(\mu)$, two applications of the uniform discrepancy bound give
\[
\begin{aligned}
\Ec_\nu(U)
&\leq\Ec_\mu(U)+\etau
\leq\Ec_\mu(\Ustar)+\rhomu+\etau\\
&\leq\Ec_\nu(\Ustar)+\rhomu+2\etau
\leq\Ec_\nu(\Ustar)+\frac{\Delta_\nu}{2}
<\Ec_\nu(\Ustar)+\Delta_\nu.
\end{aligned}
\]
By the definition of $\Delta_\nu$, this strict inequality implies $U\in\Norbit^{0}$. The strict half-gap is essential here, since using the attained gap itself would not exclude the complement's boundary points.

Finally, since $\tau<\tau_1$, Theorem~\ref{thm:strong-cvx} gives geodesic strong convexity on each $\Nc_0R$, while Corollary~\ref{thm:crit-pt} establishes criticality at every $\Ustar R$. Thus each orbit representative is the unique minimizer on its aligned component. Since every global minimizer belongs to $\Norbit^0$, the complete minimizer set is $\Ustar\Spm$.
\end{proof}

\subsection{Far-region switching margin under Gaussian convolution}

We next verify the far-region switching-margin bound. Unlike the curvature and sublevel-set estimates, this bound holds for every $\tau>0$, without a small-noise restriction.

\begin{lemma}[Gaussian convolution satisfies the switching margin]
\label{lem:convolution-switching-margin}
Consider Example~\ref{ex:gaussian-convolution}, let $\Nc$ be an aligned component with $\Ustar\in\overline{\Nc}$ and $r_{\Nc}=\sup_{U\in\Nc}d_g(U,\Ustar)$, set $\theta:=2r_{\Nc}\in(0,\pi/2)$, and let $0<\eta_0\leq1/2$. For $d\geq4$, also assume $1/\sqrt2+r_{\Nc}+\eta_0<1$ and set
\[
b_0:=\sqrt{1-(1/\sqrt2+r_{\Nc}+\eta_0)^2}\qquad\text{and}\qquad
\Xi(\theta):=\frac{\cos\theta}{\sin^2\theta}+\operatorname{artanh}(\cos\theta).
\]
Then, for every $\tau>0$, the probability inequality in Assumption~\ref{ass:switching-margin} holds with
\[
\Cmar\leq M\tau
\begin{cases}
\dfrac{4\arcsin(\eta_0)}{\eta_0\sqrt{2\pi}}\sin^{-2}(\theta),&d=2,\\[1ex]
\sqrt{2\pi}\sin^{-3}(\theta),&d=3,\\[1ex]
\sqrt{\dfrac2\pi}\dfrac{d-2}{b_0}\,\Xi(\theta),&d\geq4.
\end{cases}
\]
In particular, for fixed $r_{\Nc}$ and compatible $\eta_0$, $\Cmar=\bigO(M\tau d)$ as $d\to\infty$.
\end{lemma}

\begin{proof}
Set $V:=X/\|X\|$. Condition on $J=j$ and write $u:=\ustar_j$. Integrating the Gaussian density along $\Span(u)$ and using the density bound $M$ gives
\[
p_{X\mid J=j}(x)\leq M(2\pi\tau^2)^{-(d-1)/2}
\exp\!\left(-\frac{\|P_{u^\perp}x\|^2}{2\tau^2}\right).
\]
Consequently, the conditional angular density $q_j(v):=\int_0^\infty p_{X\mid J=j}(rv)r^{d-1}\,dr$ satisfies
\begin{equation}\label{eq:gaussian-angular-component}
q_j(v)\leq A_dM\tau
\bigl(1-|\langle v,u\rangle|^2\bigr)^{-d/2},\quad\text{with}\quad
A_d:=\frac{\Gamma(d/2)}{\sqrt2\,\pi^{(d-1)/2}}.
\end{equation}
Since $\Ustar\in\overline{\Nc}$, letting $U\to\Ustar$ in the caps defining $\Rnear(\Nc)$ shows that every $v\in\Rfar(\Nc)$ satisfies $|\langle v,u\rangle|\leq\cos\theta$. When $d=2$ or $3$, we combine~\eqref{eq:gaussian-angular-component} with the exact spherical-band areas $4\arcsin h$ on $S^1$ and $4\pi h$ on $S^2$, respectively. Since $\arcsin h/h$ is increasing, this gives the first two constants.

Now let $d\geq4$. Fix $w\in S^{d-1}$ with $\dist(w,\mathcal W_{\Nc})\leq\eta_0$ and decompose $w=\alpha u+bz$, where $z\in S^{d-1}\cap u^\perp$ and $b=\|P_{u^\perp}w\|$. By compactness, we may choose $w_0\in\mathcal W_{\Nc}$ with $\|w-w_0\|\leq\eta_0$. This $w_0$ is a signed-bisector normal of some $U\in\overline{\Nc}$, so
\[
|\langle w_0,u\rangle|
\leq|\langle w_0,u_j^{(U)}\rangle|+\|u-u_j^{(U)}\|
\leq\frac1{\sqrt2}+r_{\Nc},
\]
where the first inner product is either $0$ or $1/\sqrt2$ and the second term is at most $\|\Ustar-U\|_{\Fr}\leq r_{\Nc}$. Hence $|\langle w,u\rangle|\leq1/\sqrt2+r_{\Nc}+\eta_0$ and $b\geq b_0$. We parametrize
\[
v=tu+\sqrt{1-t^2}\,y,\qquad\text{for}\quad
-\cos\theta\leq t\leq\cos\theta\quad\text{and}\quad y\in S^{d-2}\subset u^\perp.
\]
For fixed $t$, the surface measure of those $y$ satisfying $|\alpha t+b\sqrt{1-t^2}\langle y,z\rangle|\leq h$ is at most $2h\omega_{d-2}/(b\sqrt{1-t^2})$. Indeed, for $d\geq4$, the coordinate map $y\mapsto\langle y,z\rangle$ pushes forward surface measure on $S^{d-2}$ to a measure with density $\omega_{d-2}(1-s^2)^{(d-4)/2}\leq\omega_{d-2}$. Combining this with~\eqref{eq:gaussian-angular-component} and the surface element $(1-t^2)^{(d-3)/2}dt\,dS(y)$ yields
\begin{align*}
\PP\!\left(V\in\Rfar(\Nc),|\langle V,w\rangle|\leq h\mid J=j\right)
&\leq\frac{2A_dM\tau h\omega_{d-2}}{b_0}
\int_{-\cos\theta}^{\cos\theta}(1-t^2)^{-2}\,dt\\
&=\sqrt{\frac2\pi}\frac{(d-2)M\tau}{b_0}\,\Xi(\theta)\,h.
\end{align*}
The bound is uniform in $j$, and averaging with the mixture weights $p_j$ completes the proof.
\end{proof}

\subsection{Proofs of the model-verification propositions}

\begin{proof}[Proof of Proposition~\ref{prop:gaussian-convolution-assumptions}]
Fix $d\geq2$ and $0<r\leq r_0$ as in Lemma~\ref{axis-case-exact}, set $\Nc=B_r(\Ustar)$, and let $0<\tau<\tau_0$. We verify the five assumptions with the constants stated in Proposition~\ref{prop:gaussian-convolution-assumptions}.

For Assumption~\ref{ass:concentration}, we have
$\EE\|Y\|^2=\bar a$ and $\|Y\|^2=S_J^2$. Since
$\|Y+\tau Z\|^4\leq8\|Y\|^4+8\tau^4\|Z\|^4$,
\[
\mtwo=\bar a+\tau^2d
\qquad\text{and}\qquad
\mfour\leq8\sum_{j=1}^dp_j\EE S_j^4+8\tau^4d(d+2),
\]
which is finite because $\|S_j^2\|_{\psi_1}<\infty$ gives $\EE S_j^4<\infty$. 

For Assumption~\ref{ass:subexp-norms}, the mixture bound $\EE e^{S_J^2/s}=\sum_jp_j\EE e^{S_j^2/s}\leq2$ for $s\geq\max_j\|S_j^2\|_{\psi_1}$ gives $\|S_J^2\|_{\psi_1}\leq\max_j\|S_j^2\|_{\psi_1}$. The Gaussian moment-generating function shows that $\|Z\|^2$ has finite $\psi_1$ norm for fixed $d$, and $\|X\|^2\leq2\|Y\|^2+2\tau^2\|Z\|^2$. Homogeneity and the triangle inequality for the $\psi_1$ norm thus give
\[
K\leq2\max_j\|S_j^2\|_{\psi_1}+2\tau^2\bigl\|\|Z\|^2\bigr\|_{\psi_1}.
\]

For Assumption~\ref{ass:strong-convexity}, Theorem~\ref{thm:strong-cvx} gives the aligned component $\Nc$, strong convexity with $m=a_{\min}/2$, and the intrinsic gradient-Lipschitz bound with $L=4a_{\max}$, since $\tau<\tau_1$. Corollary~\ref{thm:crit-pt} supplies criticality and strict local minimality at every orbit representative. For Assumption~\ref{ass:sublevel-containment}, Proposition~\ref{prop:sublevel-set} gives the invariant sublevel containment with $\rho=\Delta_\nu/4$, since $\tau<\tau_2$, and identifies the complete minimizer orbit.

Finally, for Assumption~\ref{ass:switching-margin}, Lemma~\ref{axis-case-exact} gives $c_{\Nc}=\cos(2r)-2r>1/\sqrt2$. For $d\geq4$, we choose $\eta_0\in(0,1/2]$ with $\eta_0<1-1/\sqrt2-r$, which is possible because $r\leq1/8$; for $d=2,3$, any $\eta_0\in(0,1/2]$ works. Lemma~\ref{lem:convolution-switching-margin} then gives the required probability bound with the stated $\Cmar$, for every $\tau>0$.
\end{proof}

\begin{proof}[Proof of Proposition~\ref{prop:balanced-gaussian-mixture-uniform}]
Assumption~\ref{ass:concentration} follows from the finite Gaussian moments. Take $d_0\geq4$ large enough for both Propositions~\ref{prop:balanced-gaussian-curvature} and~\ref{prop:balanced-gaussian-global-gap}, and let $d\geq d_0$. Set $r:=\min\{r_0,1/8\}$ and $\Nc=B_r(\Ustar)$, with $r_0$ as in the former proposition. The curvature proposition gives Assumption~\ref{ass:subexp-norms} with the stated $K$, and Assumption~\ref{ass:strong-convexity} on $\Nc$ with $m=\sigma_-^2/(2d)$ and $L=4\sigma_+^2/d$, since the Hessian bounds \eqref{eq:dimension-uniform-hessian} persist on the smaller ball and it has unique geodesics, disjoint translates, and $\dframe=d_g$. For Assumption~\ref{ass:switching-margin}, note that $r\leq1/8$ gives $c_{\Nc}>1/\sqrt2$, and that $\eta_0:=(1-1/\sqrt2-r)/2\in(0,1/2]$ satisfies $1/\sqrt2+r+\eta_0<1$. Each $\lambda_j=\Nc(0,\sigma_j^2)$ has a density bounded by $M=(2\pi\sigma_-^2)^{-1/2}$, independently of $d$. Since $\theta=2r$, $b_0$ and $\Xi(\theta)$ depend only on $r$, and hence only on $\sigma_-^2,\sigma_+^2$ and $b_+$, Lemma~\ref{lem:convolution-switching-margin} gives $\Cmar\leq M\tau\sqrt{2/\pi}\,(d-2)\,\Xi(\theta)/b_0\leq\Cgeom(r,\sigma_-)\,\tau d$, where $\Cgeom(r,\sigma_-):=M\sqrt{2/\pi}\,\Xi(\theta)/b_0$, using $d-2<d$. Finally, Proposition~\ref{prop:balanced-gaussian-global-gap} gives Assumption~\ref{ass:sublevel-containment} on this same $\Nc$ with $\rho=\sigma_-^2c_r/(8d)$, and identifies $\Ustar\Spm$ as the complete minimizer orbit.
\end{proof}

\section{Technical proofs for Theorem~\ref{thm:convergence}}\label{app:erm-proofs}

\subsection{Subdifferential analysis of the empirical objective}\label{app:subdiff}

This subsection records the standard Clarke calculus needed for the empirical risk and the loss-specific consequences used below. We use the Euclidean theory of \cite{clarke1990optimization} together with its Riemannian extension in \cite{HosseiniPouryayevali2011}.

\begin{lemma}[Lipschitz continuity]\label{lem:lipschitz-empirical-obj}
Each $f_{x_i}(U)=\|x_i-P_U(x_i)\|_2^2$ is globally Lipschitz on $\orth(d)$ with respect to the ambient Frobenius distance, with constant $2\|x_i\|_2^2$. Consequently, $\Ec_n$ is globally Lipschitz with constant $L_n:=\frac{2}{n}\sum_{i=1}^n\|x_i\|_2^2$, and $\EE[L_n]=2\mtwo$.
\end{lemma}

\begin{proof}
Write $f_{x_i}(U)=\min_j g_{i,j}(U)$, where $g_{i,j}(U)=\|x_i\|_2^2-(u_j^\top x_i)^2$. For $U,U'\in\orth(d)$,
\[
|g_{i,j}(U)-g_{i,j}(U')|
\leq 2\|x_i\|_2^2\|U-U'\|_{\Fr}.
\]
A finite minimum preserves a common Lipschitz constant. Averaging over $i$ gives $L_n$, and its expectation follows by linearity.
\end{proof}

Recall the Riemannian Clarke subdifferential $\partial_Rf(U)$ from Section~\ref{sec:background}. We now specialize this definition to our loss.

\begin{proposition}[Clarke formulas and canonical selection]
\label{prop:clarke-formulas}
Recall $J(x;U)$, $j^\star(x;U):=\min J(x;U)$, $g_{\mathrm{can}}(x,U):=\gamma_{j^\star(x;U)}(x,U)$, and $G_n(U):=n^{-1}\sum_i g_{\mathrm{can}}(x_i,U)$ from Section~\ref{subsec:empirical-nonsmooth}. Then
\[
\partial_R f_x(U)=\conv\{\gamma_j(x,U):j\in J(x;U)\}.
\]
Moreover, $g_{\mathrm{can}}$ is Borel and belongs to $\partial_Rf_x(U)$; if Assumptions~\ref{ass:concentration} and~\ref{ass:strong-convexity} hold, then for every fixed $U\in\Nc$, $g_{\mathrm{can}}(x,U)$ equals $\Rgrad f_x(U)$ almost surely. For every sample realization,
\[
\partial_R\Ec_n(U)=\frac1n\sum_{i=1}^n\partial_Rf_{x_i}(U)\qquad\text{and}
\qquad G_n(U)\in\partial_R\Ec_n(U).
\]
\end{proposition}

\begin{proof}
\emph{Pointwise formula.}\enspace
Set $\phi_j(U):=(u_j^\top x)^2$. Since
$f_x(U)=\|x\|_2^2-\max_j\phi_j(U)$, the finite-maximum Clarke rule \cite[Proposition~2.3.12]{clarke1990optimization}, applied in a local chart and transferred to the tangent space as in \cite{HosseiniPouryayevali2011}, gives
\[
\partial_R f_x(U)
=\conv\{-\Rgrad\phi_j(U):j\in J(x;U)\}
=\conv\{\gamma_j(x,U):j\in J(x;U)\}.
\]
\emph{Empirical formula.}\enspace
Set $\phi_{ij}(U):=(u_j^\top x_i)^2$ and $H_n(U):=n^{-1}\sum_i\max_j\phi_{ij}(U)$. Separability gives
\[
H_n(U)=\max_{\alpha\in[d]^n}\frac1n\sum_i\phi_{i\alpha_i}(U).
\]
The active multi-indices are $\prod_iJ(x_i;U)$, so another application of the finite-maximum rule and the finite-set Minkowski identity gives
\[
\partial_RH_n(U)
=\frac1n\sum_i\conv\{\Rgrad\phi_{ij}(U):j\in J(x_i;U)\}.
\]
Negation proves the empirical formula, which contains $G_n(U)$.

\emph{Measurability.}\enspace
The score functions $s_j(x,U):=|u_j^\top x|^2$ are continuous and
\[
\{j^\star=j\}=\bigcap_{k<j}\{s_j>s_k\}\cap\bigcap_{k>j}\{s_j\geq s_k\}.
\]
Thus $g_{\mathrm{can}}$ is Borel by finite pasting. For $U\in\Nc$, Lemma~\ref{lem:diff-implies-tie-null} shows that nonzero active ties are null. At $x=0$, both $g_{\mathrm{can}}(0,U)$ and $\Rgrad f_0(U)$ vanish.
\end{proof}

The empirical formula above also yields Lemma~\ref{lem:empirical-local-min} on empirical local minimizers.

\begin{proof}[Proof of Lemma~\ref{lem:empirical-local-min}]
We retain the notation $\phi_{ij}$, $H_n$, and active multi-indices $\prod_iJ(x_i;V)$ from the previous proof. Since $\Ec_n(U)=n^{-1}\sum_i\|x_i\|_2^2-H_n(U)$, the point $V$ is a local maximizer of $H_n$. For every active multi-index $\alpha\in\prod_iJ(x_i;V)$ and every $U$ sufficiently close to $V$,
\[
\frac{1}{n}\sum_{i=1}^n\phi_{i\alpha_i}(U)\leq H_n(U)\leq H_n(V)=\frac{1}{n}\sum_{i=1}^n\phi_{i\alpha_i}(V),
\]
so the smooth function $U\mapsto n^{-1}\sum_i\phi_{i\alpha_i}(U)$ has a local maximum at $V$. Since $\Rgrad\phi_{ij}=-\gamma_j(x_i,\cdot)$, this gives
\begin{equation}\label{eq:active-stationary}
\frac{1}{n}\sum_{i=1}^n\gamma_{\alpha_i}(x_i,V)=0\qquad\text{for every }\alpha\in\prod_{i=1}^nJ(x_i;V).
\end{equation}
In particular, the canonical multi-index $\alpha_i=j^\star(x_i;V)$ gives $G_n(V)=0$. Assume, then, that some nonzero sample $x_k$ has two distinct active indices $i,j\in J(x_k;V)$. We choose active multi-indices $\alpha,\alpha'\in\prod_{\ell=1}^n J(x_\ell;V)$ that agree except at coordinate $k$, where $\alpha_k=i$ and $\alpha'_k=j$. Subtracting the corresponding identities in~\eqref{eq:active-stationary} gives $\gamma_i(x_k,V)-\gamma_j(x_k,V)=0$. This contradicts~\eqref{eq:tie-rotation} applied at the frame $V$, whose right-hand side is nonzero whenever $x_k\neq0$ has two distinct active indices. Thus every nonzero sample has a unique active index at $V$. For zero samples, all branch gradients vanish. Consequently,
\[
\conv\{\gamma_j(x_k,V):j\in J(x_k;V)\}=\{\gamma_{j^\star(x_k;V)}(x_k,V)\}
\]
for every $k$, and Proposition~\ref{prop:clarke-formulas}, together with~\eqref{eq:active-stationary}, yields $\partial_R\Ec_n(V)=\{G_n(V)\}=\{0\}$. Finally, by continuity of the scores and finiteness of the sample, the indices $j^\star(x_i;V)$ of the nonzero samples remain unchanged on a common neighborhood of $V$. On that neighborhood,
\[
\Ec_n(U)=\frac{1}{n}\sum_{i=1}^n\Bigl(\|x_i\|_2^2-\phi_{i\,j^\star(x_i;V)}(U)\Bigr),
\]
since zero samples contribute zero regardless of the selected index. Thus $\Ec_n$ is real analytic there.
\end{proof}

\begin{proposition}[Subdifferential norm bound]
\label{prop:individual-clarke-bounds}
For $A=[a_1,\ldots,a_d]\in\RR^{d\times d}$, let
$\widetilde f_x(A):=\|x\|_2^2-\max_j(a_j^\top x)^2$ and write
$\partial_Ef_x(U):=\partial_C\widetilde f_x(U)$. Then, for every $x\in\RR^d$ and $U\in\orth(d)$, every $\xi\in\partial_Ef_x(U)\cup\partial_Rf_x(U)$ satisfies $\|\xi\|_{\Fr}\leq2\|x\|_2^2$. 
\end{proposition}

\begin{proof}
    For each $j\in J(x;U)$, the Euclidean gradient of the corresponding active branch is $-2(u_j^\top x)xe_j^\top$, whose Frobenius norm is bounded by $2\|x\|_2^2$. Since the Euclidean Clarke subdifferential is the convex hull of these gradients, every $A\in\partial_E f_x(U)$ satisfies $\|A\|_{\Fr}\leq2\|x\|_2^2$. By Proposition~\ref{prop:clarke-formulas}, $\partial_R f_x(U)=\Pi_U(\partial_E f_x(U))$. Since the orthogonal projection $\Pi_U$ is nonexpansive, $\sup_{\xi\in\partial_R f_x(U)}\|\xi\|_{\Fr}\leq2\|x\|_2^2.$
\end{proof}

\subsection{Pointwise concentration of empirical subgradients}\label{app:pointwise}

At a fixed frame $U\in\Nc$, Lemma~\ref{lem:diff-implies-tie-null} makes all subgradient selections agree almost surely. Hence the nonsmoothness of $\Ec_n$ is invisible at a single frame, and ordinary i.i.d. concentration gives the parametric rate.

\begin{proposition}[Pointwise subgradient concentration]
\label{prop:pointwise-bound}
Let $x_1,\ldots,x_n\sim_{\mathrm{i.i.d.}}\mu$, where $\mu$ satisfies Assumptions~\ref{ass:concentration} and~\ref{ass:strong-convexity}. Fix any deterministic $U\in\Nc$ before drawing the sample, and let $G_n(U)=n^{-1}\sum_i g_{\mathrm{can}}(x_i,U)$. The following statements hold.
\begin{enumerate}[label=\arabic*., leftmargin=2.2em, labelsep=0.6em]
\item \textbf{(Mean-square concentration)} 
  The positive semidefinite matrix
  \[
  \Sigma := \EE[(g_{\mathrm{can}}(x,U)-\Rgrad\Ec(U))(g_{\mathrm{can}}(x,U)-\Rgrad\Ec(U))^\top]
  \]
  is well-defined, and
  \[\EE\bigl[\|G_n(U) - \Rgrad\Ec(U)\|_{\Fr}\bigr]
    \leq \sqrt{\frac{\Tr(\Sigma)}{n}}
    \leq 2\sqrt{\frac{\mfour}{n}}.\]
  \item \textbf{(High-probability bound)} Under Assumption~\ref{ass:subexp-norms}, for every $t>0$,
  \[
  \PP\!\left(
    \|G_n(U) - \Rgrad\Ec(U)\|_{\Fr}
    \geq (\sqrt{8\mfour}+\sqrt2\,\Gamma)\sqrt{\frac{t}{n}}
    +\frac{\Gamma t}{n}
  \right)
  \leq 2e^{-t/4},
  \]
  where $\Gamma := C_1(2K + 2C_0\sqrt{\mfour})$ for absolute constants $C_1, C_0 > 0$.
\end{enumerate}
\end{proposition}
\begin{proof}
In this proof, we abbreviate the canonical selection $g_{\mathrm{can}}$ from the proposition as $g$. Since $U\in\Nc$, Lemma~\ref{lem:diff-implies-tie-null} gives
\(
\Rgrad\Ec(U)=\EE[g(x,U)],
\)
and the same identity holds for every measurable active-gradient selection at this fixed frame. The measurable selection $g$ satisfies 
$\|g(x,U)\|_{\Fr}\leq2\|x\|_2^2$
by Proposition~\ref{prop:individual-clarke-bounds}.
Assumption~\ref{ass:concentration} therefore makes $g(x,U)$
square-integrable, so $\Sigma$ is well-defined.
Jensen's inequality also gives
$\|\EE g(x,U)\|_{\Fr}\leq2\sqrt{\mfour}$.
By the definition $\Sigma=\EE[(g(x,U)-\EE g(x,U))(g(x,U)-\EE g(x,U))^\top]$ and linearity of trace,
\[
\Tr(\Sigma)=\EE\|g(x,U)-\EE g(x,U)\|_{\Fr}^2
=\EE\|g(x,U)\|_{\Fr}^2-\|\EE g(x,U)\|_{\Fr}^2
\leq4\mfour.
\]

\emph{Mean-square concentration.}\enspace
Since $\EE[G_n(U)] = \Rgrad\Ec(U)$, independence of the $x_i$ gives
\begin{align*}
\EE\left[\|G_n(U) - \Rgrad\Ec(U)\|_{\Fr}^2\right]
&= \frac{1}{n^2}\sum_{i=1}^n\EE\left[\|g(x_i,U) - \Rgrad\Ec(U)\|_{\Fr}^2\right]
= \frac{\Tr(\Sigma)}{n},
\end{align*}
where the cross-terms vanish by independence. Jensen's inequality gives
\[
\EE[\|G_n(U) - \Rgrad\Ec(U)\|_{\Fr}]
\leq \sqrt{\Tr(\Sigma)/n},
\]
which proves part~1.

\emph{High-probability bound.}\enspace
Let $Y_i := g(x_i,U) - \Rgrad\Ec(U)$, which are independent, mean-zero elements of $(T_U \orth(d), \langle\cdot,\cdot\rangle_{\Fr})$, so that $G_n(U) - \Rgrad\Ec(U) = \frac{1}{n}\sum_{i=1}^n Y_i$. The almost sure bound $\|Y_i\|_{\Fr} \leq 2\|x_i\|_2^2 + 2\sqrt{\mfour}$, together with the sub-exponentiality of $\|x_i\|_2^2$ under Assumption~\ref{ass:subexp-norms} and the Orlicz norm properties $\|aW\|_{\psi_1} = |a|\|W\|_{\psi_1}$ and $\|W + c\|_{\psi_1} \leq \|W\|_{\psi_1} + C_0|c|$, yields
\[
\|\|Y_i\|_{\Fr}\|_{\psi_1} \leq 2K + 2C_0\sqrt{\mfour}.
\]
The standard moment bound $\EE[|W|^m] \leq m!(C_1\|W\|_{\psi_1})^m$ then gives, for all $m \geq 2$,
\[
\EE[\|Y_i\|_{\Fr}^m] \leq m!\,\Gamma^m\qquad\text{with} \qquad \Gamma := C_1(2K + 2C_0\sqrt{\mfour}).
\]

Set $v^2 := n\EE[\|Y_1\|_{\Fr}^2] \leq 4n\mfour$ and $\varsigma^2 := 2\max\{v^2, n\Gamma^2\}$. The moment-growth condition~\eqref{eq:pinelis-moment-growth} holds since
\[
\sum_{i=1}^n\EE[\|Y_i\|_{\Fr}^m] \leq n\cdot m!\,\Gamma^m = m!\,\Gamma^{m-2}(n\Gamma^2) \leq m!\,\Gamma^{m-2}\frac{\varsigma^2}{2}.
\]
Lemma~\ref{lem:pinelis-hilbert-bennett} then gives, for all $r \geq 0$,
\[
\PP\!\left(\left\|\sum_{i=1}^n Y_i\right\|_{\Fr} \geq r\right) \leq 2\exp\!\left(-\frac{1}{4}\min\!\left\{\frac{r^2}{\varsigma^2}, \frac{r}{\Gamma}\right\}\right).
\]
For $t>0$, set
\(
s:=\frac{\varsigma\sqrt t+\Gamma t}{n}.
\)
For $r=ns$, both $r^2/\varsigma^2\geq t$ and $r/\Gamma\geq t$, and hence
\[
\PP\!\left(\|G_n(U)-\Rgrad\Ec(U)\|_{\Fr}\geq s\right)\leq2e^{-t/4}.
\]
Using $\varsigma\leq\sqrt{2v^2}+\sqrt{2n}\,\Gamma$ and $v^2\leq4n\mfour$ gives
\[
s\leq(\sqrt{8\mfour}+\sqrt2\,\Gamma)\sqrt{\frac{t}{n}}
+\frac{\Gamma t}{n},
\]
which proves part~2.
\end{proof}

The remaining subsections upgrade this fixed-frame bound to a uniform bound on the aligned component $\Nc$, which must account for switching.

\begin{lemma}[Maximum squared norm control]
\label{lem:max-norm}
Under Assumption~\ref{ass:subexp-norms}, for $n\geq3$, with probability at least $1-1/n$, the inequality
\(
\max_{i \leq n}\|x_i\|_2^2 \leq 3K\log n
\)
holds.
\end{lemma}

\begin{proof}
If $K=0$, the claim is immediate. Otherwise, the definition of $K$ and a union bound give
\[
\PP\!\left(\max_{i\leq n}\|x_i\|_2^2>K\log(2n^2)\right)
\leq n\,\EE\!\left[e^{\|x\|_2^2/K}\right]e^{-\log(2n^2)}
\leq\frac1n.
\]
Since $\log(2n^2)\leq3\log n$ for $n\geq2$, the result follows.
\end{proof}

\subsection{Geometry of switching sets}\label{app:switching}

For consistently labeled frames $U_1,U_2\in\orth(d)$, write $\Sswitch(U_1,U_2):=\{v\in S^{d-1}:j^\star(v;U_1)\neq j^\star(v;U_2)\}$ for the switching set used in Theorem~\ref{thm:switching-margin}. Throughout this subsection, $\dist_g$ denotes the geodesic (angular) distance on $S^{d-1}$, and, for $A\subseteq S^{d-1}$ and $a>0$, we write $\tube_a(A):= \{v\in S^{d-1}:\dist_g(v,A)\leq a\}$ for the closed geodesic tube of radius $a$ around $A$. 

\subsubsection{Voronoi stability under frame perturbations}\label{app:voronoi}

\begin{lemma}[Voronoi stability under perturbations]
\label{lem:voronoi-stability}
Let $U_1, U_2 \in \orth(d)$ with $\|U_1 - U_2\|_F \leq \delta$. If $v \in S^{d-1}$ satisfies $\dist_g(v, \Ac(U_1)) > (\pi/\sqrt{2})\delta$, then $j^\star(v;U_1) = j^\star(v;U_2)$.
\end{lemma}

\begin{proof}
Write $c_\ell:=\langle v,u_\ell^{(1)}\rangle$. The hypothesis ensures that $v\notin\Ac(U_1)$, so its winning index $j$ is unique. Set $\alpha:=|c_j|$ and $\beta:=\max_{\ell\neq j}|c_\ell|$, and choose a runner-up $k$ with $|c_k|=\beta$.

\emph{Exact margin-distance identity.} Consider the closed winner cell
\[
C_j:=\left\{z\in S^{d-1}:|\langle z,u_j^{(1)}\rangle|\geq\max_{\ell\neq j}|\langle z,u_\ell^{(1)}\rangle|\right\}.
\]
Its boundary is contained in $\Ac(U_1)$. Conversely, every ambiguity point in $C_j$ lies on its boundary, while any geodesic from $v\in\operatorname{int}(C_j)$ to an ambiguity point outside $C_j$ must cross $\partial C_j$.
Consequently,
\[
d_0:=\dist_g(v,\Ac(U_1))=\dist_g(v,\partial C_j).
\]
The boundary $\partial C_j$ lies in the union of the signed pairwise bisectors between axis $j$ and each $\ell\neq j$. For fixed $\ell$, the distance to the nearer bisector is $\arcsin((\alpha-|c_\ell|)/\sqrt2)$, by the identity
\[
\dist_g\bigl(v,\{z\in S^{d-1}:\langle z,n\rangle=0\}\bigr)=\arcsin|\langle v,n\rangle|
\]
for every unit normal $n$. Thus,
\[
d_0\geq\arcsin\!\left(\frac{\alpha-\beta}{\sqrt2}\right).
\]
To show that this bound is attained, choose $s_j,s_k\in\{-1,1\}$ with $s_jc_j=\alpha$ and $s_kc_k=\beta$ (either sign if $\beta=0$), and set
\[
q:=\frac{\alpha-\beta}{\sqrt2},\qquad n:=\frac{s_ju_j^{(1)}-s_ku_k^{(1)}}{\sqrt2}, \qquad
p:=v-qn,\quad\text{and}\quad w:=\frac{p}{\|p\|_2}.
\]
Since $\langle v,n\rangle=q$, we have $\|p\|_2=\sqrt{1-q^2}$. Moreover,
\[
s_j\langle p,u_j^{(1)}\rangle =s_k\langle p,u_k^{(1)}\rangle =\frac{\alpha+\beta}{2}\quad\text{and}\quad
|\langle p,u_\ell^{(1)}\rangle| \leq\beta\leq\frac{\alpha+\beta}{2}, \quad\text{for }\ell\notin\{j,k\}.
\]
Thus $j$ and $k$ tie for the largest absolute projection at $w$, so $w\in\partial C_j$. Since $\dist_g(v,w)=\arcsin(q)$, the preceding lower bound is attained, proving
\[
\alpha-\beta=\sqrt2\sin d_0\qquad\text{for }\quad 0<d_0\leq\frac{\pi}{4}.
\]

\emph{Stability under perturbation.} The identity above, together with $\sin t\geq2t/\pi$ on $[0,\pi/2]$ and $d_0>(\pi/\sqrt2)\delta$, gives $\alpha-\beta>2\delta$. On the other hand, for every $\ell$,
\[
\bigl| |\langle v,u_\ell^{(1)}\rangle| -|\langle v,u_\ell^{(2)}\rangle| \bigr| \leq\|u_\ell^{(1)}-u_\ell^{(2)}\|_2 \leq\delta.
\]
Consequently, the absolute projection onto axis $j$ under $U_2$ is at least $\alpha-\delta$, whereas every competing projection is at most $\beta+\delta$. Since $\alpha-\delta>\beta+\delta$, the winning index remains unique and unchanged, proving $j^\star(v;U_2)=j^\star(v;U_1)$.
\end{proof}

\begin{corollary}[Switching localization]
\label{cor:switching-localization}
If $U_1, U_2 \in \orth(d)$ satisfy $\|U_1 - U_2\|_F \leq \delta$, then
\[
\Sswitch(U_1, U_2) \subseteq \tube_{\frac{\pi}{\sqrt{2}}\delta}(\Ac(U_1)).
\]
\end{corollary}
This follows immediately as the contrapositive of Lemma~\ref{lem:voronoi-stability}. Covering this tube efficiently is the task of the next two subsubsections.

\subsubsection{Sphere packing and stripe construction}\label{app:stripes}

We denote by $\omega_d = 2\pi^{d/2}/\Gamma(d/2)$ the surface area of the unit sphere $S^{d-1}$ in $\RR^d$. The construction of the universal stripe collection rests on a maximal packing of $S^{d-1}$, whose covering and separation properties we record first.

\begin{proposition}[Voronoi cover from maximal packing]
\label{prop:packing-partition}
Let $0<\delta\leq1$. There exists a maximal $\delta$-separated set $\{v_1, \ldots, v_M\} \subset S^{d-1}$, with associated closed Voronoi cells $V_j := \{v \in S^{d-1} : \dist_g(v, v_j) \leq \dist_g(v, v_k) \text{ for all } k\}$, such that the packing number $M$ satisfies
\[
(d-1)\frac{\omega_d}{\omega_{d-1}}\delta^{-(d-1)} \;\leq\; M \;\leq\; 2(d-1)\pi^{d-2}\frac{\omega_d}{\omega_{d-1}}\delta^{-(d-1)},
\]
and each cell has diameter $\diam(V_j) \leq 2\delta$.
\end{proposition}

\begin{proof}
For the upper bound on $M$, the separation property ensures the geodesic balls $\{B_g(v_j, \delta/2)\}_{j=1}^M$ have pairwise disjoint interiors, so their volumes sum to at most $\omega_d$. The volume of a geodesic ball of radius $r$ on $S^{d-1}$ is
\[
|B_g(v, r)| = \omega_{d-1}\int_0^r\sin^{d-2}(\theta)\,d\theta.
\]
For $\theta \in [0, \pi/2]$, the inequality $\sin\theta \geq (2/\pi)\theta$ gives $\sin^{d-2}(\theta) \geq (2/\pi)^{d-2}\theta^{d-2}$, so
\[
|B_g(v, r)| \geq \omega_{d-1}\left(\frac{2}{\pi}\right)^{d-2}\frac{r^{d-1}}{d-1}.
\]
Applying this with $r = \delta/2$ and using $\sum_{j=1}^M|B_g(v_j,\delta/2)| \leq \omega_d$,
\[
M\cdot\omega_{d-1}\left(\frac{2}{\pi}\right)^{d-2}\frac{(\delta/2)^{d-1}}{d-1} \leq \omega_d,
\]
which gives $M \leq 2(d-1)\pi^{d-2}\omega_d\omega_{d-1}^{-1}\delta^{-(d-1)}$.

For the lower bound, by maximality the balls $\{B_g(v_j,\delta)\}_{j=1}^M$ cover $S^{d-1}$. Since $\sin\theta \leq \theta$, we have $|B_g(v,r)| \leq \omega_{d-1}r^{d-1}/(d-1)$, and the covering property gives $\omega_d \leq M\cdot \omega_{d-1}\delta^{d-1}/(d-1)$, hence $M \geq (d-1)\omega_d\omega_{d-1}^{-1} \delta^{-(d-1)}$.

For the diameter bound, any two points $v, v' \in V_j$ both lie within geodesic distance $\delta$ of $v_j$, so $\dist_g(v,v') \leq \dist_g(v,v_j) + \dist_g(v_j,v') \leq 2\delta$.
\end{proof}

For each packing center $v_j$, we associate a stripe on $S^{d-1}$ defined by proximity to the equator perpendicular to $v_j$.

\begin{definition}[Universal stripes]
\label{def:universal-stripes}
For each packing center $v_j\in S^{d-1}$, define $\Sigma_j := \{v \in S^{d-1} : |\langle v, v_j\rangle| \leq 4\delta\}.$
\end{definition}

\subsubsection{Geometric properties for stripe coverage}\label{app:coverage}

The stripe cover of Definition~\ref{def:universal-stripes} is built around the packing centers, not around the ambiguity set itself; the next proposition supplies the missing link between the two, needed to turn the universal cover into an actual cover of the switching set.

\begin{proposition}[Geometric properties for stripe coverage]
\label{prop:geometric-stripe-coverage}
Fix $0<\delta\leq1$ and a maximal $\delta$-separated set $\{v_1,\ldots,v_M\}\subset S^{d-1}$ as in Proposition~\ref{prop:packing-partition}. For every $U=[u_1,\ldots,u_d]\in\orth(d)$, the following three properties hold.
\begin{enumerate}[label=\arabic*., leftmargin=2.2em, labelsep=0.6em]
  \item \textbf{(Pairwise-bisector structure)} For $\sigma\in\{+1,-1\}$, let $g_{ij}^{\sigma}:=\{v\in S^{d-1}:u_i^\top v=\sigma u_j^\top v\}$. Then
  \[
  \Ac(U)\subseteq\bigcup_{i<j,\,\sigma\in\{+1,-1\}}g_{ij}^{\sigma},
  \]
  a union of $d(d-1)$ geodesic $(d-2)$-spheres.
  \item \textbf{(Normal direction approximation)} For each of the geodesic $(d-2)$-spheres $g_{ij}^\pm$ with unit normal $w \in S^{d-1}$ (i.e., $w^{\top} v =0 $ for all $v\in g_{ij}^\pm$), there exists a packing center $v_l$ from Proposition~\ref{prop:packing-partition} such that $\|w - v_l\|_2 \leq \delta$.
  \item \textbf{(Geometric containment)} Let $w,v_\ell\in S^{d-1}$ satisfy $\|w-v_\ell\|_2\leq\varepsilon$, and let $E_w:=\{v\in S^{d-1}:\langle v,w\rangle=0\}$. Then, for any $a>0$,
  \[
  \tube_a(E_w)\subseteq\{v\in S^{d-1}:|\langle v,v_\ell\rangle|\leq a+\varepsilon\}.
  \]
\end{enumerate}
\end{proposition}

\begin{proof}
Part (1) follows from $|u_i^\top v|=|u_j^\top v|$ iff $(u_i^\top v-u_j^\top v)(u_i^\top v+u_j^\top v)=0$; every active tie is pairwise.

For part (2), by the maximality of the $\delta$-separated packing, every point on $S^{d-1}$ lies within geodesic distance $\delta$ of some center $v_l$. In particular, taking $w$ to be the unit normal of any geodesic $g_{ij}^\pm$, there exists a center $v_l$ with $\dist_g(w, v_l) \leq \delta$. Converting to chordal distance gives $\|w - v_l\|_2 = 2\sin(\dist_g(w,v_l)/2) \leq 2\sin(\delta/2) \leq \delta$, where the last inequality uses $\sin(x) \leq x$ for $x \geq 0$.

For part (3), let $v\in\tube_a(E_w)$ and choose a nearest $y\in E_w$. Since chordal distance is at most geodesic distance,
\[
|\langle v,w\rangle|=|\langle v-y,w\rangle|\leq\|v-y\|_2\leq\dist_g(v,y)\leq a.
\]
Therefore $|\langle v,v_\ell\rangle|\leq|\langle v,w\rangle|+\|v_\ell-w\|_2\leq a+\varepsilon$.
\end{proof}

Combining this with the switching localization of Corollary~\ref{cor:switching-localization} gives the sparse stripe cover of the switching set we have been building toward.

\begin{lemma}[Stripe coverage of the switching region]
\label{lem:stripe-coverage-real}
Let $0<\delta\leq1$. For frames $U_1,U_2\in\orth(d)$ with $\|U_1-U_2\|_F\leq\delta$ and the stripes of Definition~\ref{def:universal-stripes}, there exists $\Ic(U_1,U_2)\subseteq\{1,\ldots,M\}$ with $|\Ic(U_1,U_2)|\leq d(d-1)$ such that
\[
\Sswitch(U_1, U_2) \subseteq \bigcup_{j \in \Ic(U_1,U_2)}\Sigma_j.
\]
Moreover, every center indexed by $\Ic(U_1,U_2)$ can be chosen within Euclidean distance $\delta$ of the corresponding signed-bisector normal of $U_1$.
\end{lemma}

\begin{proof}
By Corollary~\ref{cor:switching-localization} and Proposition~\ref{prop:geometric-stripe-coverage}(1),
\[
\Sswitch(U_1,U_2)\subseteq\bigcup_{i<j,\,\sigma\in\{+1,-1\}}\tube_{\frac{\pi}{\sqrt2}\delta}(g_{ij}^\sigma).
\]
Fix one $g_{ij}^\sigma$, let $w$ be its unit normal, and choose the packing center $v_\ell$ supplied by Proposition~\ref{prop:geometric-stripe-coverage}(2). Part (3) gives
\[
\tube_{\frac{\pi}{\sqrt2}\delta}(g_{ij}^\sigma)
\subseteq\left\{v:|\langle v,v_\ell\rangle|\leq\left(\frac{\pi}{\sqrt2}+1\right)\delta\right\}
\subseteq\{v:|\langle v,v_\ell\rangle|\leq4\delta\}=\Sigma_\ell.
\]
The last inclusion uses $1+\pi/\sqrt2<4$. Taking the corresponding stripe indices over the $d(d-1)$ signed bisectors proves the claim.
\end{proof}

This sparse cover is the input to the sample-count bounds of the next subsection.

\subsection{Simultaneous concentration across stripes}\label{app:concentration-frame}

Since switching is confined to the far region by Lemma~\ref{lem:far-confinement}, we work exclusively with far-region sample counts below.

\begin{lemma}[Switching confinement to the far region]
\label{lem:far-confinement}
Let $c_{\Nc}:=\cos(2r_{\Nc})-2r_{\Nc}$, as in Assumption~\ref{ass:switching-margin}, and suppose $c_{\Nc}>1/\sqrt{2}$. For any $U_1,U_2\in\Nc$ with
$\|U_1-U_2\|_F\leq\delta$ and $2\delta<c_{\Nc}-\sqrt{1-c_{\Nc}^2}$, we have $\Sswitch(U_1,U_2)\subseteq\Rfar(\Nc)$. In particular, for $\delta=n^{-1/2}$ this holds for all $n>4\,(c_{\Nc}-\sqrt{1-c_{\Nc}^2})^{-2}$.
\end{lemma}

\begin{proof}
Suppose $v \in \Rnear(\Nc)$, so $v \in \Capp(u_k^{(U)})$ for some $U \in \Nc$ and index $k$. 
Since both $U$ and $U_1$ lie in the geodesic ball of radius $r_{\Nc}$ around $\Ustar$, they are at most $2r_{\Nc}$ apart, and the standard embedding bound for $\orth(d)\subseteq\RR^{d\times d}$ gives $\|U-U_1\|_\Fr\leq2r_{\Nc}$. Therefore
\[
|\langle v,u_k^{(U_1)}\rangle|\;\geq\;|\langle v,u_k^{(U)}\rangle|-\|u_k^{(U)}-u_k^{(U_1)}\|_2\;\geq\;\cos(2r_{\Nc})-\|U-U_1\|_F\;\geq\;c_{\Nc}.
\]
Since $\{u_\ell^{(U_1)}\}_{\ell=1}^d$ is an orthonormal basis, completeness gives
$\max_{\ell\neq k}|\langle v,u_\ell^{(U_1)}\rangle|\leq\sqrt{1-c_{\Nc}^2}$. The margin of $v$ at axis $k$ under $U_1$ is therefore
\[
\Delta := |\langle v, u_k^{(U_1)}\rangle| - \max_{\ell \neq k}|\langle v, u_\ell^{(U_1)}\rangle| \;\geq\; c_{\Nc} - \sqrt{1-c_{\Nc}^2} > 2\delta,
\]
where the last inequality holds by assumption. For every $\ell$, the absolute projection changes by at most $\|u_\ell^{(U_1)}-u_\ell^{(U_2)}\|_2\leq\delta$. The margin $\Delta>2\delta$ therefore preserves $k$ as the unique winner, so $v\notin\Sswitch(U_1,U_2)$.
\end{proof}
For each stripe $\Sigma_j$, define the far-region sample count
\[
N_j^{\mathrm{far}} := \sum_{i:x_i\neq0} \one\!\left\{\frac{x_i}{\|x_i\|} \in \Sigma_j \cap \Rfar(\Nc)\right\}
\]
and the far-region stripe probability
\[
p_j^{\mathrm{far}} := \PP\!\left(X\neq0,\ \frac{X}{\|X\|} \in \Sigma_j \cap \Rfar(\Nc)\right).
\]
Write $m_j:=\EE[N_j^{\mathrm{far}}]=np_j^{\mathrm{far}}$.

\begin{lemma}[Simultaneous stripe concentration]
\label{lem:simultaneous-far}
Let $x_1,\ldots,x_n$ be i.i.d. samples, and construct the stripes of Definition~\ref{def:universal-stripes} with $\delta=n^{-1/2}$. Then, with probability at least $1-Mn^{-(d-1)}$,
all stripes simultaneously satisfy
\[
N_j^{\mathrm{far}} \leq m_j + \max\!\left(3\sqrt{m_j\,d\log n},\; 9d\log n\right) \qquad\text{for all } j = 1,\ldots,M.
\]
\end{lemma}

\begin{proof}
Write $N_j^{\mathrm{far}}=\sum_{i=1}^n X_i^{(j)}$, where $X_i^{(j)}:=\one\{x_i\neq0,\ x_i/\|x_i\|\in\Sigma_j\cap\Rfar(\Nc)\}$, so $\Var(N_j^{\mathrm{far}})\leq m_j$. Set $t := \max(3\sqrt{m_j\,d\log n},\, 9d\log n)$. Bernstein's inequality gives
\[
\PP(N_j^{\mathrm{far}} \geq m_j + t) \leq \exp\!\left(-\frac{t^2}{2(m_j + t/3)}\right).
\]
In the two cases $m_j\geq9d\log n$ and $m_j<9d\log n$, the exponent is at least $27d\log n/8$, so each stripe fails with probability at most $n^{-3d}$. A union bound over all $M$ stripes gives
\[
\PP(\mathsf E_{\mathrm{str}}^c)\leq M n^{-3d}\leq Mn^{-(d-1)},
\]
where $\mathsf E_{\mathrm{str}}$ is the displayed simultaneous event and the last step uses $n^{-3d}\leq n^{-(d-1)}$ for $n\geq1$ (we keep the simpler $n^{-(d-1)}$ rate rather than the sharper $n^{-3d}$, for comparability with the probability bounds elsewhere in this appendix). This is the stated bound.
\end{proof}

\begin{proof}[Proof of Theorem~\ref{thm:switching-margin}]
Set $\delta:=n^{-1/2}$ and let $\mathsf E_{\mathrm{str}}$ be the simultaneous concentration event of Lemma~\ref{lem:simultaneous-far}. We first bound the sample counts in stripes whose centers approximate admissible bisector normals, then use the geometric covering results to control switching simultaneously over all nearby frame pairs. Define the deterministic eligible index set
\[
\mathcal J_n:=\{j:\dist(v_j,\mathcal W_{\Nc})\leq\delta\}.
\]

\emph{Per-stripe tail bound.}
For $j\in\mathcal J_n$, the sample-size condition ensures $\dist(v_j,\mathcal W_{\Nc})\leq\eta_0$ and $4\delta<\eta_0$. Since $\Sigma_j=\{v:|\langle v,v_j\rangle|\leq4\delta\}$, Assumption~\ref{ass:switching-margin} gives
\[
p_j^{\mathrm{far}}\leq4\Cmar\delta
\qquad\text{and}\qquad m_j\leq4\Cmar\sqrt n.
\]
Applying Lemma~\ref{lem:simultaneous-far} and substituting this bound for $m_j$, we obtain on $\mathsf E_{\mathrm{str}}$, simultaneously for every $j\in\mathcal J_n$,
\begin{align*}
N_j^{\mathrm{far}}
&\leq m_j+3\sqrt{m_jd\log n}+9d\log n\\
&\leq4\Cmar\sqrt n+6\sqrt{\Cmar d\log n}\,n^{1/4}+9d\log n=B_n.
\end{align*}
Thus $\mathsf E_{\mathrm{str}}$ implies $\Efar:=\{\max_{j\in\mathcal J_n}N_j^{\mathrm{far}}\leq B_n\}$.

\emph{From stripe bound to switching bound.}
Fix any $U_1,U_2\in\Nc$ with $\|U_1-U_2\|_F\leq\delta$. Lemma~\ref{lem:far-confinement} confines their switching set to $\Rfar(\Nc)$, while Lemma~\ref{lem:stripe-coverage-real} covers it using an index set $\Ic(U_1,U_2)\subseteq\mathcal J_n$ of size at most $d(d-1)$. Consequently,
\[
\Sswitch(U_1,U_2)
\subseteq\bigcup_{j\in\Ic(U_1,U_2)}(\Sigma_j\cap\Rfar(\Nc)).
\]
On $\Efar$, each stripe in this cover contains at most $B_n$ samples from the far region, so $N_\delta(U_1,U_2)\leq d(d-1)B_n$. Although the covering indices depend on the frame pair, the event $\Efar$ controls every eligible stripe and is therefore independent of that choice. Hence the switching-count bound holds simultaneously for all admissible pairs. Finally, the sharper probability estimate $\PP(\mathsf E_{\mathrm{str}}^c)\leq Mn^{-3d}$ from Lemma~\ref{lem:simultaneous-far}, together with the packing bound of Proposition~\ref{prop:packing-partition} at $\delta=n^{-1/2}$, gives $\PP(\Efar^c)\leq p_{\mathrm{sw}}(n)$, completing the proof.
\end{proof}

\begin{remark}[Other margin exponents]\label{rem:margin-exponent}
Suppose the right side of Assumption~\ref{ass:switching-margin} is instead $\Cmar h^\alpha$ for some $\alpha>0$. The same proof gives
\[
B_{n,\alpha}=4^\alpha\Cmar n^{1-\alpha/2}
+3\sqrt{4^\alpha\Cmar d\log n}\,n^{1/2-\alpha/4}+9d\log n.
\]
Hence the switching term $\mathfrak J_n$ of Proposition~\ref{prop:quasi-lipschitz}, with $B_{n,\alpha}$ in place of $B_n$, vanishes for every $\alpha>0$, and it is of order at most $n^{-1/2}\log n$ exactly when $\alpha\geq1$, which preserves the near-parametric rate. For $\alpha<1$ it is of order $n^{-\alpha/2}\log n$ and dominates the rate. By contrast, for $\alpha=0$, where the stripe probabilities are only bounded by a constant, one gets $B_{n,0}\asymp n$ and $\mathfrak J_n$ does not vanish.
\end{remark}

\subsection{Uniform convergence of empirical subgradients}\label{app:uniform}

We prove Propositions~\ref{prop:quasi-lipschitz} and~\ref{prop:uniform-convergence}, using a deterministic covering lemma to pass from pointwise to uniform control. Throughout, we work on the aligned component $\Nc$; invariance extends the bounds to its signed-permutation translates.

\subsubsection{Coarsely Lipschitz property}\label{app:quasilip}

\begin{proof}[Proof of Proposition~\ref{prop:quasi-lipschitz}]
Fix $U_1,U_2\in\Nc$ and set $\delta:=\|U_1-U_2\|_F\leq n^{-1/2}$. Partition the sample indices into $\Ins:=\{i:j^\star(x_i;U_1)=j^\star(x_i;U_2)\}$ and $\Is:=\{1,\ldots,n\}\setminus\Ins$, according to whether their selected axes remain unchanged. Then
\[
G_n(U_1)-G_n(U_2)=\frac1n\sum_{i\in\Ins}\!\bigl(g_{\mathrm{can}}(x_i,U_1)-g_{\mathrm{can}}(x_i,U_2)\bigr)
+\frac1n\sum_{i\in\Is}\!\bigl(g_{\mathrm{can}}(x_i,U_1)-g_{\mathrm{can}}(x_i,U_2)\bigr).
\]

\emph{Non-switching contribution.}
For $i\in\Ins$, let $j$ be the common selected index, and set $X:=x_ix_i^\top$, $P_k:=u_j^{(k)}u_j^{(k)\top}$, and $\Delta:=U_1-U_2$. By~\eqref{eq:riemannian-gradient-fi}, we have $\gamma_j(x_i,U_k)=[P_k,X]U_k$, where $[A,B]:=AB-BA$. Subtracting these expressions gives
\[
\gamma_j(x_i,U_1)-\gamma_j(x_i,U_2)
=[P_1,X]\Delta+[P_1-P_2,X]U_2.
\]
Using $\|[P_1,X]\|_2\leq\|x_i\|_2^2/2$, $\|[P_1-P_2,X]\|_F\leq\sqrt2\|x_i\|_2^2\|P_1-P_2\|_2$, and $\|P_1-P_2\|_2\leq\|\Delta\|_F$, we obtain
\[
\|g_{\mathrm{can}}(x_i,U_1)-g_{\mathrm{can}}(x_i,U_2)\|_F
\leq\left(\tfrac12+\sqrt2\right)\|x_i\|_2^2\delta
\leq2\|x_i\|_2^2\delta.
\]
This estimate also applies when either endpoint is ambiguous, provided the canonical selected indices agree. To control the sum, set $Z:=\|x\|_2^2$. Assumption~\ref{ass:subexp-norms} and Jensen's inequality give $\mtwo\leq K\log2$, while the triangle inequality and Cauchy--Schwarz yield
\[
\EE\exp\!\left(\frac{|Z-\mtwo|}{2K}\right)
\leq e^{\mtwo/(2K)}\EE e^{Z/(2K)}\leq2.
\]
Thus $\|Z-\mtwo\|_{\psi_1}\leq2K$, and scalar Bernstein concentration gives an absolute constant $C_B>0$ such that
\[
\PP\!\left(\frac1n\sum_i\|x_i\|_2^2>2\mtwo\right)
\leq2\exp\!\left(-C_B n\frac{\mtwo^2}{K^2}\right)
\leq\frac1{4n},
\]
where the last inequality follows from the sample-size condition in the proposition. On the complementary event, summing the preceding pointwise estimate over $i\in\Ins$ bounds the non-switching contribution by $4\mtwo\delta$.

\emph{Switching contribution.}
For $i\in\Is$, Proposition~\ref{prop:individual-clarke-bounds} gives $\|g_{\mathrm{can}}(x_i,U_1)-g_{\mathrm{can}}(x_i,U_2)\|_F\leq4\|x_i\|_2^2$. Unlike the non-switching case, this bound need not vanish with $\delta$, so we instead control the number of switching indices. On the maximum-norm event of Lemma~\ref{lem:max-norm} and the simultaneous switching-count event of Theorem~\ref{thm:switching-margin},
\[
\left\|\frac1n\sum_{i\in\Is}\bigl(g_{\mathrm{can}}(x_i,U_1)-g_{\mathrm{can}}(x_i,U_2)\bigr)\right\|_F
\leq\frac{4|\Is|}{n}\max_i\|x_i\|_2^2
\leq\frac{12K\,d(d-1)B_n\log n}{n}=\mathfrak J_n.
\]
The simultaneous event makes this estimate valid for every admissible pair $U_1,U_2$, including data-dependent pairs. Combining the two contributions therefore gives the asserted coarsely-Lipschitz bound. Finally, the average-norm, maximum-norm, and switching-count events fail with probabilities at most $1/(4n)$, $1/n$, and $p_{\mathrm{sw}}(n)$, respectively. A union bound yields total failure probability at most $5/(4n)+p_{\mathrm{sw}}(n)\leq2/n+p_{\mathrm{sw}}(n)$, as claimed.
\end{proof}

\subsubsection{Deterministic covering lemma}\label{app:covering}

\begin{lemma}[Deterministic covering argument]
\label{lem:epsilon-covering}
Suppose $\Rgrad\Ec(\Ustar)=0$, the population gradient is intrinsically $L$-Lipschitz on $\Nc$, and the empirical subgradient satisfies 
\[\|G_n(U_1)-G_n(U_2)\|_{\Fr}\leq L_0d_g(U_1,U_2)+J_n\] 
for all $U_1,U_2\in\Nc$ with $d_g(U_1,U_2)\leq\epsilon$. If $\{U_1,\ldots,U_N\}$ is an $\epsilon$-net of $\Nc$, then
\begin{align*}
    \sup_{U\in\Nc}\|\Rgrad\Ec(U)-G_n(U)\|_{\Fr}
    &\leq\max_{1\leq i\leq N}\|\Rgrad\Ec(U_i)-G_n(U_i)\|_{\Fr}
    +\bigl(L(1+r_{\Nc})+L_0\bigr)\epsilon+J_n.
\end{align*}
\end{lemma}

\begin{proof}
For $U,V\in\Nc$, write $V=U\exp(A)$ along their unique minimizing geodesic, where $A^\top=-A$ and $d_g(U,V)=\|A\|_{\Fr}$.

\emph{Ambient Lipschitz bound.}
If $Z=VB\in T_V\orth(d)$, then
\[
\mathsf P_{V\to U}Z=U\exp(A/2)B\exp(-A/2)
\qquad\text{and}\qquad
\|\mathsf P_{V\to U}Z-Z\|_{\Fr}
\leq d_g(U,V)\|Z\|_{\Fr}.
\]
Indeed, the second estimate follows from
$\|\exp(\pm A/2)-I\|_{\opp}\leq\|A\|_{\opp}/2$.
Since $\Rgrad\Ec(\Ustar)=0$, intrinsic $L$-Lipschitz continuity gives
$\|\Rgrad\Ec(V)\|_{\Fr}\leq Ld_g(V,\Ustar)\leq Lr_{\Nc}$. Consequently,
\[
\|\Rgrad\Ec(U)-\Rgrad\Ec(V)\|_{\Fr}
\leq L(1+r_{\Nc})d_g(U,V).
\]

\emph{From net point to global bound.}
For any $U \in \Nc$, let $U_i$ be a net point satisfying $d_g(U, U_i) \leq \epsilon$. Inserting $\pm G_n(U_i)$ and $\pm \Rgrad\Ec(U_i)$ and applying the
triangle inequality,
\begin{align*}
\|\Rgrad\Ec(U) - G_n(U)\|_{\Fr}
&\leq \|\Rgrad\Ec(U) - \Rgrad\Ec(U_i)\|_{\Fr}\\
&\qquad + \|\Rgrad\Ec(U_i) - G_n(U_i)\|_{\Fr}\\
&\qquad + \|G_n(U_i) - G_n(U)\|_{\Fr}.
\end{align*}
The first term is at most $L(1+r_{\Nc})\epsilon$ by the ambient comparison above, and the third is at most $L_0\epsilon+J_n$. Therefore
\[
\|\Rgrad\Ec(U)-G_n(U)\|_{\Fr}
\leq\|\Rgrad\Ec(U_i)-G_n(U_i)\|_{\Fr}+\bigl(L(1+r_{\Nc})+L_0\bigr)\epsilon+J_n.
\]
Since the bound holds for every $U\in\Nc$, taking the supremum over $U$ proves the lemma.
\end{proof}

\subsubsection{Proof of uniform empirical subgradient control}\label{app:uniform-proof}

\begin{proof}[Proof of Proposition~\ref{prop:uniform-convergence}]
We combine pointwise concentration on a finite net with the coarsely-Lipschitz estimate of Proposition~\ref{prop:quasi-lipschitz} to obtain uniform control of the canonical empirical selection.

Set $D:=d(d-1)/2$, $\epsilon:=n^{-1/2}$, and $r:=r_{\Nc}$, so that $\Lambda_d=2+D/2+D\log(1+4r)$ as in Theorem~\ref{thm:convergence}. By Assumption~\ref{ass:strong-convexity}, every $U\in\Nc$ can be written as $U=\Exp_{\Ustar}(v)$, where $v:=\operatorname{Log}_{\Ustar}(U)$ and $\|v\|_{\Fr}=d_g(\Ustar,U)\leq r$. Cover the Euclidean $D$-ball of radius $r$ by balls of radius $\epsilon/2$, using at most $(1+4r/\epsilon)^D$ centers, and retain one point from each ball meeting $\operatorname{Log}_{\Ustar}(\Nc)$. By Lemma~\ref{lem:orth-exp-nonexpansive}, their exponential images form an $\epsilon$-net $\{U_1,\ldots,U_N\}\subset\Nc$. Since $\epsilon\leq1$,
\[
N\leq(1+4r/\epsilon)^D\leq(1+4r)^D\epsilon^{-D}=(1+4r)^D n^{D/2}.
\]

We next control the empirical gradient at every net point by applying Proposition~\ref{prop:pointwise-bound} to each deterministic $U_i\in\Nc$. Set $t:=4\Lambda_d\log n$, and let $\Enet$ be the event on which its high-probability estimate holds simultaneously at all net points. Since the sample-size conditions give $n>64$, a union bound and the covering-number estimate yield
\[
\PP(\Enet^c)\leq2Nn^{-\Lambda_d}
\leq2(1+4r)^Dn^{D/2-\Lambda_d}
=2n^{-2}(1+4r)^Dn^{-D\log(1+4r)}\leq\frac1n,
\]
while on $\Enet$,
\[
\max_{1\leq i\leq N}\|\Rgrad\Ec(U_i)-G_n(U_i)\|_{\Fr}
\leq2\left (\sqrt{8\mfour}+\sqrt2\,\Gamma\right )\sqrt{\Lambda_d}\sqrt{\frac{\log n}{n}}
+4\Gamma\Lambda_d\frac{\log n}{n}.
\]

To extend this estimate beyond the net, define
\[
\Ec_{\mathrm{av}}:=\left\{n^{-1}\sum_i\|x_i\|_2^2\leq2\mtwo\right\},
\]
\[
\Ec_{\max}:=\left\{\max_i\|x_i\|_2^2\leq3K\log n\right\},
\]
and
\[
\Efar:=\left\{\max_{j\in\mathcal J_n}N_j^{\mathrm{far}}\leq B_n\right\},
\]
and set $\Equasi:=\Ec_{\mathrm{av}}\cap\Ec_{\max}\cap\Efar$. The concentration estimates used in Proposition~\ref{prop:quasi-lipschitz} give
\[
\PP(\Equasi^c)\leq\frac{5}{4n}+p_{\mathrm{sw}}(n)\leq\frac2n+p_{\mathrm{sw}}(n).
\]
For any $U\in\Nc$, choose a net point $U_i$ with $d_g(U,U_i)\leq\epsilon$. Remark~\ref{rem:frob-geod} gives $\|U-U_i\|_F\leq d_g(U,U_i)$, so Proposition~\ref{prop:quasi-lipschitz} yields on $\Equasi$
\[
\|G_n(U)-G_n(U_i)\|_F\leq4\mtwo d_g(U,U_i)+\mathfrak J_n.
\]
Thus, on $\Egood:=\Enet\cap\Equasi$, the empirical selection is controlled between each frame and its nearest net point, while $\Enet$ controls its deviation from the population gradient at that net point. Since $\PP(\Egood)\geq1-3/n-p_{\mathrm{sw}}(n)$, Lemma~\ref{lem:epsilon-covering} combines these estimates to give
\begin{align*}
\sup_{U\in\Nc}\|\Rgrad\Ec(U)-G_n(U)\|_{\Fr}
&\leq2\left (\sqrt{8\mfour}+\sqrt2\,\Gamma\right )\sqrt{\Lambda_d}\sqrt{\frac{\log n}{n}}\\
&\qquad+4\Gamma\Lambda_d\frac{\log n}{n}+\bigl(L(1+r_{\Nc})+4\mtwo\bigr)n^{-1/2}+\mathfrak J_n.
\end{align*}
Using $n^{-1/2}\leq\sqrt{\log n/n}$ gives the bound in the proposition.
\end{proof}

\subsection{Localization of the empirical minimizer}\label{app:localization}

We prove Corollary~\ref{cor:localization} in two parts: Lemma~\ref{lem:deterministic-localization} shows deterministically that uniform objective control implies localization, and Proposition~\ref{prop:global-uniform-convergence} shows this uniform control holds with high probability.

\begin{lemma}[Deterministic localization via uniform deviation]
\label{lem:deterministic-localization}
Define the event
\[
\mathsf{A}_\rho := \left\{\sup_{U \in \orth(d)}|\Ec_n(U) - \Ec(U)| \leq \frac{\rho}{4}\right\}.
\]
On $\mathsf{A}_\rho$, every empirical minimizer $\Ustar_n\in\argmin_{U\in\orth(d)}\Ec_n(U)$ belongs to $\Srho\subseteq\Norbit$. Consequently there is a unique $R_n\in\Spm$ such that $\Ustarntilde:=\Ustar_nR_n\in\Nc$. By invariance, $\Ustarntilde$ is also an empirical minimizer and $d_g(\Ustarntilde,\Ustar)=\dframe(\Ustar_n,\Ustar)$.
\end{lemma}

\begin{proof}
Let $U\notin\Srho$. On $\mathsf A_\rho$,
the uniform deviation bound and the definition of $\Srho$ give
\[
\Ec_n(U)-\Ec_n(\Ustar)
\geq \Ec(U)-\Ec(\Ustar)-\frac{\rho}{2}
\geq \frac{\rho}{2}>0,
\]
so $U$ cannot minimize the empirical risk. Hence every empirical minimizer lies in $\Srho\subseteq\Norbit$. The pairwise-disjoint decomposition $\Norbit=\bigsqcup_{R\in\Spm}\Nc R$ supplies the unique alignment described in the lemma.
\end{proof}

\begin{proposition}[Global uniform objective convergence]
\label{prop:global-uniform-convergence}
Under Assumption~\ref{ass:subexp-norms}, fix $\rho>0$ and $\beta\in(0,1)$. If $K=0$, then $\PP(\mathsf A_\rho)=1$. Otherwise set $\rho_0:=\min\{\rho,K\}$. There is an absolute constant $\Cstat>1$ such that, if
\[
n\geq\Cstat\frac{K^2}{\rho_0^2}
\left[d^2\log\!\left(\frac{\Cstat\sqrt d\,K}{\rho_0}\right)+\log\!\left(\frac{\Cstat}{\beta}\right)\right],
\]
then $\PP(\mathsf A_\rho)\geq1-\beta$.
\end{proposition}

\begin{proof}
The argument is a standard covering argument: we first control the sample-average norm and the resulting Lipschitz constants of $\Ec_n$ and $\Ec$, then bound the deviation $|\Ec_n(U')-\Ec(U')|$ at finitely many net points, and finally propagate this control to all of $\orth(d)$ via the Lipschitz bound.

If $K=0$, then $\|X\|=0$ almost surely and the claim is immediate. Assume $K>0$, put $Z:=\|X\|_2^2$, and set $\rho_0:=\min\{\rho,K\}$. Since $0\leq f_X(U)\leq Z$, standard $\psi_1$ properties give, uniformly in $U$,
\[
\mtwo\leq C_\psi K,\qquad
\|Z-\mtwo\|_{\psi_1}\leq C_\psi K,\qquad\text{and}\qquad
\|f_X(U)-\Ec(U)\|_{\psi_1}\leq C_\psi K
\]
for an absolute $C_\psi$. Scalar Bernstein concentration therefore yields absolute $c,C_L>0$ such that the event
\[
\mathsf E_L:=\left\{\frac1n\sum_{i=1}^n\|x_i\|_2^2\leq C_LK\right\}
\]
satisfies $\PP(\mathsf E_L^c)\leq2e^{-cn}$. On $\mathsf E_L$, Lemma~\ref{lem:lipschitz-empirical-obj} shows that the sum of the empirical and population objective Lipschitz constants is at most $C_L'K$ for an absolute $C_L'>0$. Choose $\varepsilon:=\rho_0/(8C_L'K)$ and let $\Tc\subset\orth(d)$ be a maximal $\varepsilon$-separated set in Frobenius norm. It is an $\varepsilon$-net. Since $\orth(d)$ lies in the radius-$\sqrt d$ ball of $\RR^{d^2}$, ambient volume comparison gives
\[
|\Tc|\leq\left(1+\frac{2\sqrt d}{\varepsilon}\right)^{d^2}
\leq\left(\frac{C_{\mathrm{cov}}\sqrt d\,K}{\rho_0}\right)^{d^2}
\]
for an absolute $C_{\mathrm{cov}}>1$. At each $U'\in\Tc$, scalar Bernstein and the preceding $\psi_1$ bound imply
\[
\PP\!\left(|\Ec_n(U')-\Ec(U')|>\frac{\rho_0}{8}\right)
\leq2\exp\!\left[-cn\min\!\left\{\frac{\rho_0^2}{K^2},\frac{\rho_0}{K}\right\}\right].
\]
Choosing $\Cstat$ larger than the preceding absolute constants, the stated sample size and a union bound give
\[
\PP\!\left(\mathsf E_L\ \text{and}
\max_{U'\in\Tc}|\Ec_n(U')-\Ec(U')|\leq\frac{\rho_0}{8}\right)\geq1-\beta.
\]
On this event, invariance and the net approximation yield
\[
\sup_{U\in\orth(d)}|\Ec_n(U)-\Ec(U)|
\leq\frac{\rho_0}{8}+C_L'K\varepsilon
=\frac{\rho_0}{4}\leq\frac{\rho}{4}.
\]
This is precisely $\mathsf A_\rho$.
\end{proof}

\begin{proof}[Proof of Corollary~\ref{cor:localization}]
On the event $\mathsf A_\rho$ from Proposition~\ref{prop:global-uniform-convergence}, Lemma~\ref{lem:deterministic-localization} gives $\Ustar_n\in\Srho\subseteq\Norbit$ and a unique aligned minimizer $\Ustarntilde\in\Nc$.
\end{proof}

\section{Technical proofs for Theorem~\ref{thm:SGD_convergence_thm_prob}}\label{app:SGD}
We analyze the SGD scheme from Section~\ref{sec:main_results}. 

\subsection{SGD estimator and centered-ball setup}\label{app:sgd-setup}

The estimators $h$ and $g$ from \eqref{grad_estimator} are jointly Borel, by the finite family of continuous scores $|\langle u_j,x\rangle|^2$. Write $g_\infty(U):=\EE_z[g(U,z)]$ for the population mean of $g$, which exists for every $U\in\orth(d)$ under Assumption~\ref{ass:concentration} alone.

\begin{lemma}[Mean oracle, local identification, and equivariance]
\label{lem:sgd-unbiased}
Under Assumption~\ref{ass:concentration}, the mean oracle $g_\infty$ satisfies
\[
\EE\|g(U,z)-g_\infty(U)\|_{\Fr}^2\leq\frac{4\mfour}{B}\qquad\text{for every }U\in\orth(d).
\]
Moreover, if $\F_t:=\sigma(U_0,z_0,\ldots,z_{t-1})$, then $\EE[g(U_t,z_t)\mid\F_t]=g_\infty(U_t)$ for every $t$, regardless of whether $U_t\in\Nc$.

Under Assumption~\ref{ass:strong-convexity}, the mean oracle agrees with the population gradient at every $U\in\Nc$, 
\[
g_\infty(U)=\Rgrad\Ec(U).
\]
Furthermore, for every fixed $U\in\orth(d)$, every $R\in\Spm$, and every batch $z$, we have $g(UR,z)=g(U,z)R$, and hence $g_\infty(UR)=g_\infty(U)R$.
\end{lemma}

\begin{proof}
Write $\zeta(U,x):=\Pi_Uh(x,U)$, so that $g(U,z)=B^{-1}\sum_{i=1}^B\zeta(U,x_i)$ for i.i.d.\ samples $x_i\sim\mu$. The centered summands $\zeta(U,x_i)-g_\infty(U)$ are independent and mean zero. Moreover, nonexpansiveness of the tangent projection gives $\|\zeta(U,x)\|_{\Fr}\leq2\|x\|^2$. Consequently,
\[
\EE\|g(U,z)-g_\infty(U)\|_{\Fr}^2\leq B^{-1}\EE\|\zeta(U,X)\|_{\Fr}^2\leq4\mfour/B,
\]
using only the finite fourth moment in Assumption~\ref{ass:concentration}.

Fresh-batch independence gives $\EE[g(U_t,z_t)\mid\F_t]=g_\infty(U_t)$, since $U_t$ is $\F_t$-measurable and $z_t$ is independent of $\F_t$. For $U\in\Nc$, the quantity $\zeta(U,x)$ is the average of the active Riemannian branch gradients $\gamma_j(x,U)$ over $j\in J(x;U)$. Lemma~\ref{lem:diff-implies-tie-null} shows that, almost surely, this average coincides with the canonical selection $\gamma_{j^\star(x;U)}(x,U)$. The same lemma therefore gives $g_\infty(U)=\EE[\zeta(U,X)]=\Rgrad\Ec(U)$.

Finally, a signed permutation $R$ relabels $J(x;U)$ and changes the signs of the corresponding columns. Since $h$ averages over the entire active set, $h(x,UR)=h(x,U)R$ for every $x$. The projection identity $\Pi_{UR}(AR)=\Pi_U(A)R$ then gives $g(UR,z)=g(U,z)R$ for every $U$, $R$ and batch $z$, and taking expectations proves $g_\infty(UR)=g_\infty(U)R$.
\end{proof}

\begin{remark}[Centered geodesic balls]
\label{rem:injectivity}
Because $\Nc$ is an open neighborhood of $\Ustar$, there exists $r_c>0$ such that the closed geodesic ball
\[
\Bc(\Ustar,r_c)
:=\{U\in\orth(d):d_g(U,\Ustar)\le r_c\}
\]
is contained in $\Nc$. Shrinking $r_c$ if necessary, every pair of points in this closed ball is joined by a unique minimizing geodesic contained in the ball, and the ball is disjoint from all its nontrivial signed-permutation translates.
\end{remark}
\subsection{Dimension-free stochastic-gradient concentration}\label{app:sgd-concentration}

\begin{proof}[Proof of Lemma~\ref{lemma:SG_bound}]

Conditionally on $\F_t$, the iterate $U_t$ is fixed. Writing $\zeta(U,x):=\Pi_Uh(x,U)$, we have
\[
e_t=\frac1B\sum_{i=1}^B \bigl(\zeta(U_t,x_{t,i})-g_\infty(U_t)\bigr).
\]
The summands are conditionally independent and mean zero. Moreover, nonexpansiveness of the tangent projection and the bound on $h$ give
\(
\|\zeta(U,x)\|_{\Fr}\leq2\|x\|^2.
\)
Consequently,
\[
\EE[\|e_t\|_{\Fr}^2\mid\F_t]=\frac1B\EE\!\left[\|\zeta(U_t,X)-g_\infty(U_t)\|_{\Fr}^2\,\middle|\,\F_t\right]\leq\frac1B\EE[\|\zeta(U_t,X)\|_{\Fr}^2\mid\F_t]\leq\frac{4\mfour}{B},
\]
which proves~\eqref{eq:sg-second-moment}. Applying conditional Markov inequalities at each iteration and taking a union bound over $t<T$ then proves~\eqref{eq:eps_B_delta_poly}.

Under Assumption~\ref{ass:subexp-norms}, the stronger tail estimate follows from Hilbert-space concentration. Indeed, $\|\,\|X\|^2\,\|_{\psi_1}=K$, and
\[
\|g_\infty(U_t)\|_{\Fr}
\leq2\EE\|X\|^2
\leq2K\log2.
\]
Together with $\|\zeta(U_t,X)\|_{\Fr}\leq2\|X\|^2$, this gives, conditionally on $\F_t$, the dimension-free bound
\[
\left\|\,
\|\zeta(U_t,X)-g_\infty(U_t)\|_{\Fr}
\,\right\|_{\psi_1}\leq CK.
\]
The standard moment implication for the $\psi_1$ norm therefore verifies the hypotheses of Lemma~\ref{lem:pinelis-hilbert-bennett}, conditionally on $\F_t$. Applying Lemma~\ref{lem:pinelis-hilbert-bennett} to the $B$ centered summands and reparameterizing the resulting tail bound yields, for some absolute constant $C>0$,
\begin{equation}\label{eq:sg_bernstein}
\PP\!\left(
\|e_t\|_{\Fr}\geq
CK\left(\sqrt{\frac{u}{B}}+\frac{u}{B}\right)
\,\middle|\,\F_t\right)
\leq2e^{-u}
\qquad\text{almost surely}.
\end{equation}
Taking $u=\log(2T/\delta)$, averaging over $\F_t$, and applying a union bound over $t<T$ gives failure probability at most $2Te^{-u}=\delta$, proving~\eqref{eq:eps_B_delta_subexp}.
\end{proof}

\subsection{Sharp local contraction and invariant-ball control}\label{app:sgd-contraction}

\begin{lemma}[Integrated-Hessian representation]
\label{lem:sgd-integrated-hessian}
Suppose Assumption~\ref{ass:strong-convexity} holds and
$U\in\Bc(\Ustar,r_c)$. Let
$v:=\operatorname{Log}_U(\Ustar)$ and
$\gamma(s):=\Exp_U(sv)$. If $\mathsf P$ denotes parallel transport along
$\gamma$, define
\[
H_s:=
\mathsf P_{\gamma(s)\to U}
\circ\operatorname{Hess}\Ec(\gamma(s))
\circ\mathsf P_{U\to\gamma(s)}
\qquad\text{and}\qquad
\overline H_U:=\int_0^1H_s\,ds.
\]
Then $\overline H_U$ is self-adjoint,
$mI\preceq\overline H_U\preceq LI$, and
\begin{equation}\label{eq:integrated-hessian-gradient}
\Rgrad\Ec(U)=-\overline H_Uv.
\end{equation}
\end{lemma}

\begin{proof}
The geodesic $\gamma$ lies in $\Bc(\Ustar,r_c)\subseteq\Nc$, where Assumption~\ref{ass:strong-convexity} gives $mI\preceq\operatorname{Hess}\Ec\preceq LI$. Since parallel transport preserves the Riemannian metric, these bounds also hold for each $H_s$ and hence for their average $\overline H_U$. Moreover, differentiating the parallel-transported gradient along $\gamma$ and using $\Rgrad\Ec(\Ustar)=0$ yields
\[
-\Rgrad\Ec(U)
=\int_0^1H_sv\,ds
=\overline H_Uv,
\]
which proves~\eqref{eq:integrated-hessian-gradient}.
\end{proof}

\begin{proof}[Proof of Lemma~\ref{lem:sgd-sharp-contraction}]
For $U_t\in\Bc(\Ustar,r_c)$, Lemma~\ref{lem:sgd-unbiased} gives $g(U_t,z_t)=\Rgrad\Ec(U_t)+e_t$. Set $v_t:=\operatorname{Log}_{U_t}(\Ustar)$ and apply Lemma~\ref{lem:sgd-integrated-hessian} to write $\Rgrad\Ec(U_t)=-\overline H_{U_t}v_t$. Since $\Ustar=\Exp_{U_t}(v_t)$, the exponential-map comparison of Lemma~\ref{lem:orth-exp-nonexpansive} and the spectral bounds on $\overline H_{U_t}$ give
\[
r_{t+1}
\leq\|(I-\eta\overline H_{U_t})v_t+\eta e_t\|_{\Fr}
\leq q_\eta r_t+\eta\|e_t\|_{\Fr},
\]
which proves~\eqref{eq:sharp-linear-recursion}. The stated expressions for $q_\eta$ follow by maximizing $|1-\eta\lambda|$ over $\lambda\in[m,L]$.
\end{proof}

For the expected-error analysis, we also need a conditional squared version of the contraction inequality. Whenever $U_t\in\Bc(\Ustar,r_c)$, squaring the first comparison in the proof of Lemma~\ref{lem:sgd-sharp-contraction} and conditioning on $\F_t$ gives
\begin{equation}\label{eq:sharp-squared-recursion}
\EE[r_{t+1}^2\mid\F_t]
\leq q_\eta^2r_t^2
+\eta^2\EE[\|e_t\|_{\Fr}^2\mid\F_t].
\end{equation}
Indeed, the cross term vanishes because $v_t$ and $\overline H_{U_t}$ are $\F_t$-measurable and $\EE[e_t\mid\F_t]=0$ by Lemma~\ref{lem:sgd-unbiased}. The remaining deterministic term is bounded by $q_\eta^2r_t^2$ using the spectral bounds on $\overline H_{U_t}$.

\begin{proof}[Proof of Lemma~\ref{lem:sgd-invariant-ball}]
Suppose inductively that $U_t\in\Bc(\Ustar,r_c)$. Under the assumed bound $\|e_t\|_{\Fr}\leq\varepsilon$, Lemma~\ref{lem:sgd-sharp-contraction} and $1-q_\eta=\eta m$ give
\[
r_{t+1}\leq q_\eta r_t+\eta\varepsilon.
\]
Iterating this inequality yields~\eqref{eq:invariant-ball-recursion}. Its right-hand side is a convex combination of $r_{\mathrm{init}}$ and $\varepsilon/m$, both at most $r_c$. Hence $r_{t+1}\leq r_c$, closing the induction and proving invariance of the ball.
\end{proof}

\section{An in-expectation online-SGD convergence guarantee}\label{app:sgd-expected}

The preceding analysis establishes high-probability convergence by keeping the iterates inside a strongly convex neighborhood. To obtain an unconditional expected-error bound, we must also account for trajectories that leave this neighborhood. We first control their probability through the stochastic-gradient noise, then combine local contraction with a global diameter bound to control their contribution to the expected squared error.

\subsection{Escape probability}\label{app:escape}

Choose $r_c>0$ as in Remark~\ref{rem:injectivity}, align the deterministic initialization so that $r_{\mathrm{init}}:=d_g(U_0,\Ustar)<r_c$, and define the escape event
\[
\mathcal E_{\mathrm{esc}}
:=\bigcup_{t=1}^{T}\{U_t\notin\Nc\}.
\]

\begin{theorem}\label{thm:bad_evt_prob}
    Under Assumption~\ref{ass:strong-convexity}, consider the matrix iterates $U_t$ from~\eqref{SGD_update} with the fresh independent batches specified above and the centered-ball initialization. Run for $T$ iterations with step size $0<\eta\leq2/(m+L)$ and per-iteration batch size $B$. The following bounds hold for the escape probability $\PP\!\left(\bigcup_{t=1}^{T}\{U_t\notin\Nc\}\right)$.
    \begin{enumerate}[label=\arabic*., leftmargin=2.2em, labelsep=0.6em]
        \item (\textbf{Finite fourth moment}) Under Assumption~\ref{ass:concentration},
        \begin{equation}\label{eq:escape_poly}
            \PP\left(\bigcup_{t = 1}^{T}\{U_t\notin\Nc\}\right)
            \leq \xi_{\poly}(B,T)
            :=\frac{4\mfour T}{Bm^2r_c^2}.
        \end{equation}
        \item (\textbf{Sub-exponential squared norms}) Under Assumptions~\ref{ass:concentration} and~\ref{ass:subexp-norms},
        there exists an absolute constant $c>0$ such that
        \begin{equation}\label{eq:escape_subexp}
            \PP\left(\bigcup_{t = 1}^{T}\{U_t\notin\Nc\}\right)
            \leq \xi_{\mathrm{subexp}}(B,T)
            :=2T\exp\!\left[-cB\min\left\{
            \frac{m^2r_c^2}{K^2},
            \frac{mr_c}{K}\right\}\right].
        \end{equation}
    \end{enumerate}
\end{theorem}

\begin{proof}[Proof of Theorem~\ref{thm:bad_evt_prob}]
Define the prefix good-noise events
\[
G_s:=\bigcap_{t=0}^{s-1}
\{\|e_t\|_{\Fr}\leq mr_c\}
\quad\text{for}\quad s=0,1,\ldots,T,
\quad\text{with}\quad G_0:=\Omega.
\]
By Lemma~\ref{lem:sgd-invariant-ball}, on $G_s$ every iterate $U_0,\ldots,U_s$ remains in $\Bc(\Ustar,r_c)\subseteq\Nc$. Consequently, escape before time $T$ requires at least one noise realization to exceed $mr_c$, so $\mathcal E_{\mathrm{esc}}\subseteq G_T^c$. We therefore bound the probability of $G_T^c$ using the stochastic-gradient concentration estimates, which hold regardless of whether the iterates lie in $\Nc$.

Under Assumption~\ref{ass:concentration}, the conditional second-moment estimate~\eqref{eq:sg-second-moment} and Markov's inequality give a bound for each event $\{\|e_t\|_{\Fr}>mr_c\}$. Summing these bounds over $t<T$ yields
\[
\PP(G_T^c)
\leq\frac{4\mfour T}{Bm^2r_c^2},
\]
which proves~\eqref{eq:escape_poly}.

Under Assumption~\ref{ass:subexp-norms}, inverting the conditional Bernstein bound~\eqref{eq:sg_bernstein} at threshold $mr_c$, averaging over $\F_t$, and taking a union bound over $t<T$ yields, for some absolute constant $c>0$,
\[
\PP(G_T^c)
\leq2T\exp\!\left[-cB\min\left\{
\frac{m^2r_c^2}{K^2},
\frac{mr_c}{K}\right\}\right].
\]
This proves~\eqref{eq:escape_subexp}.
\end{proof}

\subsection{Expected-error bound}
\begin{theorem}[Expected online-SGD convergence]\label{thm:SGD_convergence_thm}
    Suppose Assumptions~\ref{ass:concentration} and~\ref{ass:strong-convexity} hold. Choose $r_c>0$ such that the closed geodesic ball
    \(
    \Bc(\Ustar,r_c):=\{U\in\orth(d)\colon d_g(U,\Ustar)\leq r_c\}
    \)
    is contained in $\Nc$. Let $r_{\mathrm{init}}:=d_g(U_0,\Ustar)<r_c$, and define
    \(
    \kappa:=\frac{L}{m}\),
    $\eta:=\frac{2}{m+L}$, and \(
    q:=\frac{\kappa-1}{\kappa+1}\),
    so that $1-q^2=4mL/(m+L)^2$. For a fixed $c_{\mathrm{iter}}\geq1/2$ and sample budget $n\geq e^2$, take
    \(
    T:=\left\lceil c_{\mathrm{iter}}\kappa\log n\right\rceil\) and $B:=\left\lfloor\frac{n}{T}\right\rfloor$.
    Assume $n\geq2T$. Then
    \begin{equation}\label{eq:sgd-expected}
        \EE[\dframe(U_T,\Ustar)^2] \leq \frac{r_{\mathrm{init}}^2}{n} +\frac{8\mfour}{(1-q^2)\,n}\left(\eta^2T+\frac{\pi^2dT^2}{m^2r_c^2}\right).
    \end{equation}
\end{theorem}

\begin{proof}
Let $G_s$ be the prefix good-noise events from the proof of Theorem~\ref{thm:bad_evt_prob}, and write $r_s:=d_g(U_s,\Ustar)$. We first bound the expected squared distance by combining the local contraction on $G_{s-1}$ with a uniform diameter bound on its complement, and then substitute the prescribed iteration horizon and batch size.

On $G_{s-1}$, Lemma~\ref{lem:sgd-invariant-ball} gives $U_{s-1}\in\Bc(\Ustar,r_c)$, so the conditional squared recursion~\eqref{eq:sharp-squared-recursion} and the noise bound~\eqref{eq:sg-second-moment} apply. Since $G_{s-1}$ is $\F_{s-1}$-measurable, we may multiply the conditional estimate by its indicator before taking expectations. On $G_{s-1}^c$, both $U_s$ and $\Ustar$ belong to the connected component containing $U_0$, since Riemannian exponential updates preserve connected components. The squared geodesic diameter of each component of $\orth(d)$ is at most $\pi^2d$, as follows from the principal rotation-angle representation. Splitting the expectation over these two events and bounding the contribution on the complement therefore gives
\[
\EE r_s^2
\leq q^2\EE r_{s-1}^2
+\frac{4\eta^2\mfour}{B}
+\pi^2d\,\PP(G_{s-1}^c).
\]
The finite-moment estimate from the proof of Theorem~\ref{thm:bad_evt_prob} controls the remaining probability:
\[
\PP(G_{s-1}^c)\leq\PP(G_T^c)
\leq\frac{4\mfour T}{Bm^2r_c^2}.
\]
Substituting this bound into the preceding recursion and summing the resulting geometric series, whose ratio is $q^2$ so that $\sum_{k\ge0}q^{2k}\le(1-q^2)^{-1}$, yields
\[
\EE r_T^2
\leq q^{2T}r_{\mathrm{init}}^2
+\frac{4\mfour}{(1-q^2)B}
\left(\eta^2+\frac{\pi^2dT}{m^2r_c^2}\right).
\]

It remains to express this estimate in terms of the total sample budget $n$. By the definition of $q$,
\[
q=\frac{\kappa-1}{\kappa+1}
\leq\exp\!\left(-\frac{2}{\kappa+1}\right).
\]
Since $T=\lceil c_{\mathrm{iter}}\kappa\log n\rceil$, $c_{\mathrm{iter}}\geq1/2$, and $\kappa\geq1$, this gives $q^{2T}\leq n^{-1}$. Moreover, the condition $n\geq2T$ ensures $B=\lfloor n/T\rfloor\geq n/(2T)$. Substituting these two inequalities into the preceding bound gives
\[
\EE r_T^2
\leq\frac{r_{\mathrm{init}}^2}{n}
+\frac{8\mfour}{(1-q^2)n}
\left(\eta^2T+\frac{\pi^2dT^2}{m^2r_c^2}\right).
\]
Finally, $\dframe(U_T,\Ustar)\leq d_g(U_T,\Ustar)=r_T$, so this estimate proves~\eqref{eq:sgd-expected}.
\end{proof}

\section{Real-data experiments}\label{app:realdata}

The NYT corpus from the UCI bag-of-words repository \cite{newman2008bow} serves as a natural test bed for evaluating $\name$ on real data. Following the degree-normalization, centering, and rank-selection pipeline of \cite{rohe2023}, and excluding localized singular vectors via a fourth-moment diagnostic ($\ell_4=(\sum_iU_{ij}^4)^{1/4}>0.15$), we retain a rank $k=8$ approximation
\[
L_k = U_k\Sigma_kV_k^\top, \qquad\text{with}\quad U_k\in\RR^{n\times 8}\quad\text{and}\quad V_k\in\RR^{p\times8},
\]
of the $n=300{,}000$-document by $p=102{,}660$-word matrix, where the rows of $U_k$ index documents and the rows of $V_k$ index words. The rotational ambiguity described in the introduction comes from the factorization of $L_k$ into document coefficients and word loadings,
\[
L_k = C_0D_0^\top,\qquad\text{where}\quad C_0=\sqrt n\,U_k\quad\text{and}\quad D_0=\tfrac1{\sqrt n}V_k\Sigma_k,
\]
and for any $R\in\orth(8)$ the pair $(C_0R,\,D_0R)$ reproduces the same fit, $(C_0R)(D_0R)^\top=C_0D_0^\top=L_k$. Factor analysis reads $C_0$ as scores and $D_0$ as loadings, while dictionary learning reads $C_0$ as sparse codes and $D_0$ as a dictionary. Rather than resolve this ambiguity once, we exploit it twice, applying $\name$ toward the dictionary-learning goal on the document side and toward the factor-analysis goal on the word side; the two rotations need not agree.
\begin{figure}[htpb!]
    \centering
    \includegraphics[width=\linewidth]{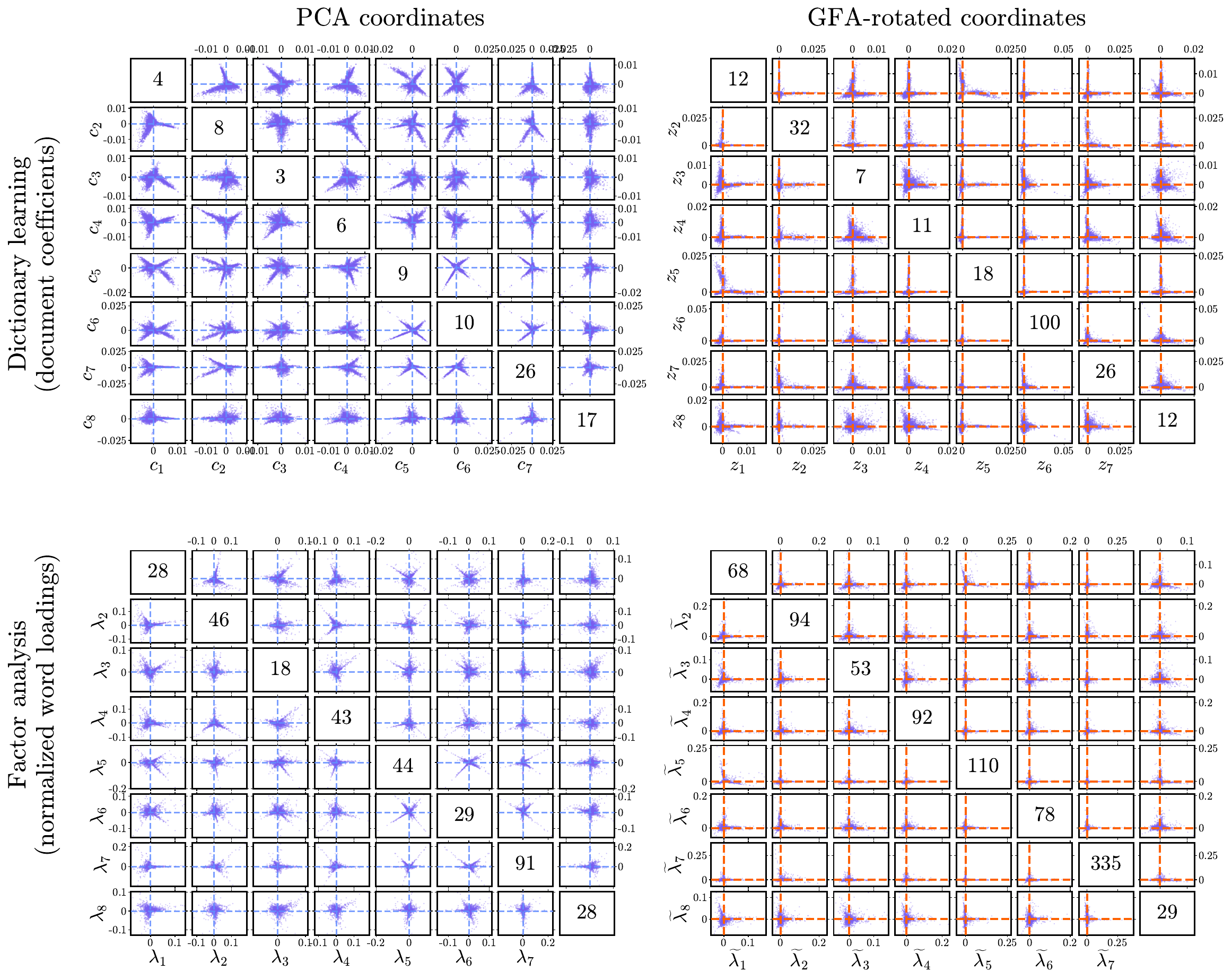}
    \caption{\small \textbf{Two rotations of one low-rank fit on the NYT corpus.} Pairwise scatter plots of a $5{,}000$-row leverage sample from each displayed population, with sample kurtosis of each coordinate on the diagonal. \textit{Top row (dictionary learning):} document coefficients $U_k$ ($c_1,\ldots,c_8$; kurtosis 3--26) and their $\name$-rotated counterparts $U_kR_{\mathrm{code}}$ ($z_1,\ldots,z_8$; kurtosis 7--100, peaking at $z_6$). \textit{Bottom row (factor analysis):} normalized word loadings $V_k$ ($\lambda_1,\ldots,\lambda_8$; kurtosis 18--91) and their $\name$-rotated counterparts $V_kR_{\mathrm{load}}$ ($\tilde\lambda_1,\ldots,\tilde\lambda_8$; kurtosis 29--335, peaking at $\tilde\lambda_7$). $R_{\mathrm{code}}$ and $R_{\mathrm{load}}$ are fit independently.}
    \label{fig:nytimes}
\end{figure}

Concretely, the top row in Figure~\ref{fig:nytimes} fits $R_{\mathrm{code}}\approx\argmin_{R\in\orth(8)}\frac1n\sum_{i=1}^n\ell\bigl((U_k)_{i\cdot},R\bigr)$ and the bottom row fits $R_{\mathrm{load}}\approx\argmin_{R\in\orth(8)}\frac1p\sum_{j=1}^p\ell\bigl((V_k)_{j\cdot},R\bigr)$ (cf.\ \eqref{eq:population_empirical_objectives}), each on the full population of documents or words rather than the displayed sample. We rotate $V_k$, not $D_0=\tfrac1{\sqrt n}V_k\Sigma_k$, because its columns have equal norm, so the bottom row isolates rotational sparsity of the loadings from scale differences already explained by the singular values. Figure~\ref{fig:nytimes} itself plots only a fixed-seed, $5{,}000$-row leverage sample from each side, using the unscaled coordinates $U_k$, $U_kR_{\mathrm{code}}$, $V_k$, and $V_kR_{\mathrm{load}}$.

The radial streaks visible before rotation, in both rows, are the signature of rotational ambiguity that motivates the Wasserstein frame approximation problem~\eqref{eq:generic_approx}. After rotation, mass concentrates near the axes in both rows, though unevenly across coordinates---sharpest at $z_6$ and $\tilde\lambda_7$---since minimizing the average transportation cost is one joint, orthogonally-constrained criterion rather than eight independent ones, and need not distribute concentration evenly. More telling is what the two rows share: the same rank-eight fit, handed to the same rotation-search procedure, concentrates document codes when pointed at $U_k$ and word loadings when pointed at $V_k$. On real unstructured text, one rotational ambiguity thus serves either of the two interpretive goals raised in the introduction, decided not by the algorithm but by which side of the factorization the practitioner cares to align.

\end{document}